\documentclass[11pt]{article}
\usepackage[margin=1in]{geometry}
\usepackage{pgfplots}
\usepackage{microtype}
\usepackage{amsmath,amsfonts,amssymb,amsthm,mathtools,bm} 
\usepackage{graphicx,xcolor}
\usepackage{float,enumitem}
\usepackage{multirow}
\usepackage{bbm}
\usepackage{float}

\usepackage{rotating}
\usepackage{booktabs,tabularx}
\usepackage{subcaption}
\usepackage{algorithm,algorithmic}
\usepackage{tikz}
\usepackage{longtable}
\usepackage{abstract}
\usepackage{hyperref}

\providecommand{\keywords}[1]{\noindent\textbf{Keywords:} #1}

\usepackage{threeparttable}
\usepgfplotslibrary{statistics}
\usetikzlibrary{patterns}
\pgfplotsset{compat=1.17}
\usepackage{siunitx}
\usepackage{placeins}

\newcolumntype{C}{>{\centering\arraybackslash}X}
\newcolumntype{R}{>{\raggedleft\arraybackslash}X}
\newcolumntype{L}{>{\raggedright\arraybackslash}X}

\newcommand{\as}[1]{}
\newcommand{\wk}[1]{}
\newcommand{\dl}[1]{}
\newcommand{\ck}[1]{}

\usepackage[disable]{todonotes}

\usepackage[numbers]{natbib}
\usepackage{appendix}
\usepackage{url}

\newif\ifsupplementary
\supplementaryfalse
\newif\ifinsupplementary

\title{Neural Embedded Mixed-Integer Optimization for Location-Routing Problems}

\author{
  Waquar Kaleem$^{1}$,
  Doyoung Lee$^{2,3}$,
  Changhyun Kwon$^{2,3}$,
  Anirudh Subramanyam$^{1}$ \\
}

\date{
  {\small
  $^{1}$Department of Industrial and Manufacturing Engineering,  
  Pennsylvania State University, University Park, PA 16802, USA \\
  $^{2}$Department of Industrial and Systems Engineering, KAIST, Daejeon 34141, Republic of Korea \\
  $^{3}$Omelet, Inc., Daejeon 34051, Republic of Korea \\
  }
}

\begin{document}

\maketitle
\begin{abstract}
We present a framework that combines machine learning with mixed-integer optimization to solve the Capacitated Location-Routing Problem (CLRP), a classical NP-hard problem that integrates strategic facility location with operational vehicle routing decisions. The proposed method trains a neural network to approximate the cost of a Capacitated Vehicle Routing Problem (CVRP) for serving any subset of customers from a candidate facility. The network is trained on an independently generated dataset of CVRP instances from the literature, entirely separate from any CLRP test instances, thereby avoiding the overfitting and information leakage that can affect learning-based methods. The trained network is then embedded as a surrogate within a mixed-integer model for location-allocation decisions, which is solved using off-the-shelf solvers, thus leveraging decades of advances in vehicle routing and the availability of mature solvers. Computational experiments across four benchmark sets show that the method delivers reasonable solution quality and scales well to large instances, where, after a one-time training cost, it reaches solutions close to the best known at a fraction of the runtime of state-of-the-art heuristics. Our results demonstrate the value of routing cost approximations from the neural surrogate in informing high-quality location-allocation decisions. Our code and data are publicly available.

\end{abstract}

\section{Introduction}

The Location-Routing Problem (LRP) involves determining where to open facilities (or depots) to serve a given set of geographically distributed customers, and identifying the vehicle routes that should be constructed to serve those customers from the opened facilities. In contrast to the classical Facility Location Problem (FLP) \citep{farahani2009facility,laporte2019introduction}, the LRP arises in applications where customers are served by vehicles operating on less-than-truckload routes. Consequently, the cost of any candidate location decision must be evaluated by solving a Vehicle Routing Problem (VRP) \citep{toth2014vehicle} for each opened facility. The LRP thus generalizes both the FLP and the VRP, making it NP-hard in general. 
The Capacitated LRP (CLRP) refers to the simplest variant of the LRP in which both facilities and vehicles have limited capacities, and the objective is to minimize the total cost of facility operations and routing while ensuring that all customer demands are satisfied without exceeding these capacities \citep{schneider2017survey,albaredasambola2020location}.

Depot location decisions directly influence routing decisions due to the integrated nature of the problem. Consequently, solving the CLRP involves making both strategic decisions, such as optimal depot locations and the number of vehicles to operate, as well as tactical decisions, such as vehicle routing plans for each opened depot. The primary objective is to minimize the sum of depot opening costs, fixed vehicle costs, and routing costs across the network.

A natural strategy to solve the CLRP is to make the location and routing decisions independently of each other.
This is partly justified in light of the fact that facility location is a strategic decision, whereas vehicle routing is operational.
In many applications, vehicle routes are often redesigned on a daily basis once the facilities are established.
However, choosing facility locations without taking into account their economic impacts on routing may result in suboptimal solutions, as the possible configurations of feasible vehicle routes are strongly influenced by the locations of open facilities. 
Any initial savings in fixed facility setup costs may not be able to compensate for large losses in distribution in the long run \citep{salhi1989effect,albaredasambola2020location}. Indeed, since distribution is a repetitive activity, any additional routing costs for having chosen a poor facility location will be incurred on a regular basis, and over time, these accumulated costs may exceed savings in setup costs.

The need to address both decisions simultaneously has led to an active area of research in developing new models, solution methods, and applications, as witnessed by the large number of surveys published on the subject \citep{drexl2015survey,schneider2017survey,albaredasambola2020location,mara2021location}.%
Variants of the CLRP can be broadly classified as either discrete or continuous, depending on whether the customer and facility locations are situated on predetermined discrete points or in continuous space.
Within the former category, which is the focus of our work, existing solution methods can be further classified as either exact or heuristic.
Exact approaches include compact mixed-integer programming (MIP) formulations \citep{contardo2013computational}, branch-and-cut techniques \citep{belenguer2011branch,contardo2013computational}, branch-and-price methods \citep{akca2009branch,berger2007location}, and extended MIP formulations based on set partitioning \citep{baldacci2011exact,contardo2014exact}.
Although exact methods can provide mathematical guarantees of suboptimality, they are limited in terms of scalability, being unable to address instances containing more than a few hundred customers or ten candidate facilities. In contrast, heuristic methods can provide solutions to large-scale instances at the expense of losing optimality certificates. 

Recently, learning-based approaches 
have been proposed to estimate routing costs \citep{varol2024neural,sobhanan2024genetic,bogyrbayeva2024machine}.
For any given input of the routing problem of interest, such as the Traveling Salesman Problem or the Capacitated Vehicle Routing Problem (CVRP), these neural surrogate models are trained to produce an estimate of the routing cost.
Once the neural surrogate is trained, it can be used in various applications that require solving the corresponding routing problems as inner optimization problems.
When these neural surrogate models are applied in the context of the CLRP, one can use any heuristic or metaheuristic algorithm to determine the location and depot-customer assignment decisions, and use the neural surrogate to evaluate candidate decisions within the algorithm.
While using neural surrogates within metaheuristic algorithms is a promising direction, this approach has an intrinsic limitation: since the metaheuristic itself governs the location-allocation decisions, its search operators must be redesigned whenever new application-specific constraints arise at this level.

To overcome this limitation, we propose a new approach that formulates an easy-to-solve MIP model in which the vehicle routing cost associated with each opened facility is approximated by a neural surrogate. 
In our approach, any custom constraints can be handled by an off-the-shelf MIP solver, while the neural surrogate enables bypassing the explicit solution of embedded routing problems, thereby making the MIP model solvable within a reasonable time. 
Our method, Neural Embedded Optimization for Location-Routing Problems (NEO-LRP), exploits the fact that the vehicle routing cost associated with a given depot must be \emph{permutation-invariant} with respect to its assigned customers. Specifically, the method uses Deep Sets \citep{zaheer2017deep}
and Graph Transformers \citep{yun2019graph} 
for this purpose, since the architecture naturally respects permutation invariance. We show that these trained models can be readily embedded within a MIP using binary variables and linear constraints, thus enabling their efficient solution using standard MIP solvers.

The NEO-LRP framework decomposes the CLRP into two hierarchical components. At the higher level, it solves a facility location and customer allocation problem to decide which depots to open and how customers are assigned. Unlike classical FLP approaches, where such location-allocation decisions are made independently of routing costs, NEO-LRP employs a neural surrogate model that approximates the routing cost on the induced subgraph formed by each open depot and its allocated customers, thereby providing routing cost information during the location-allocation phase. Once the location-allocation decisions are obtained from this surrogate-informed model, a CVRP is solved on each induced subgraph to construct the actual vehicle routes. In this way, NEO-LRP follows a locate-first-and-route-second decomposition, but differs from traditional approaches by explicitly incorporating routing cost approximation into the location phase through the neural surrogate.

A distinctive feature of our method is its modularity with respect to problem size and allocation-specific side constraints. In particular, a \emph{single pre-trained model} can be used to tackle problems consisting of any number of customers and facility locations, as well as side constraints such as customer incompatibilities and depot-specific assignment rules, all without requiring any fine-tuning on specific benchmarks or significant implementation effort.
Importantly, the neural surrogate is trained to predict the cost of a classical single-depot CVRP %
using an independently generated training dataset from the literature~\citep{uchoa2017new,queiroga202110}
that is completely disjoint from any LRP test instances. 
We caution, however, that this modularity comes at the price of an inability to compute the actual vehicle routes; indeed, NEO-LRP only provides a set of locations where facilities should be opened, along with the allocations of customers to those facilities.
This price is not too steep, however, since the actual vehicle routes can be readily computed \textit{a posteriori} using any available (exact or heuristic) CVRP solver, leveraging decades of advances in vehicle routing algorithms and the availability of mature off-the-shelf solvers.

The main contributions of our work are as follows.
First, we propose NEO-LRP, a novel solution framework that trains graph neural networks to approximate the  cost of a CVRP at each candidate depot. The trained neural networks are then embedded within an easy-to-solve MIP model to determine location-allocation decisions. The neural networks are trained on an independently generated dataset of vehicle routing instances from the literature, entirely separate from any LRP test instances, thereby avoiding overfitting and information leakage.
Second, we show through extensive experiments that NEO-LRP delivers reasonable solution quality and scales well to large instances, where it reaches solutions close to the best known. After a one-time training cost of approximately 6~hours, per-instance solve times range from under a second to a few minutes even for instances with more than 600~customers and 30~depots. Our ablation analyses further demonstrate that the quality of the location-allocation decisions, which are informed by accurate routing cost approximations from the neural surrogate, is central to the overall solution quality and is not diminished by the a~posteriori routing step.

Our paper capitalizes on recent developments in the MIP representation of trained neural networks \citep{tjeng2017evaluating,fischetti2018deep,grimstad2019relu,ceccon2022omlt}. 
To the best of our knowledge, ours is the first approach to employ a neural-embedded MIP framework for solving integrated location and routing problems. A preliminary version of this paper appeared as an extended conference abstract~\citep{kaleem2024neural}. This version did not include any detailed experimental or ablation analyses, relying on a single (and benchmark-specific) method to train a specific neural architecture. In contrast, this paper provides the most comprehensive results to-date, presenting extensive analyses of training sample requirements, neural architectures, target normalization schemes, the choice of routing solver, and the value of routing-aware location-allocation decisions across varying problem sizes.

The remainder of the paper is organized as follows. Section~\ref{sec:literature} provides an overview of the relevant literature;  Section~\ref{sec:formulation} introduces the problem definition and notation;
Section~\ref{sec:neos_clrp} presents the main idea and ingredients of NEO-LRP;
Section~\ref{sec:experiments} presents computational experiments and their findings; and finally, 
Section~\ref{sec:conclusions} concludes the paper with discussions and key takeaways.

\section{Literature Review}\label{sec:literature}

The LRP has been extensively studied for several decades with
numerous comprehensive surveys and literature reviews devoted to the subject \citep{drexl2015survey,prodhon2014survey,cuda2015survey,albaredasambola2020location,schneider2017survey,mara2021location,min1998combined}. %
We refer the reader to the aforementioned surveys for detailed reviews.

\subsection{Exact Methods} 
Exact methods for solving the CLRP aim to find provably optimal solutions but are generally limited to small and medium-sized instances due to their computational complexity. 
These methods typically employ techniques such as branch(-price)-and-cut.
For example, \citet{berger2007location} and \citet{akca2009branch} develop branch-and-price algorithms to solve variants of the LRP with additional constraints on the maximum length of each route and with capacity restrictions, respectively.
Similarly, \citet{belenguer2011branch} proposes a branch-and-cut method that uses binary variables to determine which facilities to open and which arcs to traverse using a vehicle flow formulation.
A computational comparison of various vehicle flow formulations can be found in~\citet{contardo2013computational}, who find that three-index formulations can offer computational advantages over their two-index counterparts.
In contrast to branch-and-cut methods, the work of \citet{baldacci2011exact} develops a set partitioning model and an exact solution strategy that
incorporates lower bounds from solving a relaxed CLRP and the linear programming relaxation of a Multi-Depot Capacitated Vehicle Routing Problem (MDCVRP).
Finally, a solution strategy that combines techniques from both two-index vehicle flow and set partitioning ideas is presented in~\citet{contardo2014exact}. More recently, \citet{liguori2023nonrobust} introduced strong knapsack cuts for the CLRP and related problems, extending the range of instances that can be solved exactly. Despite these advances, exact methods remain computationally demanding and are generally applicable only to small and medium-sized instances, making them unsuitable for large-scale problems.

\subsection{Heuristic Methods} 

Common heuristic methods for the CLRP include simulated annealing (SA), local and tabu search algorithms, population-based algorithms, and savings and insertion methods. For example, \citet{prins2007solving} develops a method that constructs an initial feasible solution using a greedy heuristic, which is then improved by solving the MDCVRP using a guided tabu search that penalizes infeasible solutions. %
Similarly, \citet{prins2006solving} employs a Greedy Randomized Adaptive Search Procedure (GRASP) combined with savings-based heuristics to generate initial solutions that are further refined using local search techniques such as insertion, swap, and 2-opt moves.

A more recent method \citet{loffler2023conceptually} integrates GRASP with Variable Neighborhood Search (VNS). The GRASP provides a diverse set of high-quality initial solutions, while VNS refines them by systematically exploring various neighborhood structures. This method is a representative state-of-the-art baseline against which we compare the performance of our proposed NEO-LRP method. Other works \citep{chan2005multiple,bouhafs2006combination,sahraeian2009using,ting2013multiple,schneider2019large} explore solving the FLP first and then constructing vehicle routes based on those location-allocation decisions; we choose a representative baseline from this class of methods as well.

Methods based on savings and insertion can be found in~\citet{jarboui2013variable}, \citet{prins2006solving}, and \citet{jabal2011variable}. These aim to reduce the total distance traveled by iteratively inserting or merging customers into a set of initial routes.
In contrast, memetic algorithms \citep{prins2006memetic,duhamel2008memetic,derbel2012genetic,duhamel2010grasp}, which are variants of genetic algorithms, construct a single vehicle tour for each facility that is then iteratively improved and made feasible using crossover and repair operations on appropriately defined `chromosomes'. 
A similar two-phase method can be found in~\citet{escobar2013two}, where a giant tour serving all customers is initially constructed and then split according to vehicle capacities, resulting in customer-to-facility allocations that are then individually optimized by solving a Traveling Salesman Problem (TSP). Other heuristics attempt to incorporate mathematical programming techniques. For example, \citet{alvim2013popmusic} attempt to solve the CLRP using an algorithm called Partial Optimization Metaheuristic Under Special Intensification Conditions, which involves solving the TSP on clusters of customers, opening a facility for each cluster, and greedily merging them to improve the overall solution. In another example, \citet{ozyurt2007solving} develop a mixed-integer programming formulation for the nested CLRP and decompose it using Lagrangian relaxation methods. More recent heuristics include the progressive filtering framework of \citet{arnold2021progressive}, 
which reduces the search space by bounding the number of open depots and filtering unpromising configurations through iterative VRP evaluations, 
and the hybrid genetic algorithm of \citet{lopes2016simple}, which combines a standard GA framework with local search in the mutation phase. 
Both methods have been shown to produce competitive results on benchmark CLRP instances.

\subsection{Machine Learning Methods}
The use of machine learning for discrete optimization and algorithmic decision-making has been gaining significant attention in recent years. Our work specifically falls within the broader category of machine learning methods for combinatorial optimization \citep{cappart2023combinatorial,bengio2021machine,NIPS2017_d9896106}.
A recent review of machine learning methods to solve routing problems can be found in \citet{bogyrbayeva2024machine}.
Of particular relevance to our work, however, are those approaches that employ some form of machine learning to build surrogate models or tackle nested optimization tasks \citep{larsen2024fast,NEURIPS2022_9793671e}.
For example, \citet{varol2024neural} suggest the use of a neural network for approximating TSP tour lengths; the trained network is then used to estimate routing costs within a genetic algorithm (GA).
Along similar lines, \citet{sobhanan2024genetic} solve the CLRP by training a graph neural network (GNN) to approximate routing costs and integrating its predictions within a tailored GA.

Both of the aforementioned approaches leverage neural network predictions within heuristics, specifically genetic algorithms, to guide the solution process. Although integrating trained neural networks within heuristics like GA can be beneficial, the algorithms may lack modularity, requiring significant routing-specific knowledge, especially in incorporating additional constraints, while also relying heavily on parameter tuning. 
In contrast, we propose a simple MIP-based framework that can be readily solved by off-the-shelf MIP solvers while also accommodating side constraints. Moreover, the method uses a simple feed-forward permutation-invariant neural network that does not require large amounts of training data. Our method thus strikes a good balance between modeling flexibility, solution quality, and computational efficiency.

\section{Mathematical Formulation of CLRP}
\label{sec:formulation}

Formally, the CLRP can be defined on a directed graph \(\mathcal{G} = (\mathcal{V}, \mathcal{A})\), where the node set \(\mathcal{V} = \mathcal{V}_D \cup \mathcal{V}_C\) consists of potential depot nodes \(\mathcal{V}_D\) and customer nodes \(\mathcal{V}_C\), with $N_D = |\mathcal{V}_D|$ and $N_C = |\mathcal{V}_C|$. The arcs in the set \(\mathcal{A} \subseteq \{(i,j) \in \mathcal{V} \times \mathcal{V} \mid i \neq j,  (i, j) \notin \mathcal{V}_D \times \mathcal{V}_D\}\) represent possible direct travel paths between nodes. Let \(C\) denote the (homogeneous) vehicle capacity, \(Q_d\) the capacity of depot \(d \in \mathcal{V}_D\), and \(F\) the per-vehicle fixed cost. Each customer \(i \in \mathcal{V}_C\) has an associated demand \(q_i\), and opening depot \(d \in \mathcal{V}_D\) incurs a fixed opening cost \(O_d\). Additionally, \(c_{ij}\) represents the travel cost along arc \((i, j) \in \mathcal{A}\). For each node \(i \in \mathcal{V}\), the sets of arcs leaving and entering node \(i\) are denoted by \(\delta^+(i) = \{(i,j)\mid(i,j)\in \mathcal{A}\}\) and \(\delta^-(i) = \{(j,i)\mid(j,i)\in \mathcal{A}\}\), respectively.

The binary variables \(y_d \in \{0, 1\}\) indicate whether depot $d \in \mathcal{V}_D$ is opened; $w_{id} \in \{0, 1\}$ represent allocation decisions indicating if customer $i \in \mathcal{V}_C$ will be served by depot $d \in \mathcal{V}_D$; and \(x_{ijd} \in \{0, 1\}\) represent routing decisions, indicating whether a vehicle originating from depot \(d \in \mathcal{V}_D\) traverses arc \((i, j) \in \mathcal{A}\) (naturally, $x_{ijd} = 0$ if $(i,j)$ involves a depot $d' \neq d$). The continuous variables \(v_{ijd} \geq 0\) represent the vehicle load on arc \((i,j)\in \mathcal{A}\) for vehicles originating from depot \(d \in \mathcal{V}_D\) and are used to model vehicle capacity constraints. We use the arc-flow and commodity-flow variables, $x$ and $v$, purely for ease of presentation; other approaches, such as set partitioning formulations or subtour elimination constraints for modeling capacity limits, can also be used.

A three-index commodity flow formulation for the CLRP, which shall be useful for illustrating our ideas, is as follows. 
\begin{align}
\min\quad
& \mathmakebox[0pt][l]{%
    \sum_{d \in \mathcal{V}_D} O_{d} y_{d}
    + \sum_{d \in \mathcal{V}_D}\sum_{(i,j)\in \mathcal{A}} c_{ij} x_{ijd}
    + F\sum_{d \in \mathcal{V}_D}\sum_{(d,j)\in\delta^{+}(d)} x_{djd}
  } 
  \label{clrp_obj}  \\
\text{s.t.}\quad
& w_{id}\le y_d
		&& \forall i\in \mathcal{V}_C, d\in \mathcal{V}_D, \label{clrp_open}\\ 
&\sum_{d\in \mathcal{V}_D} w_{id}=1
		&& \forall i\in \mathcal{V}_C, \label{clrp_visit}\\ 
& \sum_{i\in \mathcal{V}_C} q_i w_{id} \le Q_d y_d
		&& \forall d\in \mathcal{V}_D, \label{clrp_depot_cap}\\ 
& \sum_{(j,i)\in\delta^{-}(i)} x_{jid}
  = \sum_{(i,j)\in\delta^{+}(i)} x_{ijd}
  = w_{id}
		&& \forall i \in \mathcal{V}_{C}, d \in \mathcal{V}_D, \label{clrp_flow_cust} \\ 
& \sum_{d \in \mathcal{V}_D}\sum_{(i,j)\in\delta^{+}(i)} v_{ijd}
  = \sum_{d \in \mathcal{V}_D}\sum_{(j,i)\in\delta^{-}(i)} v_{jid} -\,q_i
		&& \forall i \in \mathcal{V}_{C}, \label{clrp_load_bal} \\ 
& v_{ijd} \le C x_{ijd}
		&& \forall d \in \mathcal{V}_D, (i,j)\in \mathcal{A}, \label{clrp_vehicle_cap} \\ 
& v_{ijd} \ge 0
		&& \forall d \in \mathcal{V}_D, (i,j)\in \mathcal{A}, \label{clrp_dom_f} \\ 
& x_{ijd} \in \{0,1\}
		&& \forall d \in \mathcal{V}_D, (i,j)\in \mathcal{A}, \label{clrp_dom_x} \\
& w_{id} \in \{0,1\}
		&& \forall i \in \mathcal{V}_{C}, d \in \mathcal{V}_D, \label{clrp_dom_w} \\ 
& y_{d} \in \{0,1\} 
		&& \forall d \in \mathcal{V}_D. \label{clrp_dom_y}
\end{align}
The objective function~\eqref{clrp_obj} minimizes the sum of depot opening costs, travel costs, and fixed vehicle costs. Constraints~\eqref{clrp_open} ensure that no customers are assigned to depots unless they are open. Constraints~\eqref{clrp_visit} ensure that each customer is assigned to exactly one depot. Constraints~\eqref{clrp_depot_cap} impose capacity restrictions for depots. Constraints~\eqref{clrp_flow_cust} enforce flow conservation at each customer node and that it lies along exactly one vehicle route. Constraints~\eqref{clrp_load_bal} enforce vehicle load conservation at customer nodes. Constraints~\eqref{clrp_vehicle_cap} impose capacity restrictions along individual vehicle routes. Finally, constraints~\eqref{clrp_dom_f}--\eqref{clrp_dom_y} specify the domain of the decision variables.

\section{Neural Embedded Optimization for Location-Routing Problems}
\label{sec:neos_clrp}

In this section, we present the NEO-LRP framework. We begin by describing a mathematical decomposition that motivates the use of a neural surrogate for routing costs, followed by an overview of the framework and detailed descriptions of its individual components.

\subsection{Problem Decomposition and Neural Surrogate Formulation}
Exact formulations based on arc-flow variables, such as the one presented above, face severe scalability issues, since the number of binary variables  grows quadratically with the number of customers. As a result, even moderate-sized instances with a few hundred customers and multiple depots can become intractable for standard MIP solvers. Although this is partly alleviated by set-partitioning formulations, these also remain limited in terms of their scalability. This can be primarily attributed to the coupled nature of the depot location and vehicle routing decisions.

To address this, we first propose a reformulation that \emph{decomposes} the CLRP into two hierarchical components: a high-level FLP for depot locations and customer allocations, and a low-level CVRP for constructing vehicle routes. To elaborate, note that the original CLRP objective \eqref{clrp_obj}, reproduced below, consists of depot opening, routing, and vehicle usage costs; of these, only the latter two components involve the $O(N_C^2)$ arc-level routing variables \( x_{ijd} \):
\[
\min \sum_{d \in \mathcal{V}_D} O_d y_d 
+ \sum_{d \in \mathcal{V}_D} \sum_{(i,j)\in \mathcal{A}} c_{ij} x_{ijd} 
+ F \sum_{d \in \mathcal{V}_D} \sum_{(d,j)\in\delta^{+}(d)} x_{djd}.
\]
We therefore approximate the routing and vehicle usage cost components using a depot-level cost function, defined as follows:
\begin{equation}\label{eq:cvrp_cost_function}
R_d(w) \coloneqq 
\min_{x, v \in \mathbb{R}^{|\mathcal{A}|}} \left\{ \sum_{(i,j)\in \mathcal{A}} c_{ij} x_{ijd} + F \hspace{-1em} \sum_{(d,j)\in\delta^{+}(d)} \hspace{-1em} x_{djd} :
\text{\eqref{clrp_flow_cust}--\eqref{clrp_dom_x} with} \ \mathcal{V}_D \ \text{replaced by} \ \{d\}
\right\}.
\end{equation}
If $V \coloneqq \big\{(y, w) \in \mathbb{R}^{N_D} \times \mathbb{R}^{N_C N_D}  :  \text{\eqref{clrp_open}--\eqref{clrp_depot_cap}, \eqref{clrp_dom_w}--\eqref{clrp_dom_y}} \big\}$, then for any depot $d \in \mathcal{V}_D$, we can interpret
$R_d: \mathop{\mathrm{Proj}}_w(V) \to \mathbb{R}$ as the depot-level cost function that maps customers allocated to depot~$d$ to the routing cost associated with serving them using vehicles originating from $d$.
This simplifies the objective~\eqref{clrp_obj} to yield the following reformulation of the CLRP:
\begin{align}
\min_{(y,w) \in V} \quad
& g(y,w) \coloneqq \sum_{d \in \mathcal{V}_D} O_d y_d + \sum_{d \in \mathcal{V}_D} R_d(w). \label{eq:clrp_reform}
\end{align}

A key observation underlying our methodology is that, for any fixed depot $d \in \mathcal{V}_D$, the cost $R_d$ depends only on the induced subgraph
\begin{equation}\label{eq:depot_subgraph}
\mathcal G_d \coloneqq \mathcal G[\{d\} \cup S_d(w)],
\end{equation}
where $S_d(w) \coloneqq \{i \in \mathcal{V}_C: w_{id} = 1\}$ is the subset of customers allocated to depot~$d$.
To illustrate this, consider Figure~\ref{fig:tikz_diagram1} which depicts a feasible CLRP solution, with depots represented as squares and customers as circles. If $\hat{w}$ is the allocation vector for the solution shown in Figure~\ref{fig:tikz_diagram1}, then we have $S_1(\hat{w}) = \emptyset$,
$S_2(\hat{w}) = \{\text{a}, \text{b}\}$,
$S_3(\hat{w}) = \{\text{c}, \text{d}, \text{i}\}$,
$S_4(\hat{w}) = \{\text{e}, \text{f}, \text{g}, \text{h}\}$,
and $S_5(\hat{w}) = \emptyset$.

\begin{figure}
    \centering
    \begin{tikzpicture}[>=latex, scale=0.9]

        \begin{scope}
            \node[draw, thick, minimum size=0.4cm] (n1) at (0,0) {$4$};
            \node[draw, thick, minimum size=0.4cm] (n3) at (2,2) {$3$};
            \node[draw, thick, minimum size=0.4cm] (n4) at (-1.5,3) {$2$};
            \node[draw, pattern=north east lines, pattern color=gray, minimum size=0.4cm] (n5) at (-3,0.5) {$1$};
            \node[draw, pattern=north east lines, pattern color=gray, minimum size=0.4cm] (n5) at (1.3 ,4) {$5$};

            \node[circle, draw, minimum size=0.3cm] (a) at (-3.0,3.5) {};
            \node[anchor=east, xshift=-0.5pt, inner sep=0pt] at (a.west) {a};
            \node[circle, draw, minimum size=0.3cm] (b) at (-0.5,4) {}; \node[anchor=west] at (-0.4,4) {b}; 
            \node[circle, draw, minimum size=0.3cm] (c) at (0,2) {}; \node[anchor=west] at (0.1,2.2) {c}; 
            \node[circle, draw, minimum size=0.3cm] (d) at (3,3) {}; \node[anchor=west] at (3.1,3) {d}; 
            \node[circle, draw, minimum size=0.3cm] (e) at (2,-1) {}; \node[anchor=west] at (2.1,-1) {e}; 
            \node[circle, draw, minimum size=0.3cm] (f) at (0,-2) {}; \node[anchor=west] at (0.1,-2.1) {f}; 
            \node[circle, draw, minimum size=0.3cm] (g) at (-1,1) {}; \node[anchor=west] at (-0.9,1) {g}; 
            \node[circle, draw, minimum size=0.3cm] (h) at (-1,-1) {}; \node[anchor=west] at (-0.9,-1) {h}; 
            \node[circle, draw, minimum size=0.3cm] (i) at (1,0.8) {}; \node[anchor=west] at (1.1,0.8) {i}; 

            \draw[->, thick] (n4) -- (a); \draw[->, thick] (a) -- (b); \draw[->, thick] (b) -- (n4);
            \draw[->, thick, bend left=30] (n3) to (d);
            \draw[->, thick, bend left=30] (d) to (n3);
            \draw[->, thick] (n3) -- (c); \draw[->, thick] (c) -- (i); \draw[->, thick] (i) -- (n3);
            \draw[->, thick] (n1) -- (g); \draw[->, thick] (g) -- (h); \draw[->, thick] (h) -- (n1);
            \draw[->, thick] (n1) -- (f); \draw[->, thick] (f) -- (e); \draw[->, thick] (e) -- (n1);
        \end{scope}
        
        \begin{scope}[xshift=4cm, yshift=1cm]
            \small
            \matrix[matrix anchor=west, column sep=0.1cm, row sep=0.4cm, nodes={anchor=west}] {
                \node[draw, minimum size=0.3cm] {}; & \node {Opened Depot}; \\
                \node[draw, pattern=north east lines, pattern color=gray, minimum size=0.3cm] {}; & \node {Candidate Depot (not opened)}; \\
                \node[circle, draw, minimum size=0.25cm] {}; & \node {Customer}; \\
                \node[inner sep=0pt] 
                    {\tikz\draw[->, thick] (0,0.15) -- (0.45,0.15);}; & \node {Vehicle Route}; \\
            };
        \end{scope}

    \end{tikzpicture}
    \caption{An illustrative CLRP solution with 5 candidate depots and 9 customers. The solution opens depots 2, 3, and 4, while depots 1 and 5 remain closed. The customers allocated to each depot are as follows: $S_1(\hat{w}) = \emptyset$, $S_2(\hat{w}) = \{\text{a}, \text{b}\}$, $S_3(\hat{w}) = \{\text{c}, \text{d}, \text{i}\}$, $S_4(\hat{w}) = \{\text{e}, \text{f}, \text{g}, \text{h}\}$, and $S_5(\hat{w}) = \emptyset$. The vehicle routes are depicted as directed arrows between nodes.
    }
    \label{fig:tikz_diagram1}
\end{figure}

The above observation might seem trivial, but it is crucial for our subsequent development.
In particular, it implies that $R_d$ is a function of the graph $\mathcal G_d$, and hence (with slight abuse of notation), we can equivalently use $R_d(\mathcal G_d)$ to refer to $R_d(w)$.
Therefore, $R_d(\mathcal G_d)$ can be interpreted as the optimal cost of the CVRP defined on the input graph $\mathcal G_d$ and depot node $d$.
This motivates the use of a neural network to learn a surrogate $\hat{R}_d(\cdot)$ that can replace $R_d(\cdot)$ in the MIP model~\eqref{eq:clrp_reform}.
By replacing the large number of arc-flow (or set partitioning) variables with a small number of scalar terms, the surrogate allows decomposing the CLRP into an FLP and a CVRP. The new challenge, however, is that cost evaluation now depends on a {decision-dependent predictor}. The following sections show how we mitigate this dependence without sacrificing decomposability.

\subsection{Overall Framework}

\begin{figure}[t]
    \centering
    \includegraphics[width=\linewidth]{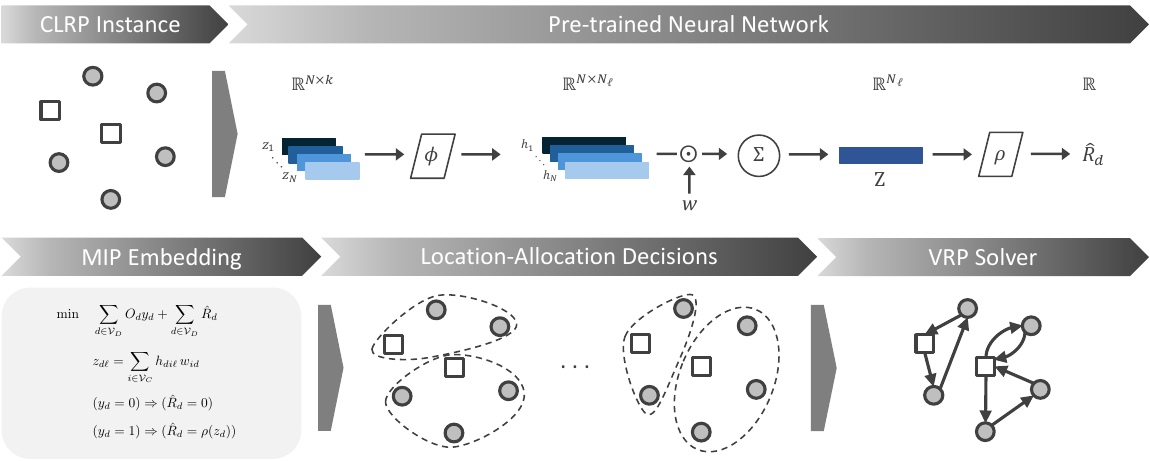}
    \caption{Overview of the NEO-LRP framework. \emph{Top:} A pre-trained neural network processes node-level input features through a feature extractor ($\phi$) to produce embeddings, which are combined with the allocation variables $w$ and aggregated to form a depot-level representation $Z$; a regressor ($\rho$) then predicts the routing cost $\hat{R}_d$. The neural network architecture and notation are detailed in Section~\ref{sec:learning}. \emph{Bottom:} The regressor is embedded within a MIP that determines location-allocation decisions, which are then completed by an off-the-shelf VRP solver to obtain the final routes. %
    }
    \label{fig:Overall-Framework}
\end{figure}

Figure~\ref{fig:Overall-Framework} provides a high-level overview of the NEO-LRP framework. Given a CLRP instance, the method proceeds as follows.
First, a \emph{pre-trained neural network} (Section~\ref{sec:learning}) processes the input data %
to produce a compact numerical representation that approximates the routing cost for each depot.
This network is trained offline, in a one-time preprocessing step, on a dataset of single-depot CVRP instances (Section~\ref{sec:data_generation})
after appropriate normalization (Section~\ref{sec:feature_normalization}).
Notably, these instances are generated independently of any LRP benchmarks using established instance generators from the literature, ensuring that the surrogate generalizes to unseen LRP instances.
The trained neural network is then incorporated into a \emph{MIP formulation} (Section~\ref{sec:NEO-CLRP-formulation}), in which the intractable routing cost function $R_d(\cdot)$ is replaced by its neural surrogate $\hat{R}_d(\cdot)$. Solving this MIP yields {location-allocation decisions} of which depots to open and which customers to assign to each depot. Finally, for each open depot, an off-the-shelf {VRP solver} (Section~\ref{sec:computing_final_routes}) constructs the actual vehicle routes on the induced subgraph.

\subsection{Neural Surrogate Architecture}\label{sec:learning}

Our goal in this section is to build a neural surrogate $\hat{R}_d(\cdot)$ for the depot-level CVRP cost function $R_d(\cdot)$.
\emph{Graph neural networks} are a natural choice of surrogate, since $R_d$ operates on graphs.
For ease of exposition, we temporarily drop the subscript $d$ 
and let $R(\mathcal G')$ denote the optimal cost of the CVRP instance defined on an arbitrary input graph $\mathcal G' = (\mathcal V', \mathcal A')$. Here, the node set $\mathcal V' = S \cup \{d_0\}$ consists of the set $S$ of customers and a single depot node $d_0$.
Each customer $i \in S$ has demand $q_i$ and all customers are served by a homogeneous vehicle fleet of capacity $C$ and per-vehicle fixed cost $F$. 
Each node $u \in \mathcal V'$ is located at coordinates $(s_u,t_u) \in \mathbb R^2$. 

Since $\hat{R}$ operates on graphs, we can model it as a graph neural network of the form \citep{velivckovic2023everything}:
\begin{equation}
	\hat{R}(\mathcal G') = \rho \biggl(\sum_{u \in \mathcal{V}'} h_u \biggr),
\end{equation}
where we can interpret $h_u \in \mathbb{R}^{N_\ell}$ as some \emph{latent representation} of node $u \in \mathcal{V}'$ in latent dimension $N_\ell$, and $\rho:\mathbb{R}^{N_\ell} \to \mathbb{R}$ as a \emph{regressor} that operates on their (component-wise) aggregation.
Let $\mathcal{N}_u := \{ v \in \mathcal V' \mid (u,v) \in \delta^+(u)  \text{ or } (v,u) \in \delta^-(u)\}$ denote the neighbors of $u$ in $\mathcal G'$ and $Z_{\mathcal{N}_u}$ be the (multi)set consisting of neighbor features \{$z_v$ : $v \in {\mathcal{N}_u}$ \}.
Then, the latent vector is computed by a local \emph{feature extractor} 
$h_u = \phi(z_u, Z_{\mathcal{N}_u})$,
where $z_u \in \mathbb{R}^k$ is the raw feature vector associated with each node $u \in \mathcal{V}'$, with $k$ the feature dimension. Here, $z_u$ encodes information such as spatial coordinates $(s_u, t_u)$ and other node-specific attributes like customer demand $q_u$, appropriately scaled and normalized as described in Section~\ref{sec:feature_normalization}.

The feature extractor $\phi$ determines how local neighborhood information is encoded, while the regressor $\rho$ respects any permutation invariance of $\hat{R}$ with respect to nodes in $\mathcal V'$.
In fact, it can be readily seen that $\hat{R}$ will be permutation invariant in $\mathcal V'$ if $\phi$ is permutation invariant in $Z_{\mathcal{N}_u}$.
Defining $\phi$ satisfying these permutation invariance properties is an active area of machine learning research.
In practice, both $\phi$ and $\rho$ are feedforward neural networks. Existing choices for $\phi$ can be grouped into one of three broad classes, depending on how information from the neighborhood $\mathcal{N}_u$ is aggregated into the latent representation $h_u$. These include convolutional, attentional, and message passing models \citep{bronstein2021geometric}.
Graph convolutional and attention networks are relevant to our method.

A \emph{convolutional} choice of $\phi$ applies a fixed kernel to neighbor features:
\begin{equation}
    h_u  =  \phi \biggl(z_u, \sum_{v \in \mathcal{N}_u} c_{uv} \psi(z_v)\biggr), \label{eq:gcn}
\end{equation}
where $\psi$ is a neural network and $c_{uv}$ are predetermined weights (e.g., based on 
travel costs or adjacency).
An \emph{attentional} variant of $\phi$ weight neighbor features using attention scores:
\begin{equation}
    h_u  =  \phi \biggl(z_u, \sum_{v \in \mathcal{N}_u} a(z_u,z_v) \psi(z_v)\biggr), \label{eq:gan}
\end{equation}
where $a(z_u,z_v)$ are learned attention coefficients. This class includes graph transformers and attention-based message passing. 

To better understand the role of local neighborhood information, consider two extreme cases. First, assume that only node features $z_u$ are available and that no adjacency information is exploited: $\mathcal{N}_u=\emptyset$ for every node $u$.
In this case, both~\eqref{eq:gcn} and~\eqref{eq:gan} reduce to $h_u = \phi(z_u)$ so that
\begin{equation} 
    \hat{R}_d(\mathcal G') = \rho \biggl(\sum_{u \in \mathcal{V}'} \phi(z_u) \biggr). \label{eq:ds}
\end{equation}
This is the so-called Deep Sets model~\citep{zaheer2017deep}, where no adjacency structure is exploited. 
At the other extreme, when $\mathcal G'$ is a fully connected complete graph, it can be shown \citep{velivckovic2023everything} that convolutional graph neural networks~\eqref{eq:gcn} again reduce to Deep Sets, losing their ability to exploit adjacency information. Hence, in such cases attention-based Graph Transformers based on~\eqref{eq:gan} offer a more expressive alternative, as the learned attention coefficients $a(z_u, z_v)$ can capture pairwise node interactions even in fully connected graphs. %
In this case, the surrogate reduces to:
\begin{equation}
    \hat{R}_d(\mathcal G') = \rho \biggl(\sum_{u \in \mathcal{V}'} \phi\bigl(z_u, \sum_{v \in \mathcal{V}'} a(z_u,z_v) \psi(z_v)\bigr) \biggr). \label{eq:gt}
\end{equation}
This increased expressiveness, however, comes at the cost of increased training data and difficulty, as the number of pairwise interactions grows quadratically with the number of nodes.

In light of the above discussion, we consider the Deep Sets (DS) surrogate defined in~\eqref{eq:ds} as our primary neural architecture. Indeed, DS provides a simple and interpretable approximation that can handle customer sets $S$ of arbitrary size and allows pre-computation of node latent vectors $h_u$, so that only the regressor $\rho$ is embedded in the optimization model. It also allows for efficient training with limited data, as it does not require learning pairwise interactions between nodes. To explore the potential of exploiting adjacency information, we also explore the Graph Transformers (GT) surrogate defined in~\eqref{eq:gt} as an alternative architecture in our ablation studies. %

During training, each node feature vector $z_u \in \mathbb{R}^k$ is transformed into a latent embedding $h_u = \phi(z_u, Z_{\mathcal N_u}) \in \mathbb{R}^{N_\ell}$, where $\phi : \mathbb{R}^k \times \mathbb{R}^{k|\mathcal N_u|} \to \mathbb{R}^{N_\ell}$ is the feature extractor. These node embeddings are then aggregated into a graph-level representation $Z = \sum_{u \in \mathcal V'} h_u \in \mathbb{R}^{N_\ell}$. Finally, the regressor $\rho : \mathbb{R}^{N_\ell} \to \mathbb{R}$ maps $Z$ to the predicted routing cost $\hat{R}(\mathcal G') \in \mathbb{R}$. See Figure~\ref{fig:Overall-Framework} for an illustration. In the final MIP (Section~\ref{sec:NEO-CLRP-formulation}), only the regressor $\rho$ is embedded, as the feature extractor $\phi$ can be applied offline to pre-compute the node embeddings $h_u$.

\subsection{Data Collection for Surrogate Training}
\label{sec:data_generation}

To train the surrogate $\hat R(\cdot)$, we need to generate a collection of VRP instances 
$\{\mathcal G'^{(m)}\}_{m=1}^M$ using different sampling schemes and computing their routing costs $R(\mathcal G'^{(m)})$.
The goal is to generate a dataset $\mathcal D = \{ (\mathcal G'^{(m)}, R(\mathcal G'^{(m)})) : m = 1,\dots,M \}$ for supervised model training.

A guiding principle in our data generation procedure is to enable the use of a single pre-trained model that can generalize to input graphs $\mathcal G'$ with varying numbers of customers and across different depots. Crucially, the training data consists entirely of single-depot CVRP instances generated using established instance generators from the VRP literature, completely independent of any LRP test instances or the distributions from which they were drawn. This eliminates the risk of overfitting and information leakage from test to training data, a concern that can affect learning-based methods which extract training instances from the very LRP benchmarks used for evaluation.

Specifically, we generate CVRP instances following the methodology of \citet{uchoa2017new} and \citet{queiroga202110}\footnote{Code available at \url{http://vrp.galgos.inf.puc-rio.br/index.php/en/updates}}, with node coordinates scaled to the range $[0,100]$ (compared to $[0,1000]$ in the original code). The number of customers $N$ is varied in the range $\{5,10,\ldots,100\}$. Depot positions are drawn from different spatial configurations (random, centered, cornered), customer positions follow either uniform, clustered, or mixed distributions, demand values are sampled from multiple distributions, and average route sizes are varied from very short to ultra-long. This produces a diverse collection of CVRP instances spanning a wide range of spatial and demand patterns.  

Each VRP instance $\mathcal G'$ in the dataset is solved using the VROOM solver~\citep{vroom_v1.14} with a 5-second time limit and exploration level set to $5$. 
The solver returns the routing cost $R(\mathcal G')$, which includes both the variable travel cost and the fixed vehicle usage cost. 
We generate $100{,}000$ instances for training, $10{,}000$ instances for validation, and a fixed set of $10{,}000$ instances for testing.

\subsection{Feature Normalization, Sparsity Control, and Hyperparameter Optimization}
\label{sec:feature_normalization}
Each node $u \in \mathcal V'$ is associated with raw features of spatial coordinates $(s_u, t_u)$ and demand $q_u$.
To ensure that these features are consistently scaled across different instances, we normalize both the node features and the routing cost.
The normalized feature vector for node $u$ is defined as
$
z_u = \big(P^{-1}(s_u - s_{d_0}),\; P^{-1}(t_u - t_{d_0}),\; \mathbf{1}_{\{u=d_0\}},\; C^{-1}q_u \big),
$
where $P = P(\mathcal G') > 0$ denotes the maximum spread of all customer coordinates (across both  spatial dimensions) in instance
$\mathcal G'$ after centering the depot $d_0$ at $(0,0)$, and is thus fixed for the instance,
and $C$ denotes the vehicle capacity.
For customer nodes $u$, the indicator entry $\mathbf{1}_{\{u=d_0\}}$ equals $0$, while for the depot it equals $1$. 
In the DS architecture, this indicator term is dropped, yielding a 3-dimensional feature vector $z_u \in [-1,1]^3$. The full 4-dimensional feature vector $z_u \in [-1,1]^4$ is used in the GT architecture explored in our ablation studies. %
The same scale factor %
is also employed to normalize the routing cost. 
Specifically, if $R'$ denotes the raw routing cost of instance $\mathcal G'$, then we define the normalized cost as $P^{-1} R'$.
We also explore the effect of alternative normalization schemes in our ablation studies. %

We use the mean squared error (MSE) loss to train the surrogate $\hat{R}$ on the generated training dataset $\mathcal D$. The trained surrogate is then embedded within the MIP formulation described in the next section.
We select hyperparameters based on validation performance and with a goal of controlling the sparsity of the regressor $\rho$ to ensure tractability of the MIP formulation.
More specifically, recall from Section~\ref{sec:neos_clrp} that the surrogate predictor is given by 
$\hat{R}(\mathcal G') = \rho\!\left(\sum_{u \in \mathcal V'} h_u \right)$,
where $h_u = \phi(z_u, Z_{\mathcal N_u})$ are node-level embeddings computed by the feature extractor $\phi$ and $\rho$ is a regressor. 
As we shall see in the next section, embedding $\rho$ into a MIP introduces one new binary variable per neuron.
Therefore, we restrict $\rho$ to a single hidden layer with ReLU activations and explicitly control the latent dimension $N_\ell$. 
The feature extractor $\phi$ may have multiple layers without affecting the final MIP size. Final hyperparameters are provided in
\ifsupplementary
\else
Appendix~\ref{app:hyperparameters}.
\fi

\subsection{Neural Embedded MIP Formulation}
\label{sec:NEO-CLRP-formulation}
The neural-embedded formulation for the location-allocation component of the CLRP integrates latent node embeddings, generated using the DS architecture, into the optimization model. For each depot $d \in \mathcal{V}_D$, the $N_\ell$-dimensional embedding $h_{di} = \phi(z_{di}) \in \mathbb{R}^{N_\ell}$ is pre-computed offline for all customers $i \in \mathcal{V}_C$ using the trained feature extractor $\phi$. The depot's embedding is denoted $h_{d0}$. 

In addition to the location-allocation variables $(y, w)$ defined in~\eqref{eq:clrp_reform}, the formulation defines the following decision variables:  continuous variables \( z_{d\ell} \in \mathbb{R} \) represent aggregated latent embeddings for each depot \( d \in \mathcal{V}_D \) and latent dimension \(\ell\in \mathcal{L}\), and continuous variables \(\hat{R}_d \in \mathbb{R}_{+}\) represent the predicted routing cost for depot \( d \in \mathcal{V}_D \). %
The regressor $\rho$ is a feedforward neural network with a single ReLU hidden layer; %
see Section~\ref{sec:feature_normalization}.

The neural embedded MIP formulation is as follows:
\begin{align}
    \min\quad &
      \sum_{d \in \mathcal{V}_D} O_d y_d
      +\sum_{d \in \mathcal{V}_D} \hat{R}_d \label{neuro_obj}\\
\text{s.t.}\quad &
      \sum_{d\in \mathcal{V}_D} w_{id}=1
      &&\forall i\in \mathcal{V}_C, \label{neuro_demand}\\ 
    & w_{id}\le y_d
      &&\forall d\in \mathcal{V}_D,  i\in \mathcal{V}_C, \label{neuro_assign_open}\\ 
    & \sum_{i\in \mathcal{V}_C} q_i w_{id} \le Q_d y_d
      &&\forall d\in \mathcal{V}_D, \label{neuro_capacity}\\ 
    & z_{d\ell}= 
    \sum_{i\in \mathcal{V}_C} h_{di\ell} \, w_{id}
          &&\forall d\in \mathcal{V}_D, \ell\in \mathcal{L}, \label{neuro_z_def}\\
    & (y_d = 0) \Rightarrow (\hat{R}_d = 0)
          &&\forall d\in \mathcal{V}_D, \label{neuro_indicator_closed}\\
    & (y_d = 1) \Rightarrow \big(\hat{R}_d = \rho(z_d)\big)
          &&\forall d\in \mathcal{V}_D, \label{neuro_indicator_open}\\
    & w_{id}\in\{0,1\}
    &&\forall d\in \mathcal{V}_D, i\in \mathcal{V}_C, \label{neuro_dom_w}\\
    & y_{d}\in\{0,1\}
      &&\forall d\in \mathcal{V}_D, \label{neuro_dom_y}\\
    & z_{d\ell}\in\mathbb{R}
      &&\forall d\in \mathcal{V}_D, \ell\in \mathcal{L}, \label{neuro_dom_z}\\
    & \hat{R}_d\geq 0
      &&\forall d\in \mathcal{V}_D. \label{neuro_dom_c}
\end{align}

The objective function \eqref{neuro_obj} minimizes the sum of depot opening costs and the predicted costs \(\hat{R}_d\). Constraints \eqref{neuro_demand}--\eqref{neuro_assign_open} ensure that each customer is assigned to exactly one open depot. Constraints \eqref{neuro_capacity} impose capacity restrictions for depots. Constraints \eqref{neuro_z_def} aggregate the pre-computed node embeddings \(h\) into depot-specific embedding vectors by summing the embeddings of assigned customers. Constraints \eqref{neuro_indicator_closed}--\eqref{neuro_indicator_open} are indicator constraints that activate the surrogate cost prediction only when the depot is open. The regressor $\rho$ is a feedforward neural network with a single ReLU hidden layer, which introduces additional binary variables when embedded into the MIP~\citep{fischetti2018deep}.
Constraints \eqref{neuro_dom_w}--\eqref{neuro_dom_c} specify the domain of the decision variables.

When using Deep Sets as the underlying neural network architecture, we refer to the resulting formulation \eqref{neuro_obj}--\eqref{neuro_dom_c} as NEO-DS, and when using Graph Transformers, we refer to it as NEO-GT.

\subsubsection{Reduction in Model Complexity} \label{sec:var_reduction}
We compare the number of decision variables in the exact arc-flow formulation from Section~\ref{sec:formulation} and the neural-embedded MIP to assess model complexity. The former includes arc-level routing and flow variables, which scale quadratically with the number of customers. The neural-embedded formulation eliminates these arc-dependent variables and introduces a small set of latent variables \( z_{d\ell} \) and binary variables for the ReLU activations in $\rho$, whose sizes depend only on the number of depots, the embedding dimension \(N_\ell\), and the number of hidden neurons. As noted in Section~\ref{sec:feature_normalization} we can adjust these parameters to directly control model size and explore the trade-off between expressiveness and computational cost.

As shown in Table~\ref{tab:count_vars_v3}, the number of binary variables is reduced by a factor of $45$ for instances with $50$ customers and $5$ depots, increasing to $194$ for instances with $200$ customers and $10$ depots. The continuous variables are similarly reduced, from $364\times$ to over $5{,}700\times$ over the same range of instance sizes. The reduction in model size is a key factor enabling NEO-LRP to scale to larger instances. However, this scalability comes at a cost. Since $\hat{R}_d$ is a learned surrogate, it may over- or under-predict the optimal routing cost, and we cannot guarantee whether the objective~\eqref{neuro_obj} is a lower or an upper bound on the original objective~\eqref{clrp_obj}. Nevertheless, the location-allocation decisions returned by the neural-embedded model always define a feasible CLRP solution. Its total cost, obtained by recomputing the routing costs a posteriori (see Section~\ref{sec:computing_final_routes}), is a valid upper bound on the optimal CLRP cost.

\begin{table}
  \centering
    \caption{Variable counts in the arc-flow and NEO-LRP formulations.}
  \label{tab:count_vars_v3}
  \scriptsize
  \setlength{\tabcolsep}{3pt}
  \resizebox{\ifdim\width>\linewidth \linewidth\else \width\fi}{!}{%
  \begin{tabular}{ll*{4}{rr}}
  \toprule
  \multirow{2}{*}{Variable} & \multirow{2}{*}{Description} &
    \multicolumn{2}{c}{$(N_C, N_D)$} &
    \multicolumn{2}{c}{(50, 5)} &
    \multicolumn{2}{c}{(100, 10)} &
    \multicolumn{2}{c}{(200, 10)}\\
  \cmidrule(lr){3-4}\cmidrule(lr){5-6}\cmidrule(lr){7-8}\cmidrule(lr){9-10}
  & & Arc flow & NEO & Arc flow & NEO & Arc flow & NEO & Arc flow & NEO\\
  \midrule
  $y_d$   & Depot location          & $N_D$ & $N_D$ & 5 & 5 & 10 & 10 & 10 & 10 \\
  $x_{ijd}$ & Arc-level routing     & $O(N_C^2 N_D)$ & -- & 12,750 & -- & 101,000 & -- & 402,000 & -- \\
  $w_{id}$  & Customer allocation   & -- & $N_C N_D$ & -- & 250 & -- & 1,000 & -- & 2,000 \\
  $a_{dh}$  & ReLU activation       & -- & $N_{\rho} N_D$ & -- & 30 & -- & 60 & -- & 60 \\
  \cmidrule(lr){1-10}
  \multicolumn{2}{l}{\textit{Total binary variables}} & $O(N_C^2 N_D)$ & $O(N_C N_D)$ & 12,755 & 285 & 101,010 & 1,070 & 402,010 & 2,070 \\
  \multicolumn{2}{l}{\textit{Reduction factor ($\times$)}} & \multicolumn{2}{r}{$O(N_C)$} & \multicolumn{2}{r}{45} & \multicolumn{2}{r}{94} & \multicolumn{2}{r}{194} \\
  \midrule
  $v_{ijd}$ & Arc load              & $O(N_C^2 N_D)$ & --  & 12,750 & --   & 101,000 & --  & 402,000 & -- \\
  $z_{d\ell}$ & Latent aggregation  & -- & $N_\ell N_D$   & -- & 30 & -- & 60 & -- & 60 \\
  $\hat R_d$   & Surrogate cost     & -- & $N_D$   & -- & 5   & -- & 10  & -- & 10 \\
  \cmidrule(lr){1-10}
  \multicolumn{2}{l}{\textit{Total continuous variables}} & $O(N_C^2 N_D)$ & $O(N_\ell N_D)$ & 12,750 & 35 & 101,000 & 70 & 402,000 & 70 \\
  \multicolumn{2}{l}{\textit{Reduction factor ($\times$)}} & \multicolumn{2}{r}{$O(N_C^2/N_\ell)$} & \multicolumn{2}{r}{364} & \multicolumn{2}{r}{1,443} & \multicolumn{2}{r}{5,743} \\
  \bottomrule
  \end{tabular}%
  }

  \vspace{4pt}
  \parbox{\linewidth}{\scriptsize\textit{Note.}
$N_C = |\mathcal{V}_C|$;
$N_D = |\mathcal{V}_D|$;
$N_\ell$ is the latent dimension;
and $N_\rho$ is the number of ReLU neurons in the regressor network $\rho$.
The number of arc routing and load variables is $N_D N_C(N_C+1)$.
We suppose that $\max\{N_\rho, N_\ell\} \ll N_D \ll N_C$. For the specific instance sizes considered in the table, we set $N_\rho = N_\ell = 6$. Variables common to both formulations (i.e., $w_{id}$) are not presented in the Arc flow column.}
\end{table}

\subsection{Obtaining the Final Routes}
\label{sec:computing_final_routes}

The neural embedded MIP formulation (Section~\ref{sec:NEO-CLRP-formulation}) determines depot openings ($y_d$) and customer assignments ($w_{id}$), with routing costs approximated during optimization by the surrogate $\hat{R}_d(\cdot)$. The surrogate provides an approximation of routing costs, and it enables the solver to explore the space of location-allocation decisions without constructing routes explicitly. Since no actual vehicle routes are generated in the location-allocation stage, the actual routes and routing costs have to be explicitly computed in a post-processing step after solving the neural embedded MIP. 

Let $(\hat{w}, \hat{y})$ denote the location-allocation decisions obtained from solving the neural embedded MIP.
We construct for each open depot $d$ with $\hat{y}_d=1$ the set of assigned customers $S_d(\hat{w}) = \{ i \in \mathcal{V}_C : \hat{w}_{id} = 1 \}$ and form the induced subgraph $\mathcal{G}_d = \mathcal{G}[\{d\} \cup S_d(\hat{w})]$. We then solve a CVRP on each $\mathcal{G}_d$ to obtain feasible vehicle routes (either via an exact or a heuristic solver) and the corresponding routing cost $R_d(S_d(\hat{w}))$. The final CLRP cost is obtained by combining depot opening costs with these post-processed routing costs 
$\sum_{d \in \mathcal{V}_D} O_d \hat{y}_d + \sum_{d \in \mathcal{V}_D} R_d(S_d(\hat{w}))$. 

\section{Experimental Results and Discussion}
\label{sec:experiments}

All experiments were conducted on an Intel(R) Xeon(R) Gold 6248R CPU (3.00GHz) and an NVIDIA Tesla P100-PCIE 12GB GPU. 
We use Gurobi 10.0.3~\citep{Gurobi2021} as the MIP solver, the Gurobi machine learning package~\citep{gurobi_ml} for implementing neural network constraints, and
PyTorch 2.7.1~\citep{paszke2019pytorch} for training. We use DeepHyper~\citep{Egele2025} for hyperparameter search in DS. For generating routing cost labels in our training dataset (Section~\ref{sec:data_generation}) and for computing final vehicle routes (Section~\ref{sec:computing_final_routes}), we use the heuristic VRP solver VROOM~\citep{vroom_v1.14}.  Our codes, including pre-trained models, are available at \url{https://github.com/Subramanyam-Lab/NEO-LRP}.

NEO-LRP incorporates supervised machine learning and incurs a \emph{one-time} computational cost for data generation and neural network training. We generated a dataset of $120{,}000$ VRP instances using the sampling methodology described in Section~\ref{sec:data_generation}. The training data consists entirely of single-depot CVRP instances generated using the established instance generator of \citet{uchoa2017new}, and thus \emph{completely independent} of any LRP benchmark instances or the distributions from which they were drawn. In terms of computational times, data sampling ($T_{\text{sampling}}$) took less than one minute, and solving these instances with VROOM (5-second time limit) to obtain routing cost labels  $T_{\text{solving}}$ took approximately 25 minutes with parallelization across 180 CPU cores.
We allocate a time budget of approximately 6 hours per setting for the DS architecture, %
including hyperparameter search and model training, resulting in a total training time of approximately $T_{\text{train}}^{\text{DS}} = 12$ hours. %

\subsection{Benchmark Instances}
\label{sec:benchmark_instances}
To evaluate the performance of our NEO-LRP approach, we consider the classical instances of \citet{tuzun1999two}, \citet{prins2004nouveaux}, and \citet{barreto_2004}, which comprise a total of 79 small and medium-sized test instances. These three sets are referred to as $\mathbb{T}$, $\mathbb{P}$, and $\mathbb{B}$ respectively. More recently, \citet{schneider2019large} introduced a substantially larger and more diverse benchmark set, denoted by $\mathbb{S}$, which contains 202 instances spanning a broad range of characteristics and sizes.

The $\mathbb{B}$ set consists of 13 instances (7 with proven optimal solutions), featuring 21 to 150 customers and 5 to 14 capacitated depots. The $\mathbb{P}$ set includes 30 instances (20 with known optimum values), with customer counts ranging from 20 to 200 and 5 to 10 capacitated depots. The $\mathbb{T}$ set comprises 36 instances (6 with optimal solutions), characterized by 100-200 customers and 10-20 uncapacitated depots. The more extensive $\mathbb{S}$ dataset contains 202 instances (64 with proven optimal values) that vary in customer size, number of depots, vehicle capacities, and depot costs. These instances span 100 to 600 customers and 5 to 30 capacitated depots.

The benchmark instances differ in how routing costs are computed. For benchmarks $\mathbb{P}$ and $\mathbb{S}$, travel costs are given by the Euclidean distance multiplied by a factor of~100 and rounded to the nearest integer, with a per-vehicle fixed cost of $F = 1000$\footnote{See \url{http://prodhonc.free.fr/Instances/instances_us.htm} for instance details and Section~3.1 of~\citet{schneider2019large} for benchmark $\mathbb{S}$.}. For benchmarks $\mathbb{T}$ and $\mathbb{B}$, travel costs are given by the Euclidean distance, with $\mathbb{T}$ including a per-vehicle fixed cost of $F = 10$ and $\mathbb{B}$ setting $F = 0$. To account for this, we train two variants of the surrogate: NEO-DS$^{\text{S}}$, trained on labels that include the scaled distances and fixed cost $F = 1000$, used for $\mathbb{P}$ and $\mathbb{S}$; and NEO-DS$^{\text{U}}$, trained on labels with unscaled distances and $F = 0$, used for $\mathbb{T}$ and $\mathbb{B}$.

\subsection{Comparison with Baselines}
\label{comparison-sota}

We evaluate the performance of NEO-DS against state-of-the-art baselines across four benchmark sets: $\mathbb{P}$ of~\citet{prins2004nouveaux}, $\mathbb{T}$ of \citet{tuzun1999two}, $\mathbb{B}$ of \citet{barreto_2004}, and $\mathbb{S}$ of \citet{schneider2019large} (Section~\ref{sec:benchmark_instances}).

For benchmarks $\mathbb{P}$, $\mathbb{T}$, and $\mathbb{B}$, the baselines are HCC-500K~\citep{hemmelmayr2012adaptive}, which applies an adaptive large neighborhood search combining destroy-and-repair mechanisms with adaptive scoring strategies; TSBA$_{\text{speed}}$~\citep{schneider2019large}, which adopts a fast construction-based approach using sequential assignment and location strategies;  GRASP/VNS~\citep{loffler2023conceptually}, which uses a matheuristic framework integrating randomized construction with variable neighborhood search; and HALNS~\citep{voigt2022hybrid}, a hybrid adaptive large neighborhood search. For benchmark $\mathbb{S}$, the baselines are TSBA$_{\text{basic}}$~\citep{schneider2019large} and HGAMP~\citep{he2025hybrid}, a genetic algorithm featuring multi-depot edge assembly crossover and a multi-population scheme organized by depot configurations.

We recall that NEO-DS generates the routes \textit{a posteriori} using the heuristic solver VROOM~\citep{vroom_v1.14} as discussed in Section~\ref{sec:computing_final_routes}. This is in contrast to the baseline methods that simultaneously determine both location and routing decisions. To evaluate solution quality, we measure the relative gap between the cost of the NEO-LRP solution $(\hat{y}, \hat{w})$ and the cost of the best-known solution $(y^{*}, w^{*})$ reported in the literature:
\begin{align}
    E_{\text{BKS}}^{\text{gap}}
    = \frac{\big|\, g(\hat{y}, \hat{w}) - g(y^{*}, w^{*})\,\big|}{g(y^{*}, w^{*})},
\end{align}
where $g(\cdot)$ is the total CLRP cost of a feasible location-allocation solution $(y, w) \in V$, defined in~\eqref{eq:clrp_reform}. We emphasize that $g(\hat{y}, \hat{w})$ is computed using the \emph{true} routing costs $R_d(\hat{w})$ obtained by solving the CVRP for each open depot a posteriori (Section~\ref{sec:computing_final_routes}), not the surrogate predictions~$\hat{R}_d$.
Smaller values of $E_{\text{BKS}}^{\text{gap}}$ indicate that the solutions are closer to the best-known benchmark.

The BKS values for Set $\mathbb{P}$ are taken from \citet{loffler2023conceptually} and \citet{he2025hybrid}.
For Sets $\mathbb{B}$ and $\mathbb{T}$, the BKS values are taken from \citet{sobhanan2024genetic}, \citet{he2025hybrid}, and \citet{loffler2023conceptually}.
For Set $\mathbb{S}$, the BKS values are taken from \citet{he2025hybrid}. For NEO-DS, we report both the time required for location-allocation decisions ($T_{\text{LA}}$) and the total computational time, including routing ($T_{\text{total}}$). For baseline methods, we report the total computational time ($T_{\text{total}}$) as reported in their respective papers~\citep{hemmelmayr2012adaptive,schneider2019large,loffler2023conceptually,he2025hybrid}. All results are averaged over five runs. We note that computation times may not be directly comparable across baseline methods due to hardware variations.

\begin{figure}[!ht]
    \centering
    \begin{subfigure}{0.48\linewidth}
        \centering
        \includegraphics[width=\linewidth]{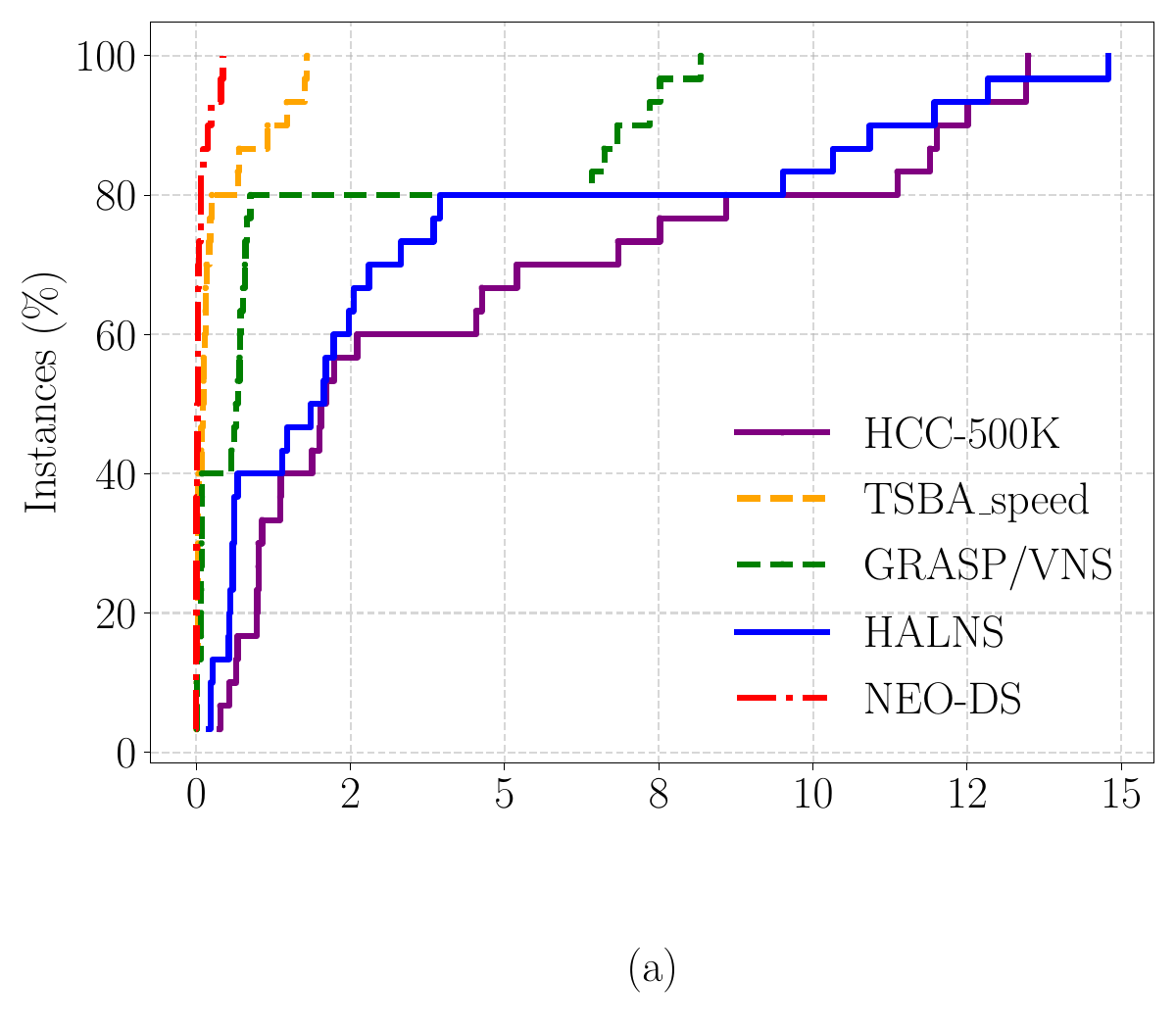}
        \label{fig:ecdf_runtime_prodhon}
    \end{subfigure}
    \hfill
    \begin{subfigure}{0.48\linewidth}
        \centering
        \includegraphics[width=\linewidth]{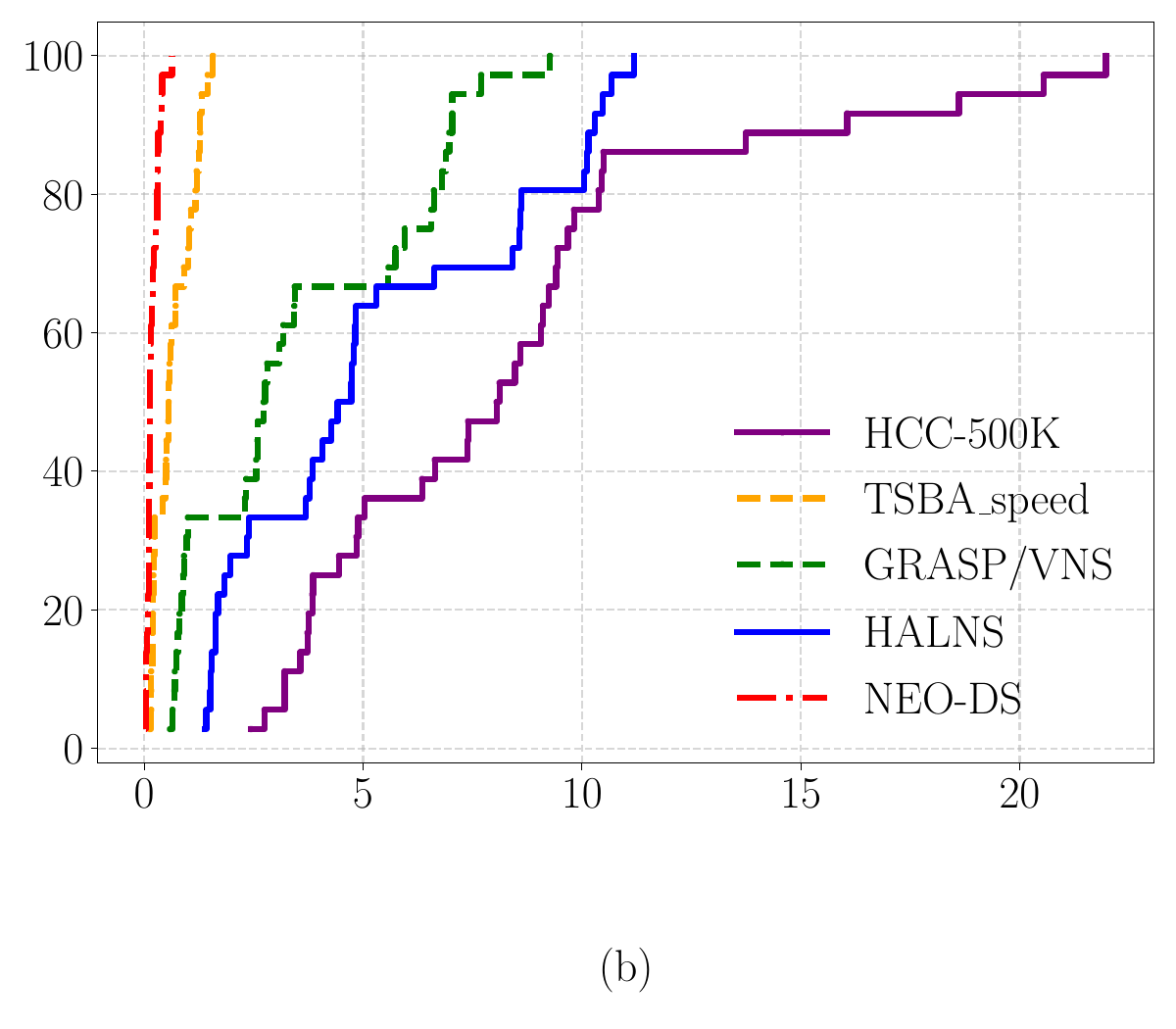}
        \label{fig:ecdf_runtime_tuzun}
    \end{subfigure}
    \begin{subfigure}{0.48\linewidth}
        \centering
        \includegraphics[width=\linewidth]{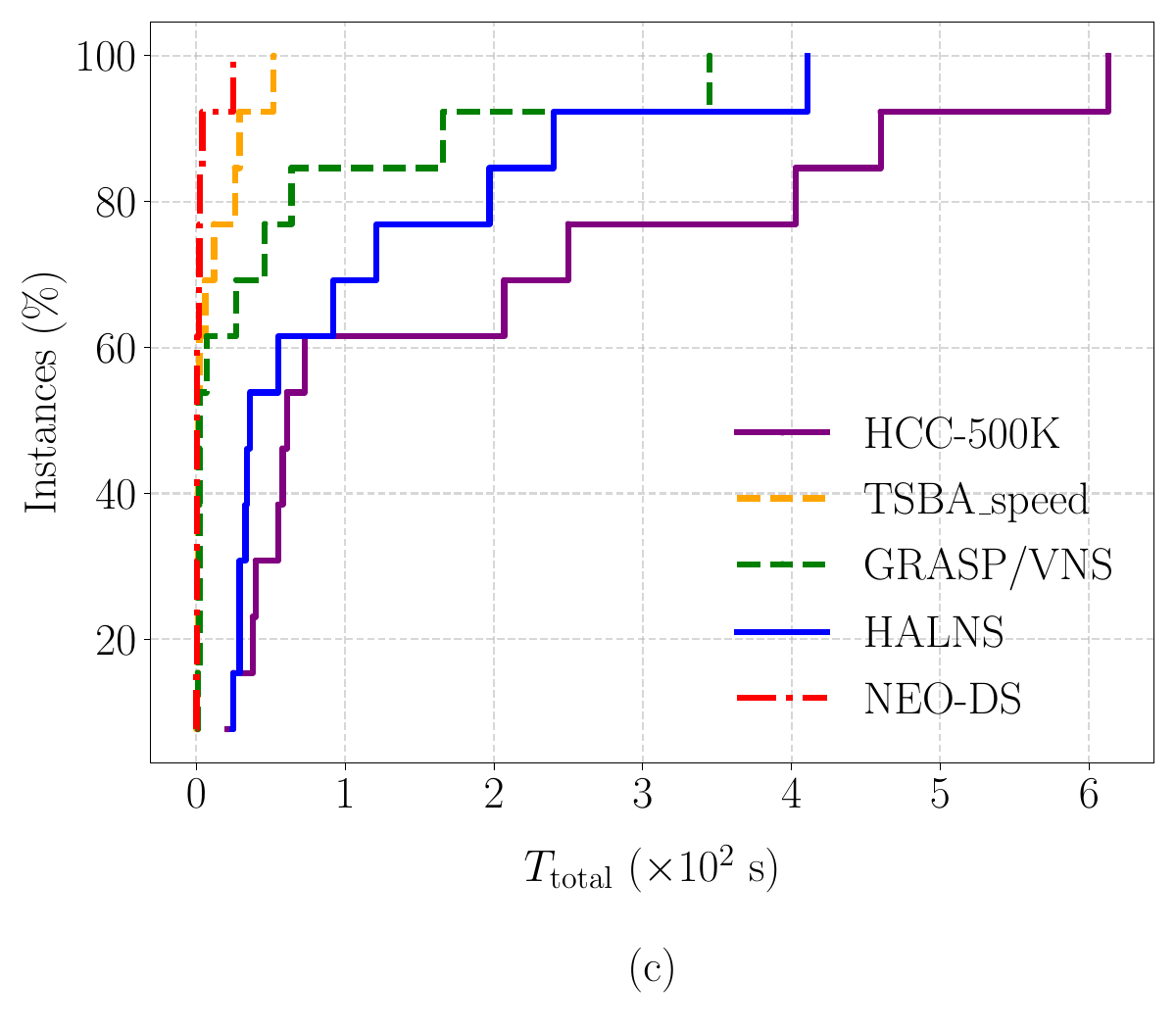}
        \label{fig:ecdf_runtime_barretto}
    \end{subfigure}
    \hfill
    \begin{subfigure}{0.48\linewidth}
        \centering
        \includegraphics[width=\linewidth]{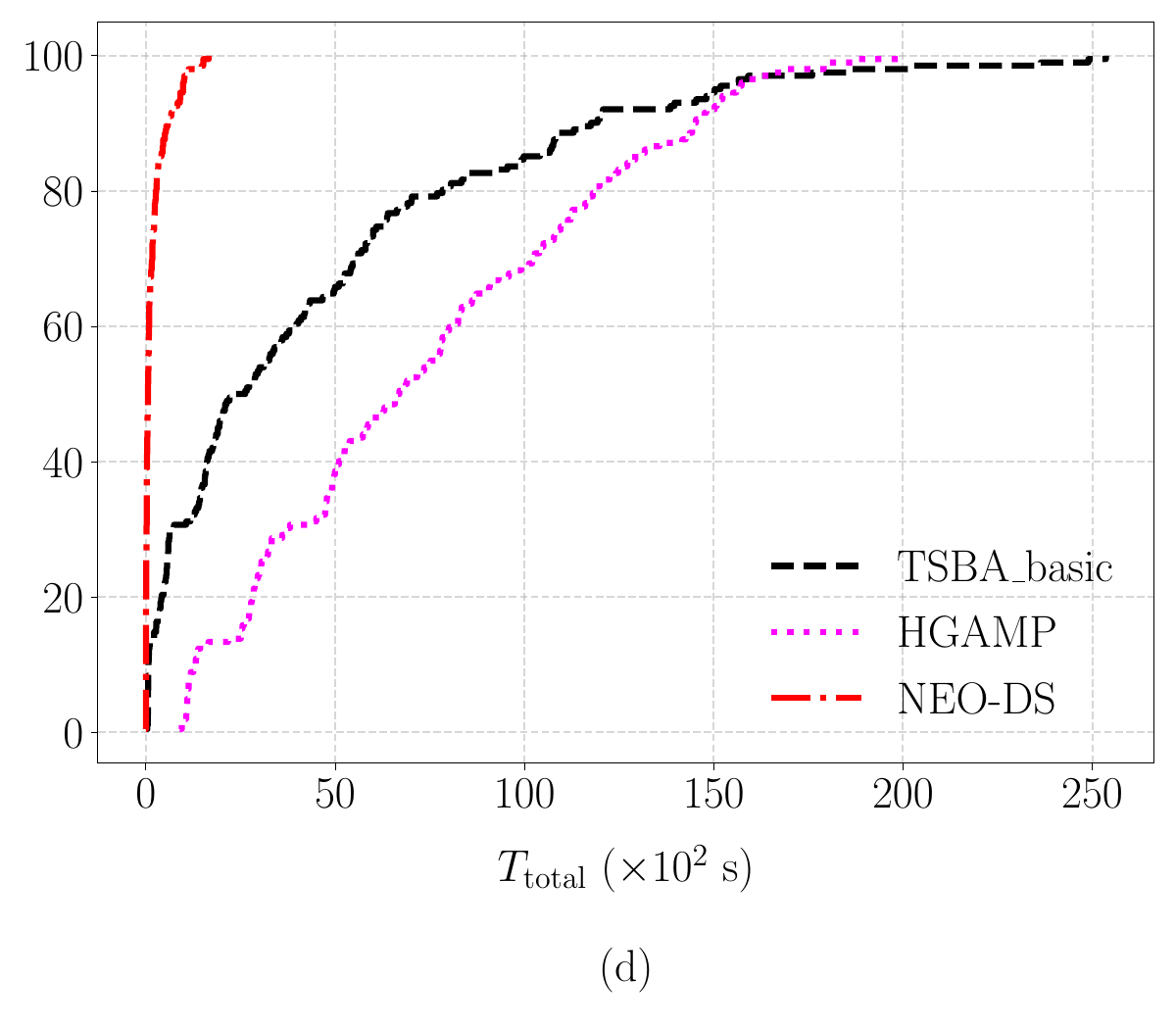}
        \label{fig:ecdf_runtime_schneider}
    \end{subfigure}
    \vspace{-1.5em}
    \caption{Empirical cumulative distribution of computational time ($T_{\text{total}}$): (a)~$\mathbb{P}$ benchmark of \citet{prins2004nouveaux}, (b)~$\mathbb{T}$ benchmark of \citet{tuzun1999two}, (c)~$\mathbb{B}$ benchmark of \citet{barreto_2004}, and (d)~$\mathbb{S}$ benchmark of \citet{schneider2019large}. Higher and to the left is better.}
    \label{fig:ecdf_time}
\end{figure}

\begin{figure}[!ht]
    \centering
    \begin{subfigure}{0.48\linewidth}
        \centering
        \includegraphics[width=\linewidth]{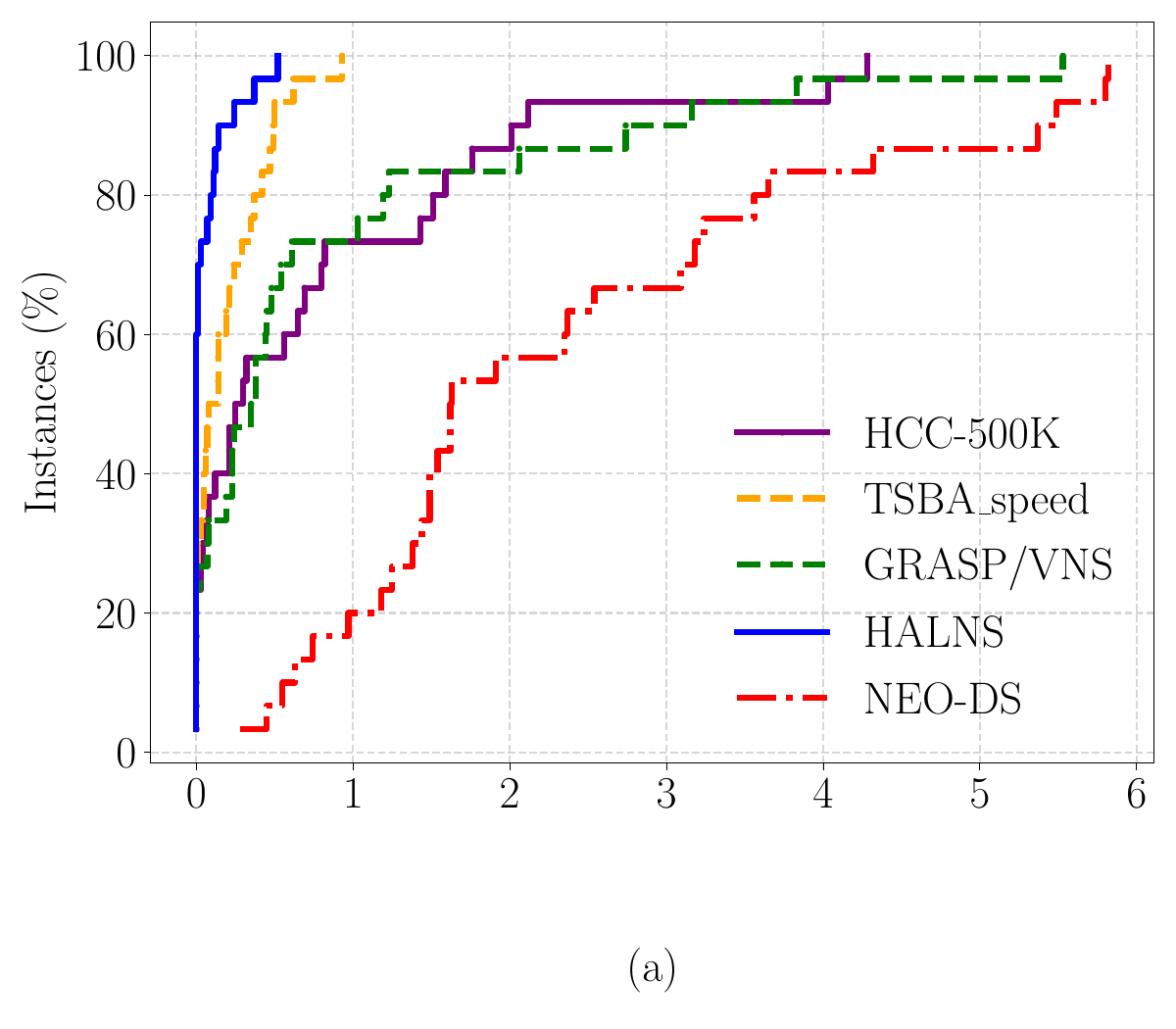}
        \label{fig:ecdf_gap_prodhon}
    \end{subfigure}
    \hfill
    \begin{subfigure}{0.48\linewidth}
        \centering
        \includegraphics[width=\linewidth]{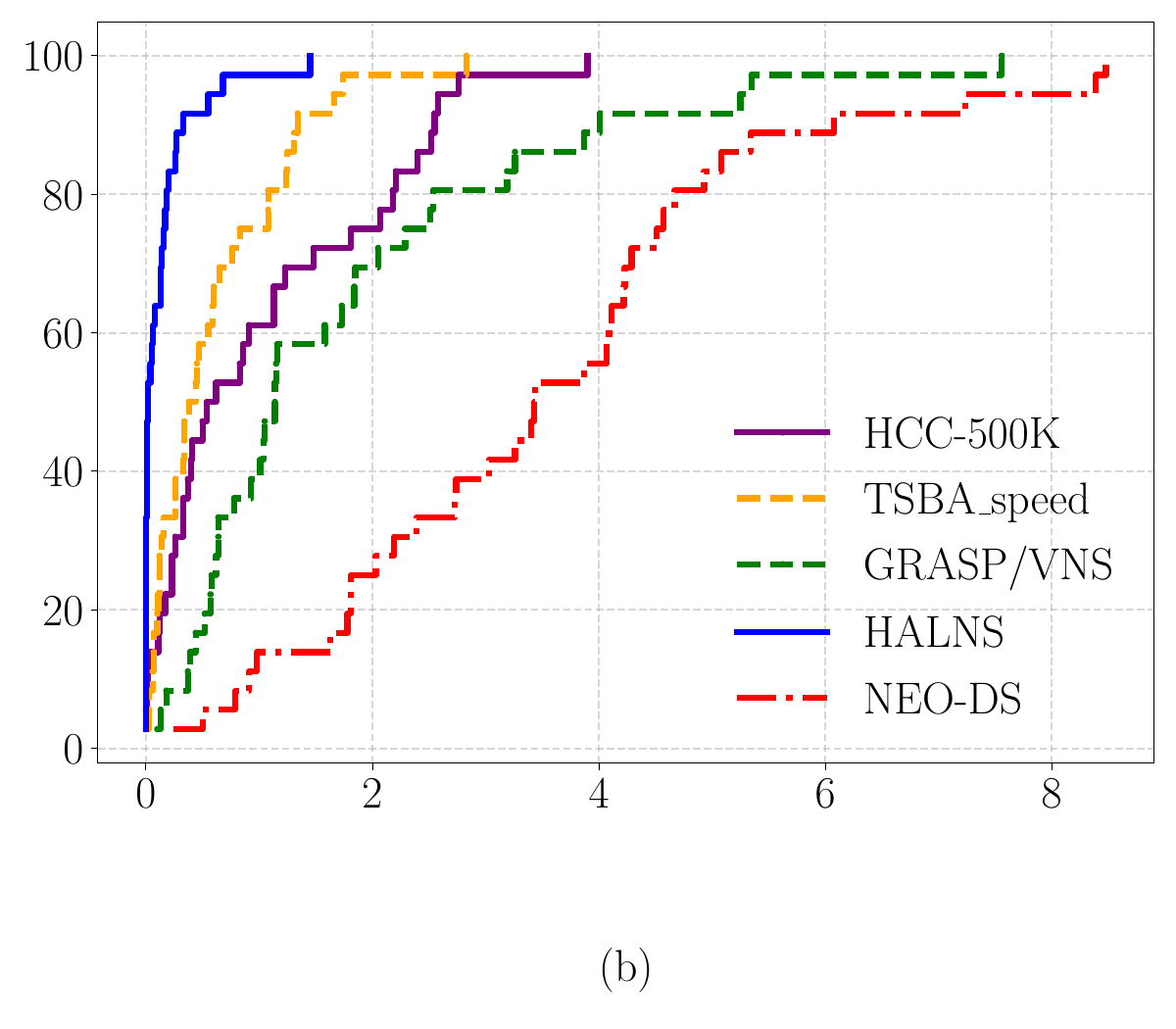}
        \label{fig:ecdf_gap_tuzun}
    \end{subfigure}
    \begin{subfigure}{0.48\linewidth}
        \centering
        \includegraphics[width=\linewidth]{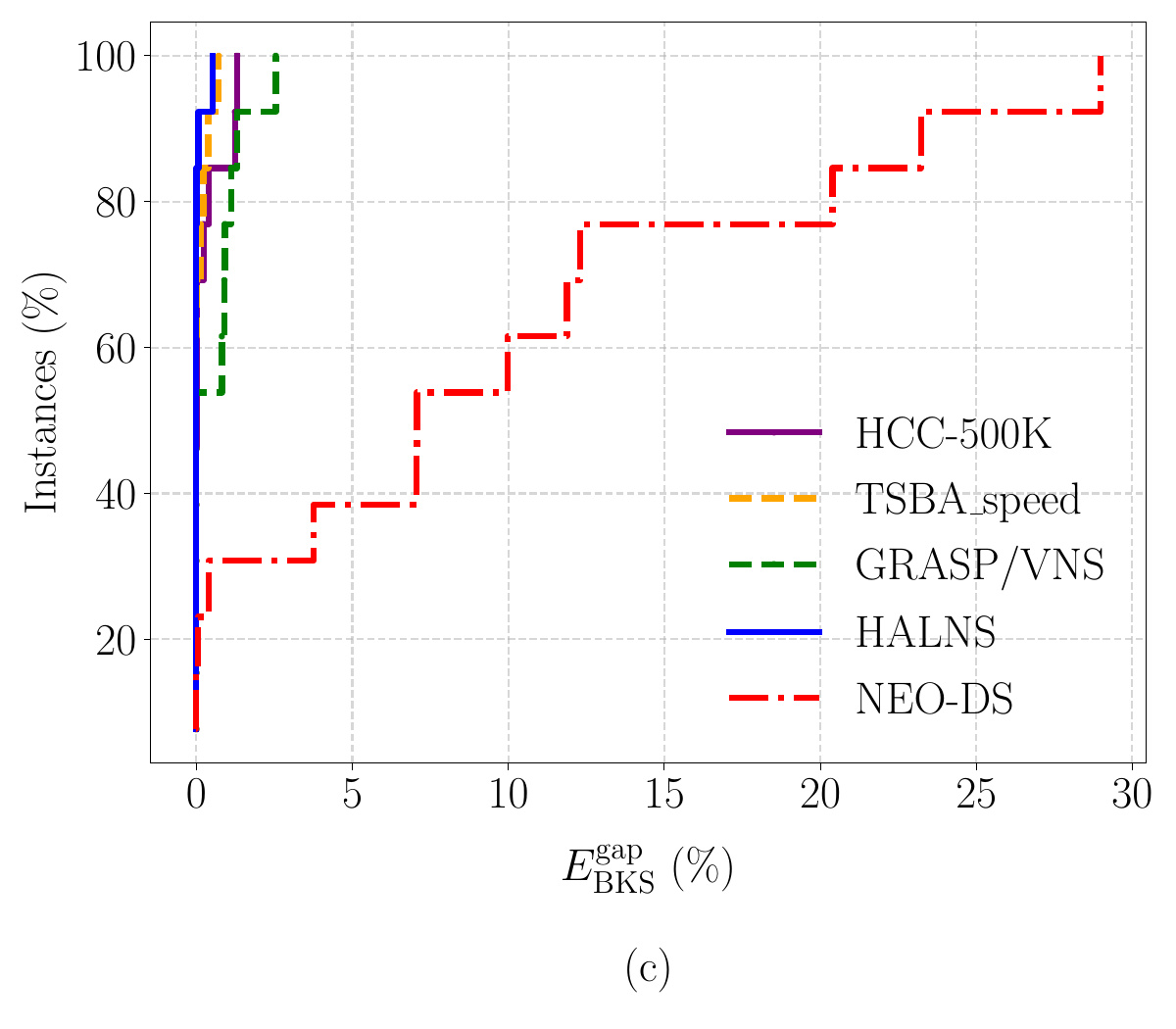}
        \label{fig:ecdf_gap_barretto}
    \end{subfigure}
    \hfill
    \begin{subfigure}{0.48\linewidth}
        \centering
        \includegraphics[width=\linewidth]{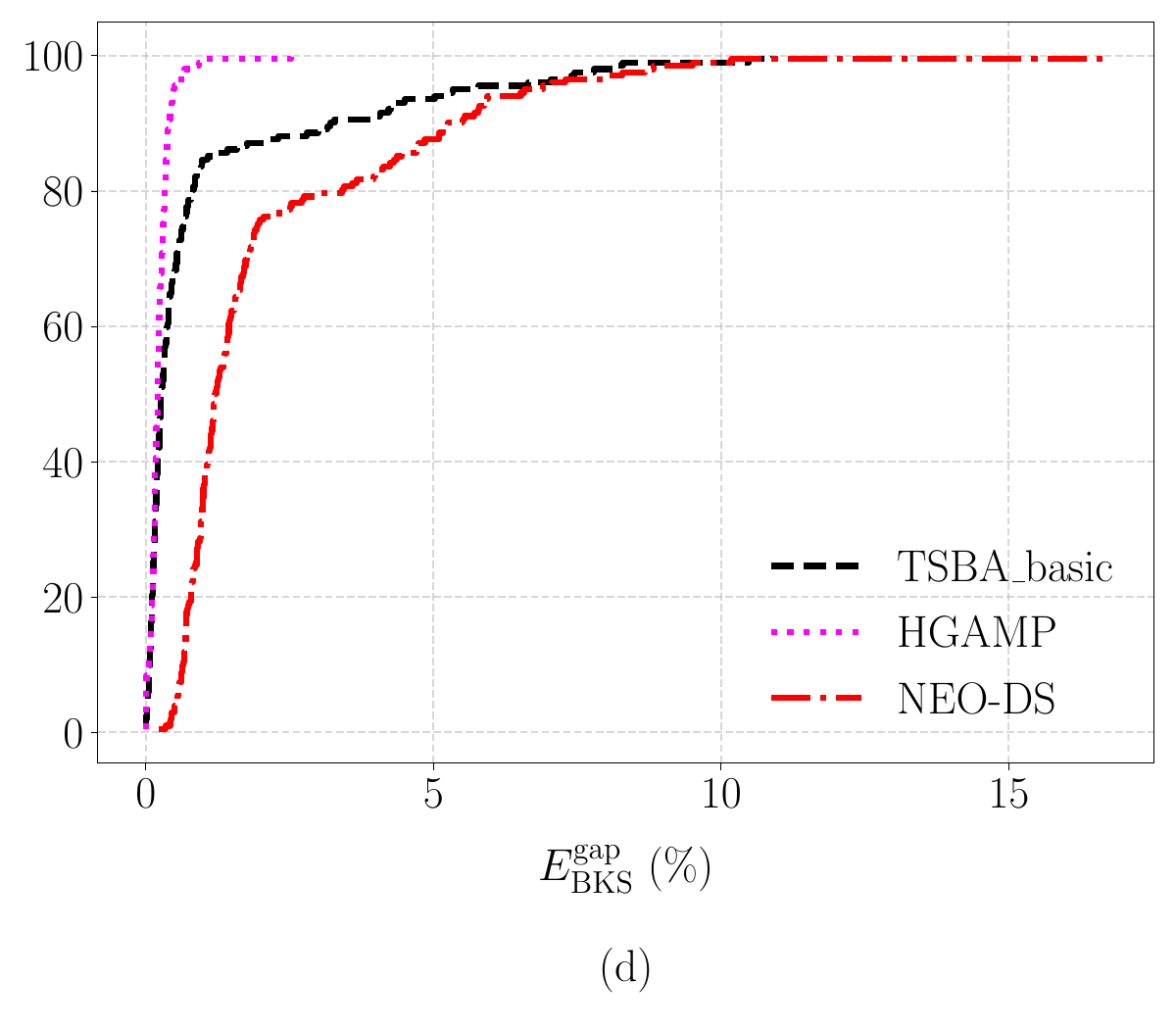}
        \label{fig:ecdf_gap_schneider}
    \end{subfigure}
    \vspace{-1.5em}
    \caption{Empirical cumulative distribution of the gap to the best-known solution ($E^{\text{gap}}_{\text{BKS}}$): (a)~$\mathbb{P}$ benchmark of \citet{prins2004nouveaux}, (b)~$\mathbb{T}$ benchmark of \citet{tuzun1999two}, (c)~$\mathbb{B}$ benchmark of \citet{barreto_2004}, and (d)~$\mathbb{S}$ benchmark of \citet{schneider2019large}. Higher and to the left is better.}
    \label{fig:ecdf_gap}
\end{figure}

Figures~\ref{fig:ecdf_time} and~\ref{fig:ecdf_gap} present the empirical cumulative distribution functions (ECDFs) of computational time and solution quality, respectively, across all four benchmark sets. Furthermore, detailed instance-level results are provided in
\ifsupplementary
\else
Tables~\ref{tab:detailed_prodhon}--\ref{tab:detailed_schneider_500_600} in Appendix~\ref{app:detailed_results}.
\fi

\subsubsection{Prins Benchmark}%

Figures~\ref{fig:ecdf_gap}a and~\ref{fig:ecdf_time}a present the ECDFs for the $\mathbb{P}$ benchmark set. In terms of $E^{\text{gap}}_{\text{BKS}}$ (Figure~\ref{fig:ecdf_gap}a), TSBA$_{\text{speed}}$ achieves a median gap of 0.11\% and all instances within 1\% of the BKS. HALNS achieves the lowest median gap at 0.00\%, solving all instances within 1\% of the BKS. NEO-DS achieves a median gap of 1.62\%, achieving gaps within 2\% for 60.0\% of instances and within 5\% for 90.0\% of instances. In terms of computational times (Figure~\ref{fig:ecdf_time}a), NEO-DS has a median solution time of 2.2s, solving all instances within 100s. TSBA$_{\text{speed}}$ requires a median time of 11.3s, while HCC-500K and GRASP/VNS require 206.0s and 65.9s, respectively, with HCC-500K taking up to 1349.0s for the larger instance with 200 customers and 10 depots. HALNS requires a median time of 196.7s.

\subsubsection{Tuzun-Burke Benchmark}%

Figures~\ref{fig:ecdf_gap}b and~\ref{fig:ecdf_time}b present the ECDFs for the $\mathbb{T}$ benchmark set. In terms of $E^{\text{gap}}_{\text{BKS}}$ (Figure~\ref{fig:ecdf_gap}b), HALNS achieves the lowest median gap at 0.02\%, with 97.2\% of instances within 1\% of the BKS. TSBA$_{\text{speed}}$ achieves a median gap of 0.41\%, with 75.0\% of instances within 1\% of the BKS. HCC-500K achieves a median gap of 0.58\%, with 61.1\% of instances within 1\% of the BKS and all instances within 5\%. GRASP/VNS also demonstrates strong performance with a median gap of 1.14\%, achieving gaps within 2\% for 69.4\% of instances. NEO-DS shows a median gap of 3.33\%, with 86.1\% of instances within 5\% of the BKS.  In terms of computational times (Figure~\ref{fig:ecdf_time}b), NEO-DS achieves the fastest performance with a median time of 10.3s, solving all instances within 100s. TSBA$_{\text{speed}}$ takes 56.0s median time, while GRASP/VNS and HCC-500K require longer times of 274.5s and 808.5s, respectively, with HCC-500K taking up to 2197.0s for instance \texttt{132212}. HALNS requires a median time of 457.5s.

\subsubsection{Barreto Benchmark}%

Figures~\ref{fig:ecdf_gap}c and~\ref{fig:ecdf_time}c present the ECDFs for the $\mathbb{B}$ benchmark set. In terms of $E^{\text{gap}}_{\text{BKS}}$ (Figure~\ref{fig:ecdf_gap}c), TSBA$_{\text{speed}}$ and GRASP/VNS achieve strong performance with median gaps of 0.00\%, with TSBA$_{\text{speed}}$ achieving all instances within 1\% of the BKS and GRASP/VNS achieving 76.9\% within 1\%. HCC-500K  also demonstrates strong performance with a median gap of 0.01\% and 84.6\% of instances within 1\%. HALNS also achieves a median gap of 0.00\%, solving all instances within 1\% of the BKS. In contrast, NEO-DS shows a weaker solution quality with a median gap of 7.06\%, though achieving 38.5\% of instances within 1\% of the BKS. In terms of computational times (Figure~\ref{fig:ecdf_time}c), NEO-DS demonstrates fast computation with a median time of 0.6s, solving all instances within 40s. TSBA$_{\text{speed}}$ and GRASP/VNS require median times of 1.0s and 2.0s, respectively, while HCC-500K requires a median time of 61.0s and up to 613.0s for larger instances with 150 customers and 10 depots. HALNS requires a median time of 35.9s. 

\subsubsection{Schneider Benchmark}%

Figures~\ref{fig:ecdf_gap}d and~\ref{fig:ecdf_time}d present the ECDFs for the $\mathbb{S}$ benchmark set, which contains the largest and most diverse instances. In terms of $E^{\text{gap}}_{\text{BKS}}$ (Figure~\ref{fig:ecdf_gap}d), HGAMP achieves the best performance with a median gap of 0.21\% and 99.5\% of instances within 1\% of the BKS. TSBA$_{\text{basic}}$ demonstrates strong performance with a median gap of 0.27\%, achieving 84.7\% of instances within 1\%. NEO-DS achieves a median gap of 1.21\%, achieving 75.7\% of instances within 2\% and 87.6\% within 5\% of the BKS. In terms of computational times (Figure~\ref{fig:ecdf_time}d), NEO-DS requires a median time of 55.4s, with 63.9\% of instances solved within 100s. In contrast, TSBA$_{\text{basic}}$ and HGAMP require substantially longer times with median values of 2432.0s and 6677.5s, respectively, with HGAMP taking up to 19931.0s for the largest instances with 600 customers and 30 depots.

\subsubsection{Summary of Baseline Comparisons}
Across all four benchmark sets, NEO-DS consistently demonstrates a strong balance between solution quality and computational efficiency. NEO-DS achieves median gaps ranging from 1.21\% ($\mathbb{S}$) to 7.06\% ($\mathbb{B}$), with all instances solved within 100s on benchmarks $\mathbb{P}$, $\mathbb{T}$, and $\mathbb{B}$, and 63.9\% of instances on $\mathbb{S}$. On each benchmark, NEO-DS runs faster than the best-performing baseline, with the reduction in solve time ranging from $2\times$ on $\mathbb{B}$ (1.0s to 0.6s) to $120\times$ on $\mathbb{S}$ (6677.5s to 55.4s), and $5\times$ on $\mathbb{P}$ (11.3s to 2.2s) and $78\times$ on $\mathbb{T}$ (808.5s to 10.3s) in between. The faster solve times come with some loss in solution quality. On the largest and most diverse benchmark $\mathbb{S}$, NEO-DS reaches within $2\%$ of the best known solution on $75.7\%$ of instances. HALNS reaches the best solution quality on the classical sets, but at median runtimes of 35.9s to 457.5s, one to two orders of magnitude slower than NEO-DS.

\subsection{Ablation Studies}
\label{sec:ablation_studies}
We perform detailed ablation studies to understand the contribution of each component of NEO-LRP and to identify which design choices are essential for solution quality. A central question is whether the neural surrogate's routing-aware location-allocation decisions provide genuine value over simpler alternatives, or whether their effect is diminished by the ex-post VRP solution step. Section~\ref{sec:val_loc_alloc} directly addresses this by comparing NEO-DS against a natural baseline that replaces the neural surrogate with a standard facility location model. We further examine the effect of problem size (Section~\ref{sec:problem_size}),
\ifsupplementary
\else
training sample size and generalization (Section~\ref{sec:sample_size}), the choice of routing solver for computing final routes (Section~\ref{sec:routing_solver}), the neural architecture (Section~\ref{sec:ablation_architecture}), and the target normalization scheme (Section~\ref{sec:target_normalization}).
\fi

Unless otherwise stated, the default configuration is NEO-DS trained using 110{,}000 VRP instances sampled using the methodology of Section~\ref{sec:data_generation} with scaled labels, evaluated on the $\mathbb{P}$ benchmark set, %
and the VROOM solver (Section~\ref{sec:computing_final_routes}) for generating final routes. All results are averaged over 5~runs.

\subsubsection{Value of Neural Surrogate}
\label{sec:val_loc_alloc}

We compare NEO-DS against a baseline approach that we call FLP-VRP, which replaces the neural surrogate with a standard capacitated facility location problem (FLP) based on direct travel costs, followed by the same CVRP routing stage used by NEO-DS. Since both methods use the same VRP solver to compute final routes, any difference in total CLRP cost is attributable solely to the location-allocation decisions. We highlight that the direct travel cost in FLP-VRP is a simple, interpretable \emph{distance-based} routing-cost proxy. The FLP stage uses binary variables $w_{id}\in\{0,1\}$ indicating whether customer $i\in\mathcal{V}_C$ is assigned to depot $d\in\mathcal{V}_D$, and is formulated as follows:

\begin{align}
\min_{y,w}\quad & 
   \sum_{d\in\mathcal{V}_D} O_d y_d 
   + \sum_{d\in\mathcal{V}_D}\sum_{i\in\mathcal{V}_C} c_{id} w_{id} \label{flpvrp_obj} \\
\text{s.t.}\quad
& \sum_{d\in\mathcal{V}_D} w_{id} = 1 
   && \forall i\in\mathcal{V}_C, \label{flpvrp_assign}\\
& w_{id} \le y_d 
   && \forall d\in\mathcal{V}_D, i\in\mathcal{V}_C \label{flpvrp_open}\\
& \sum_{i\in\mathcal{V}_C} q_i w_{id} \le Q_d y_d 
   && \forall d\in\mathcal{V}_D, \label{flpvrp_capacity}\\
& w_{id}\in\{0,1\}, y_d\in\{0,1\}
   && \forall d\in\mathcal{V}_D, i\in\mathcal{V}_C. \label{flpvrp_dom}
\end{align}
The objective~\eqref{flpvrp_obj} minimizes depot opening and direct travel costs. Constraints~\eqref{flpvrp_assign} ensure each customer is assigned to exactly one depot, \eqref{flpvrp_open} enforce that customers can only be assigned to open depots, and \eqref{flpvrp_capacity} impose depot capacity limits.

Let $S_d = \{ i\in\mathcal{V}_C : w_{id}=1 \}$ denote the set of customers assigned to depot $d$. In the routing stage, each open depot $d$ with assigned set $S_d$ induces a CVRP on the subgraph $\mathcal{V}_d(S_d)=\{d\}\cup S_d$ with arc set $\mathcal{A}_d(S_d)=\mathcal{A}\cap(\mathcal{V}_d(S_d)\times\mathcal{V}_d(S_d))$. Let $R_d(S_d)$ denote the CVRP cost for depot $d$, including travel and fixed vehicle costs, computed using the heuristic solver VROOM \citep{vroom_v1.14}. The total FLP-VRP objective is therefore $\sum_{d\in\mathcal{V}_D} O_d y_d + \sum_{d\in\mathcal{V}_D} R_d(S_d)$.

In Table~\ref{tab:flp_vrp} we report instance-wise location, routing, and allocation differences between FLP-VRP and NEO-DS on the $\mathbb{P}$ benchmark set of \citet{prins2004nouveaux}. Let $\mathcal{D}^{\text{FLP}}$ and $\mathcal{D}^{\text{NEO}}$ denote the sets of opened depots, $n^{\text{open}} = |\mathcal{D}|$ the number of opened depots, and $F(\mathcal{D}) = \sum_{d \in \mathcal{D}} O_d$ the total fixed facility cost associated with a set~$\mathcal{D}$. We further denote by $R^{\text{FLP}}$ and $R^{\text{NEO}}$ the routing costs returned by FLP-VRP and NEO-DS, respectively. For each instance we report the opened-depot counts $n^{\text{open}}_{\text{FLP}}$ and $n^{\text{open}}_{\text{NEO}}$, the fixed facility costs $F_{\text{FLP}} = F(\mathcal{D}^{\text{FLP}})$ and $F_{\text{NEO}} = F(\mathcal{D}^{\text{NEO}})$, and the routing costs $R^{\text{FLP}}$ and $R^{\text{NEO}}$. We also report the differences (always defined as NEO minus FLP), namely $\Delta n^{\text{open}} = n^{\text{open}}_{\text{NEO}} - n^{\text{open}}_{\text{FLP}}$, $\Delta F = F_{\text{NEO}} - F_{\text{FLP}}$, and $\Delta R = R^{\text{NEO}} - R^{\text{FLP}}$, where negative values indicate that NEO-DS opens fewer depots, has lower fixed costs, or incurs lower routing costs. To quantify differences in the location decisions, we report the normalized symmetric difference of the opened-depot sets,
\[
\|\Delta \mathcal{D}\|_{\text{norm}} = \frac{|\mathcal{D}^{\text{NEO}} \triangle \mathcal{D}^{\text{FLP}}|}{|\mathcal{D}^{\text{NEO}}| + |\mathcal{D}^{\text{FLP}}|}.
\]
To capture changes in location-allocations (customer-to-depot assignments), we compute the normalized Hamming distance between the two assignment vectors, denoted $\|\Delta A\|_{\text{norm}}$, but only when both methods open exactly the same depot set ($\mathcal{D}^{\text{NEO}} = \mathcal{D}^{\text{FLP}}$); otherwise, allocation differences are not directly comparable and are therefore omitted. The Identical column in Table~\ref{tab:flp_vrp} indicates whether the two assignment vectors coincide exactly in those comparable cases.
\begin{table}[htbp]
\centering
\caption{Instance-wise comparison of FLP-VRP and NEO-DS on the $\mathbb{P}$ benchmark set of 
\citet{prins2004nouveaux}. 
For each instance: 
$n^{\text{open}}$ is the number of opened depots; 
$F$ is the total fixed facility cost; 
$R$ is the total routing cost. 
Differences are defined as NEO minus FLP: 
$\Delta n^{\text{open}} = n^{\text{open}}_{\text{NEO}} - n^{\text{open}}_{\text{FLP}}$, 
$\Delta F = F_{\text{NEO}} - F_{\text{FLP}}$, 
$\Delta R = R^{\text{NEO}} - R^{\text{FLP}}$.
$\|\Delta \mathcal{D}\|_{\text{norm}}$ measures the normalized symmetric difference between opened-depot sets. 
$\|\Delta A\|_{\text{norm}}$ measures the normalized Hamming distance between customer assignments, computed only when both methods open the same depot set.}
\label{tab:flp_vrp}
\scriptsize
\setlength{\tabcolsep}{2pt}
\resizebox{\ifdim\width>\linewidth \linewidth\else \width\fi}{!}{%
\begin{tabular}{l r r r  r r r  r r r r  r c}
\toprule
\textbf{Instance} & \multicolumn{3}{c}{\textbf{FLP--VRP}} & \multicolumn{3}{c}{\textbf{NEO--DS}} & \multicolumn{4}{c}{\textbf{Diff.\ (NEO$-$FLP)}} & \multicolumn{2}{c}{\textbf{Alloc (same $\mathcal{D}$)}} \\
\cmidrule(lr){2-4}
\cmidrule(lr){5-7}
\cmidrule(lr){8-11}
\cmidrule(lr){12-13}
& $n^{open}$ & $F$ & $R$ & $n^{open}$ & $F$ & $R$ & $\Delta n^{open}$ & $\Delta F$ & $\Delta R$ & $\|\Delta \mathcal{D}\|_{\text{norm}}$ & $\|\Delta A\|_{\text{norm}}$ & Identical \\
\midrule
20-5-1 & 3 & 25549 & 31019 & 3 & 25549 & 31019 & 0 & 0 & 0 & 0.000 & 0.000 & Yes \\
20-5-1b & 2 & 15497 & 26072 & 2 & 15497 & 25034 & 0 & 0 & -1038 & 0.000 & 0.100 & No \\
20-5-2 & 3 & 24196 & 24998 & 3 & 22769 & 28984 & 0 & -1427 & 3986 & 0.333 & — & — \\
20-5-2b & 2 & 21739 & 21280 & 2 & 13911 & 25693 & 0 & -7828 & 4413 & 0.500 & — & — \\
50-5-1 & 3 & 25442 & 66138 & 3 & 25442 & 66138 & 0 & 0 & 0 & 0.000 & 0.000 & Yes \\
50-5-1b & 3 & 25442 & 44174 & 2 & 15385 & 48603 & -1 & -10057 & 4429 & 0.200 & — & — \\
50-5-2 & 3 & 29319 & 62947 & 3 & 29319 & 61389 & 0 & 0 & -1558 & 0.000 & 0.060 & No \\
50-5-2BIS & 3 & 19785 & 65684 & 3 & 19785 & 65631 & 0 & 0 & -53 & 0.000 & 0.060 & No \\
50-5-2b & 3 & 29319 & 43009 & 3 & 29319 & 41896 & 0 & 0 & -1113 & 0.000 & 0.080 & No \\
50-5-2bBIS & 3 & 18763 & 35088 & 3 & 18763 & 35575 & 0 & 0 & 487 & 0.000 & 0.160 & No \\
50-5-3 & 3 & 37954 & 57418 & 2 & 18961 & 68073 & -1 & -18993 & 10655 & 0.600 & — & — \\
50-5-3b & 3 & 37954 & 38762 & 2 & 18961 & 45538 & -1 & -18993 & 6776 & 0.600 & — & — \\
100-10-1 & 3 & 165068 & 129713 & 3 & 154942 & 137102 & 0 & -10126 & 7389 & 0.333 & — & — \\
100-10-1b & 3 & 165068 & 77765 & 3 & 154942 & 81482 & 0 & -10126 & 3717 & 0.333 & — & — \\
100-10-2 & 3 & 149586 & 98656 & 3 & 149586 & 97631 & 0 & 0 & -1025 & 0.000 & 0.060 & No \\
100-10-2b & 3 & 149586 & 59538 & 3 & 145956 & 63216 & 0 & -3630 & 3678 & 0.333 & — & — \\
100-10-3 & 3 & 136123 & 117243 & 3 & 136123 & 119247 & 0 & 0 & 2004 & 0.000 & 0.120 & No \\
100-10-3b & 3 & 136123 & 71312 & 3 & 136123 & 73253 & 0 & 0 & 1941 & 0.000 & 0.110 & No \\
100-5-1 & 3 & 132890 & 145433 & 3 & 132890 & 146014 & 0 & 0 & 581 & 0.000 & 0.050 & No \\
100-5-1b & 3 & 132890 & 83349 & 3 & 132890 & 83349 & 0 & 0 & 0 & 0.000 & 0.000 & Yes \\
100-5-2 & 2 & 102246 & 92250 & 2 & 102246 & 91963 & 0 & 0 & -287 & 0.000 & 0.010 & No \\
100-5-2b & 2 & 102246 & 56091 & 2 & 102246 & 56091 & 0 & 0 & 0 & 0.000 & 0.000 & Yes \\
100-5-3 & 2 & 88287 & 113313 & 2 & 88287 & 114667 & 0 & 0 & 1354 & 0.000 & 0.010 & No \\
100-5-3b & 2 & 88287 & 67616 & 2 & 88287 & 68867 & 0 & 0 & 1251 & 0.000 & 0.060 & No \\
200-10-1 & 3 & 266151 & 217943 & 3 & 253840 & 225619 & 0 & -12311 & 7676 & 0.333 & — & — \\
200-10-1b & 3 & 266151 & 119823 & 3 & 253840 & 124124 & 0 & -12311 & 4301 & 0.333 & — & — \\
200-10-2 & 3 & 280370 & 169822 & 3 & 280370 & 170620 & 0 & 0 & 798 & 0.000 & 0.055 & No \\
200-10-2b & 3 & 280370 & 94577 & 3 & 280370 & 95974 & 0 & 0 & 1397 & 0.000 & 0.060 & No \\
200-10-3 & 3 & 272528 & 201601 & 3 & 272528 & 204726 & 0 & 0 & 3125 & 0.000 & 0.055 & No \\
200-10-3b & 3 & 272528 & 111465 & 3 & 234660 & 132379 & 0 & -37868 & 20914 & 0.333 & — & — \\
\bottomrule
\end{tabular}%
}
\end{table}

From Table~\ref{tab:flp_vrp} we observe several consistent patterns. FLP-VRP and NEO-DS usually open the same number of depots, and in many cases, the opened sets are identical. For example, in instances \texttt{20-5-1} and \texttt{50-5-2}, we have $\|\Delta \mathcal{D}\|_{\text{norm}} = 0$. Even when the depot sets coincide, the customer assignments may differ, leading to improved routing costs. This can be observed in instances \texttt{20-5-1b} and \texttt{100-10-2}, where $\|\Delta \mathcal{D}\|_{\text{norm}} = 0$ but $\|\Delta A\|_{\text{norm}} > 0$, and NEO-DS achieves lower routing costs ($\Delta R < 0$) through better allocations. This demonstrates that the neural surrogate captures routing structure beyond simple direct travel costs, which translates to more informed allocation decisions even when depot locations are identical. When the opened depots do differ, NEO-DS typically has lower cost, resulting in negative fixed-cost differences $\Delta F$, as seen in instances \texttt{50-5-3} and \texttt{200-10-1}. In some cases, such as \texttt{50-5-3}, NEO-DS consolidates by opening fewer depots ($\Delta n^{\text{open}} = -1$); in this case, there is increased routing cost with substantial fixed cost reduction ($\Delta F = -18993$), illustrating the surrogate's ability to capture the trade-off between depot opening and routing costs.

Furthermore,
\ifsupplementary
\else
in Figures~\ref{fig:inst_20_5_1}--\ref{fig:inst_50_5_2b} in Appendix~\ref{app:flpvrp_figures} we illustrate
\fi
four representative cases that highlight how the two methods differ in depot openings, customer-to-depot allocations and the resulting routing structure. These results establish that the neural surrogate within NEO-DS provides routing-aware cost approximation that leads to materially different and better location and allocation decisions compared to the direct-distance objective used in FLP-VRP. The improvements persist even after both methods use the same VRP solver to compute final routes, confirming that the value of NEO-LRP lies in the quality of its location-allocation decisions and is not ``washed out'' by the ex-post routing step.

\subsubsection{Effect of Problem Size}
\label{sec:problem_size}
\begin{figure}[!ht]
    \centering
    \begin{subfigure}{0.48\linewidth}
        \centering
        \includegraphics[width=\linewidth]{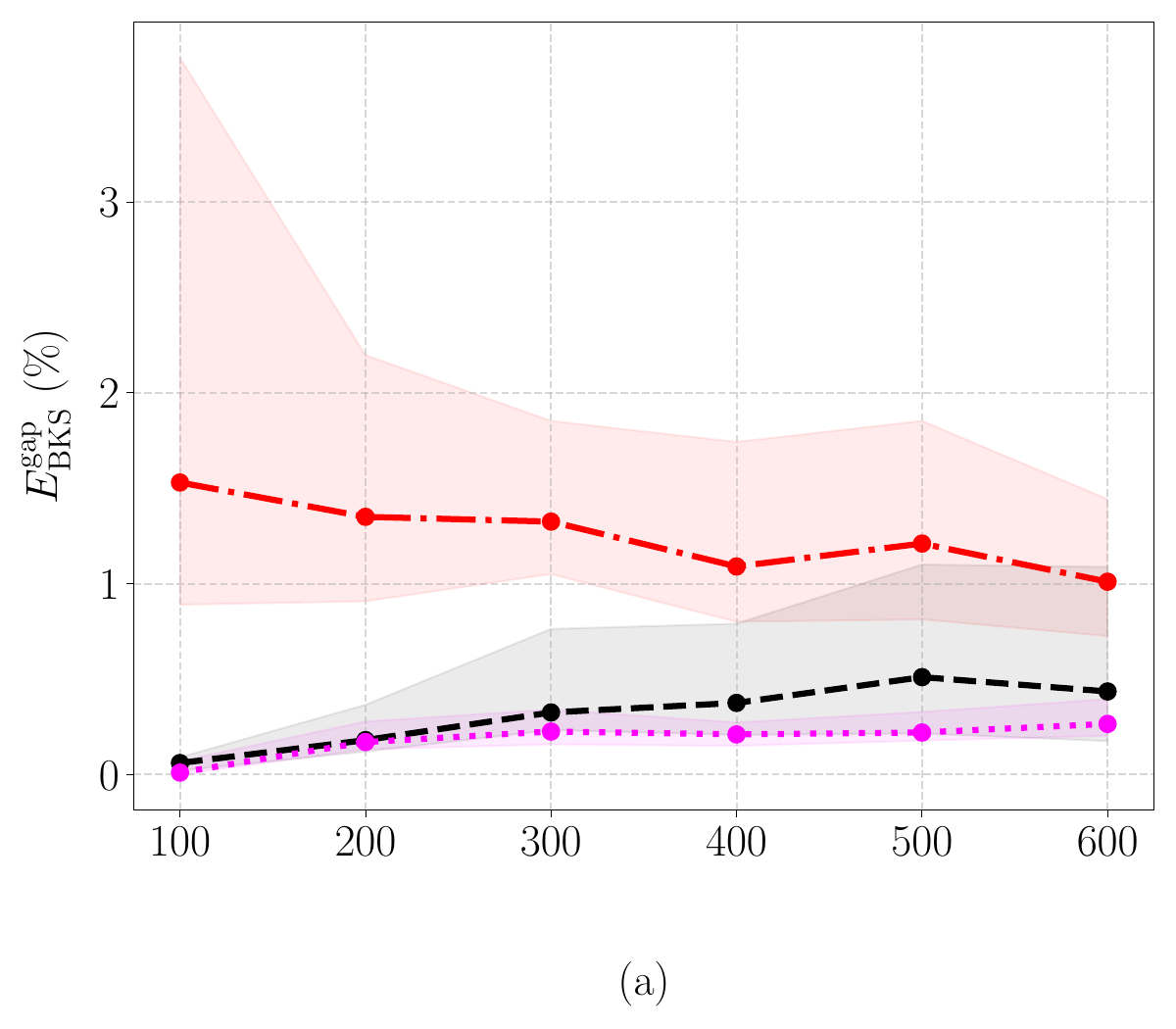}
        \label{fig:scalability_gap_schneider}
    \end{subfigure}
    \hfill
    \begin{subfigure}{0.48\linewidth}
        \centering
        \includegraphics[width=\linewidth]{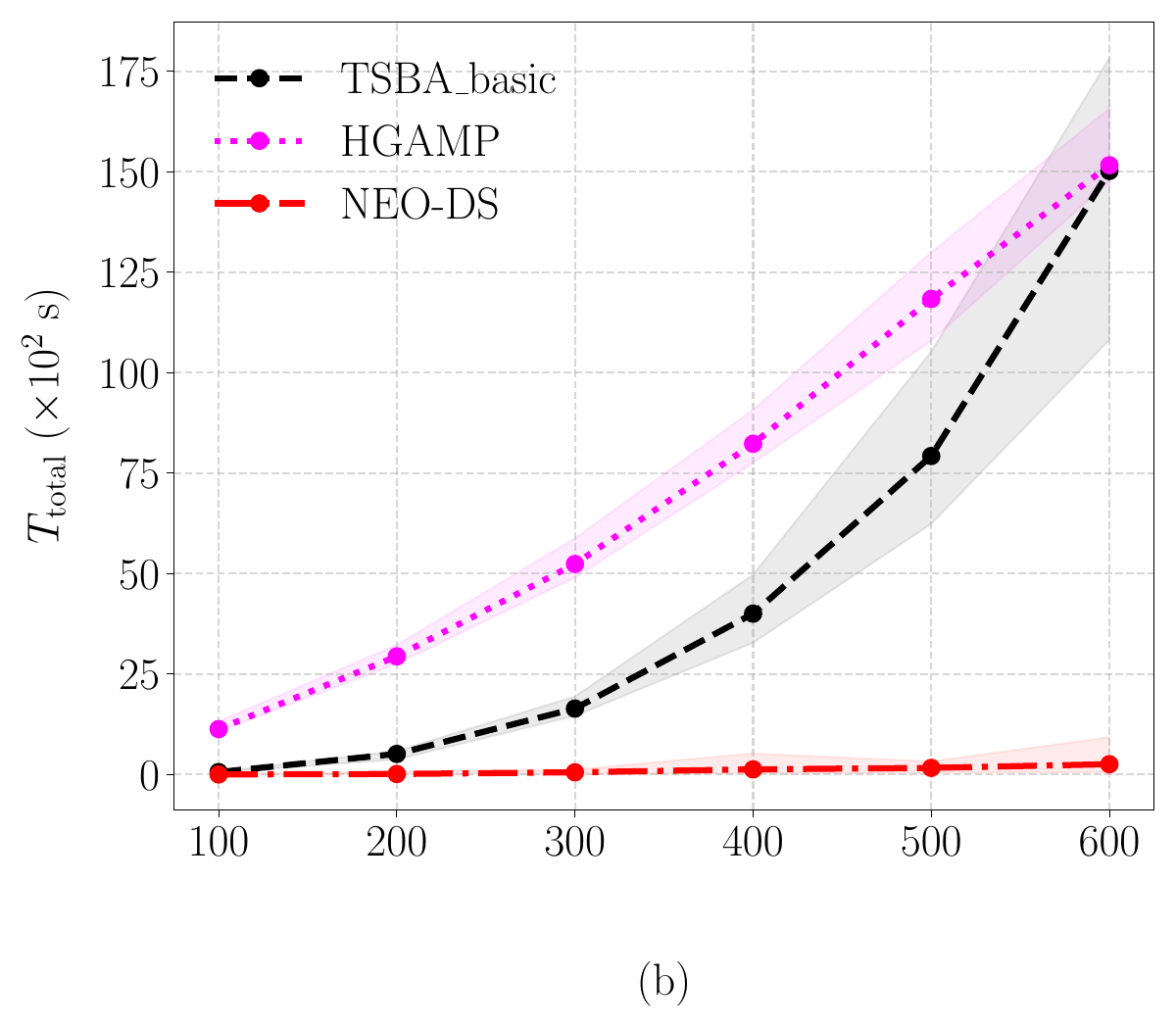}
        \label{fig:scalability_time_schneider}
    \end{subfigure}
    \vspace{-2.5em}
    \caption{Effect of problem size on the $\mathbb{S}$ benchmark set of \citet{schneider2019large}. (a) displays the median $E^{\text{gap}}_{\text{BKS}}$, while (b) shows the median $T_{\text{total}}$, both as a function of problem size (number of customers). Shaded regions represent the interquartile range.}
    \label{fig:scalability_schneider}
\end{figure}

To study the effect of problem size, we consider NEO-DS trained on data sampled per Section~\ref{sec:data_generation} with scaled labels, evaluated on the $\mathbb{S}$ benchmark set, as it contains the largest and most diverse instances. Figure~\ref{fig:scalability_schneider} illustrates the scalability of NEO-DS compared to baseline methods across different problem sizes, with shaded regions representing the interquartile range.

Figure~\ref{fig:scalability_schneider}a shows an improving $E^{\text{gap}}_{\text{BKS}}$ trend for NEO-DS as problem size increases, with median gaps decreasing from $1.53\%$ at 100 customers to $1.01\%$ at 600 customers. The interquartile range also narrows from [0.89\%, 3.75\%] to [0.72\%, 1.44\%], indicating an improving performance at larger sizes. This trend suggests that our neural embedded optimization approach becomes more effective as the problem scale increases. While the exact reason for this improvement is not fully clear, we hypothesize that larger problems provide a richer solution space, giving the surrogate more opportunities to guide the optimization toward regions with lower costs. In contrast, TSBA$_{\text{basic}}$ shows degrading quality, with median gaps increasing from $0.06\%$ at 100 customers to $0.43\%$ at 600 customers, while HGAMP maintains consistently good quality with median gaps ranging from $0.01\%$ to $0.27\%$ across all sizes.

Figure~\ref{fig:scalability_schneider}b shows NEO-DS's computational scalability. Median solution times grow from $1.7\,\text{s}$ at 100 customers to $256.3\,\text{s}$ at 600 customers. In contrast, both TSBA$_{\text{basic}}$ and HGAMP exhibit substantially longer times, with median times reaching $15{,}017\,\text{s}$ and $15{,}159\,\text{s}$ respectively for 600-customer instances. This represents computational speedups of $58.6\times$ over TSBA$_{\text{basic}}$ and $59.2\times$ over HGAMP, highlighting the value of our neural-embedded approach for large-scale location-routing problems. Our findings suggest that NEO-DS is particularly suited to large-scale instances, where it produces solutions close to the best known (a $1.01\%$ median gap) at a fraction of the runtime of these heuristics.

\ifsupplementary
\else
\ifinsupplementary
\section{Effect of Sample Size}
\label{app:sample_size}
\else
\subsubsection{Effect of Sample Size}
\label{sec:sample_size}
\fi
To study the effect of training sample size on the performance of NEO-DS, we consider the $\mathbb{P}$ benchmark set of \citet{prins2004nouveaux} and vary the amount of training data. Specifically, we train models using $\{100, 1{,}000, 10{,}000, 100{,}000\}$ training instances, with an additional 10\% of each training set allocated for validation.
\ifinsupplementary
All data are sampled with scaled labels generated using the VROOM solver.
\else
All data are sampled using the methodology of Section~\ref{sec:data_generation} with scaled labels generated using the VROOM solver (see Section~\ref{sec:computing_final_routes}).
\fi
We maintain a fixed test set of 10,000 instances across all experiments. After training, we evaluate the trained models by solving the embedded optimization problem and assessing four metrics: training error $E^{\mathrm{train}}$, test error $E^{\mathrm{test}}$, prediction error $E^{\mathrm{pred}}$, and optimization gap $E^{\mathrm{gap}}_{\mathrm{BKS}}$.
Here, $E^{\mathrm{pred}}$ quantifies the accuracy of the neural surrogate in predicting routing costs under the obtained location-allocation decisions. Letting $\mathcal{D}_{\text{open}} = \{ d \in \mathcal{V}_D : \hat{y}_d = 1 \}$ denote the set of open depots, the prediction error is defined as $E^{\text{pred}} = \sum_{d \in \mathcal{D}_{\text{open}}} |\hat{R}_d(\hat{w}) - R_d(\hat{w})| / R_d(\hat{w}) \times 100\%$, where $R_d(\hat{w})$ is the true routing cost computed by the VRP solver and $\hat{R}_d(\hat{w})$ is the surrogate prediction. Similarly, $E^{\mathrm{train}}$ and $E^{\mathrm{test}}$ are defined as the mean absolute percentage error of the surrogate, $\frac{1}{|\mathcal{D}|}\sum_{m} |\hat{R}(\mathcal{G}'^{(m)}) - R(\mathcal{G}'^{(m)})| / R(\mathcal{G}'^{(m)}) \times 100\%$, evaluated over the training and test datasets, respectively.

Figure~\ref{fig:samplesize_prodhon} illustrates how these metrics change with increasing training set size. Figure~\ref{fig:samplesize_prodhon}a and Figure~\ref{fig:samplesize_prodhon}b show that both $E^{\mathrm{train}}$ and $E^{\mathrm{test}}$ decrease as sample size increases. Specifically, the median $E^{\mathrm{train}}$ drops from 7.97\% at 110 samples to 2.54\% at 110{,}000 samples, while the median $E^{\mathrm{test}}$ decreases from 8.52\% to 2.50\% over the same range. The close alignment between training and test errors across all sample sizes provides empirical evidence that the model generalizes well without overfitting.
\ifinsupplementary\else
This is consistent with our use of benchmark-independent training data (Section~\ref{sec:data_generation}).
\fi
Additionally, the interquartile range narrows from [3.76\%, 13.40\%] to [1.14\%, 4.83\%] for training error and from [3.98\%, 15.01\%] to [1.11\%, 4.74\%] for test error.

Figure~\ref{fig:samplesize_prodhon}c shows that the median $E^{\mathrm{pred}}$ remains relatively stable across sample sizes, ranging from 11.52\% to 15.07\%, suggesting that the prediction error on individual instances is less sensitive to training set size. Finally, Figure~\ref{fig:samplesize_prodhon}d demonstrates that $E^{\mathrm{gap}}_{\mathrm{BKS}}$ improves with sample size, with the median decreasing from 3.52\% at 110 samples to 1.82\% at 11,000 samples. However, the improvement plateaus beyond 11,000 samples, with the median gap remaining at 1.62\% at 110{,}000 samples. The interquartile range also narrows from [2.13\%, 6.46\%] to [1.20\%, 3.09\%].
\ifinsupplementary
This observation suggests that increasing the training sample size beyond a certain point has a limited effect on the final solution quality.
\else
This observation suggests that increasing the training sample size beyond a certain point has a limited effect on the final solution quality in terms of the optimization gap $E^{\mathrm{gap}}_{\mathrm{BKS}}$.
\fi

Finally, we also highlight that a high $E^{\mathrm{pred}}$ does not translate into a large optimization gap $E^{\mathrm{gap}}_{\mathrm{BKS}}$. The reason is the a posteriori step of obtaining the final routes. NEO-LRP optimizes over the surrogate predictions $\hat{R}_d$ to select the location-allocation decisions $(\hat{y}, \hat{w})$, and then computes the true routing cost $R_d(S_d(\hat{w}))$ by solving a CVRP for each open depot (see Section~\ref{sec:computing_final_routes}). Therefore, a surrogate with a higher $E^{\mathrm{pred}}$ can still identify high-quality solutions with a lower $E^{\mathrm{gap}}_{\mathrm{BKS}}$.

\begin{figure}[!ht]
    \centering
    \includegraphics[width=0.95\linewidth]{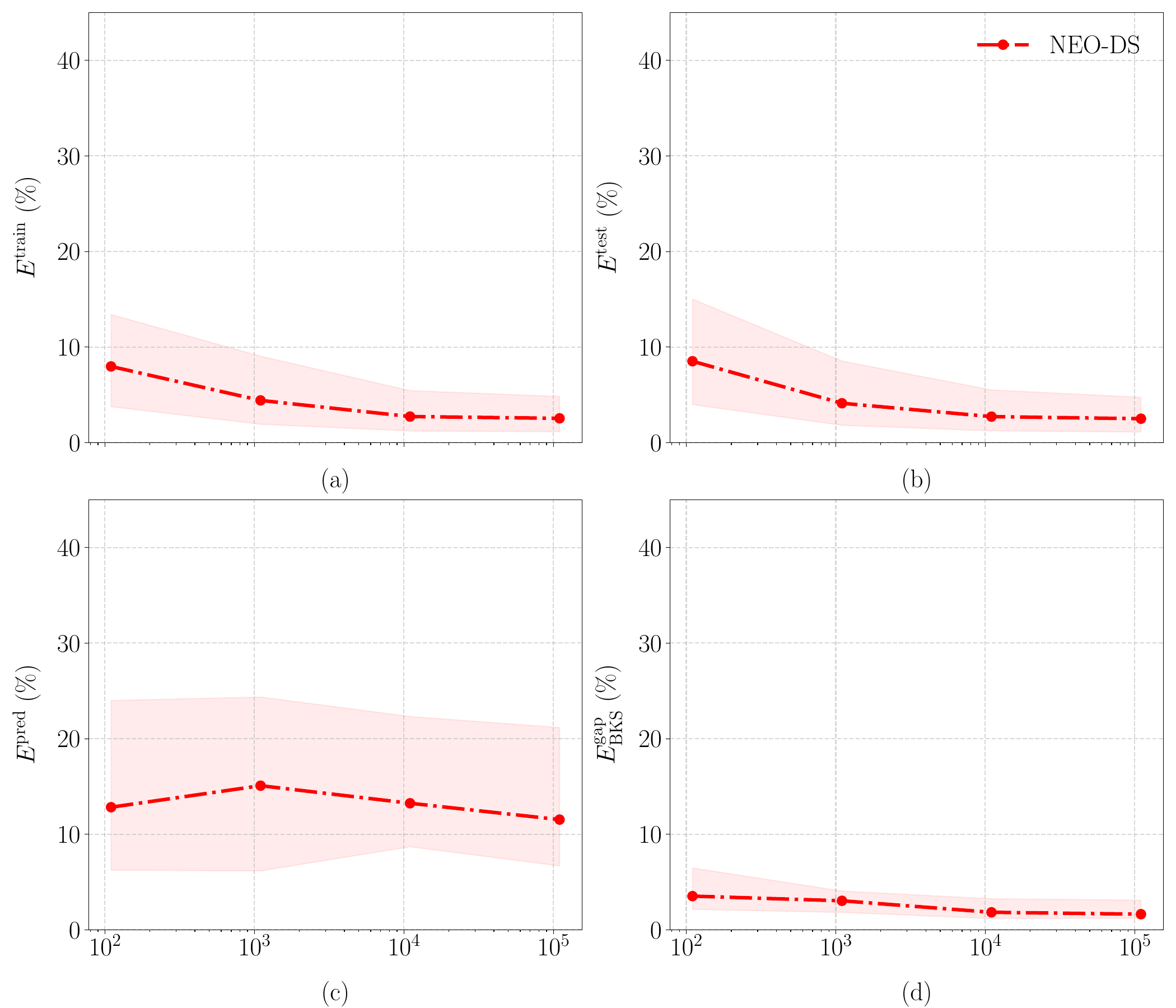}
    \vspace{-0.5em}
    \caption{Effect of training sample size on NEO-DS performance using the $\mathbb{P}$ benchmark set of \citet{prins2004nouveaux}. (a) shows training error $E^{\mathrm{train}}$, (b) shows test error $E^{\mathrm{test}}$, (c) shows prediction error $E^{\mathrm{pred}}$, and (d) shows optimization gap $E^{\mathrm{gap}}_{\mathrm{BKS}}$, all as functions of training set size. Shaded regions represent the interquartile range.}
    \label{fig:samplesize_prodhon}
\end{figure}

\ifinsupplementary
\section{Effect of Routing Solver}
\label{app:routing_solver}
\else
\subsubsection{Effect of Routing Solver}
\label{sec:routing_solver}
\fi
\ifinsupplementary
We evaluate the impact of the routing solver used for computing final routes \emph{a posteriori}.
\else
We evaluate the impact of the routing solver used for computing final routes \emph{a posteriori} (see Section \ref{sec:computing_final_routes}).
\fi
We compare three solvers: VROOM~\citep{vroom_v1.14} (heuristic), OR-Tools~\citep{ortools} (heuristic), and VRPSolverEasy~\citep{errami2024vrpsolvereasy} (exact branch-cut-and-price). VROOM and OR-Tools are each given a 5-second time limit per depot, while VRPSolverEasy is given a 1-hour time limit. We solve our neural-embedded MIP model to get the location-allocation decisions, followed by routing with the respective solver. Since the location-allocation MIP embeds the same trained neural surrogate across all runs, this study isolates the effect of routing solver choice on the framework's overall performance.

Table~\ref{tab:routing_solver_ablation} reports the average $E^{\text{gap}}_{\text{BKS}}$ and computation times averaged over all $N_C$ customer instances in the $\mathbb{P}$ benchmark. The location-allocation times ($T_{\text{LA}}$) are nearly identical across all solvers; this is expected as the same neural-embedded MIP is solved in each case. All three solvers produce identical solution quality on small instances ($n=20$), with gaps of 4.55\%. As problem size increases, VRPSolverEasy achieves better gaps (0.87\% vs.\ 1.01\% at $N_C=200$), but at higher computational times (94s vs.\ 18s). OR-Tools yields slightly higher gaps than VROOM, which shows the best trade-off between solution quality and speed.

We observe that solver choice has minimal impact on solution quality in terms of average $E^{\text{gap}}_{\text{BKS}}$, %
with pairwise differences of at most $1.19$ percentage points across all three solvers and instance sizes.
\ifinsupplementary
This finding reinforces that the primary driver
\else
This finding reinforces the conclusion from Section~\ref{sec:val_loc_alloc}: the primary driver
\fi
of NEO-LRP's performance is the quality of its location-allocation decisions, not the routing solver used to construct final routes. In practice, the exact solver may be preferred for smaller instances, whereas the heuristic solver is more appropriate for larger instances.
\ifinsupplementary
Detailed instance-level results are provided in Table~\ref{tab:detailed_routing_solver}.
\else
Detailed instance-level results are provided in
\ifsupplementary
the Supplementary Material.
\else
Table~\ref{tab:detailed_routing_solver} in Appendix~\ref{app:detailed_results}.
\fi
\fi

\begin{table*}[ht!]
\centering
\small
\caption{Summary of routing solver comparison on the $\mathbb{P}$ benchmark set of \citet{prins2004nouveaux}, averaged by instance size.}
\label{tab:routing_solver_ablation}
\begin{tabular}{l ccc ccc ccc}
\toprule
& \multicolumn{3}{c}{VROOM} & \multicolumn{3}{c}{OR-Tools} & \multicolumn{3}{c}{VRPSolverEasy} \\
\cmidrule(lr){2-4} \cmidrule(lr){5-7} \cmidrule(lr){8-10}
$N_C$ & $E^{\text{gap}}_{\text{BKS}}$ & $T_{\text{LA}}$ & $T_{\text{total}}$ & $E^{\text{gap}}_{\text{BKS}}$ & $T_{\text{LA}}$ & $T_{\text{total}}$ & $E^{\text{gap}}_{\text{BKS}}$ & $T_{\text{LA}}$ & $T_{\text{total}}$ \\
\midrule
20  & 4.55 & 0.17 & 0.18 & 4.55 & 0.19 & 12.70 & 4.55 & 0.17 & 0.26 \\
50  & 2.89 & 0.33 & 0.48 & 3.02 & 0.33 & 13.46 & 2.79 & 0.33 & 7.04 \\
100 & 1.76 & 3.74 & 4.81 & 2.26 & 3.35 & 16.68 & 1.61 & 3.39 & 56.69 \\
200 & 1.01 & 13.38 & 17.76 & 2.06 & 13.01 & 28.01 & 0.87 & 13.07 & 93.87 \\
\bottomrule
\end{tabular}
\end{table*}

\ifinsupplementary
\begin{table}[H]
\centering
\scriptsize
\caption{Detailed routing solver comparison results on the $\mathbb{P}$ benchmark of \citet{prins2004nouveaux}. BKS values are taken from \citet{loffler2023conceptually} and \citet{he2025hybrid}.}
\label{tab:detailed_routing_solver}
\renewcommand{\arraystretch}{0.9}%
\resizebox{\ifdim\width>\linewidth \linewidth\else \width\fi}{!}{%
\begin{tabular}{l r  ccc  ccc  ccc}
\toprule
Instance & BKS &
\multicolumn{3}{c}{VROOM (5s)} &
\multicolumn{3}{c}{OR-Tools (5s)} &
\multicolumn{3}{c}{VRPSolverEasy (3600s)} \\
\cmidrule(lr){3-5} \cmidrule(lr){6-8} \cmidrule(lr){9-11}
 & &
$E^{\text{gap}}_{\text{BKS}}$ & $T_{\text{LA}}$(s) & $T_{\text{total}}$(s) &
$E^{\text{gap}}_{\text{BKS}}$ & $T_{\text{LA}}$(s) & $T_{\text{total}}$(s) &
$E^{\text{gap}}_{\text{BKS}}$ & $T_{\text{LA}}$(s) & $T_{\text{total}}$(s) \\
\midrule
20-5-1a & 54793 & 3.24 & 0.19 & 0.20 & 3.24 & 0.28 & 15.30 & 3.24 & 0.21 & 0.36 \\ 
20-5-1b & 39104 & 3.65 & 0.17 & 0.19 & 3.65 & 0.18 & 10.18 & 3.65 & 0.17 & 0.23 \\ 
20-5-2a & 48908 & 5.82 & 0.15 & 0.16 & 5.82 & 0.15 & 15.15 & 5.82 & 0.14 & 0.22 \\ 
20-5-2b & 37542 & 5.49 & 0.16 & 0.19 & 5.49 & 0.17 & 10.18 & 5.49 & 0.16 & 0.26 \\ 
\midrule
\textbf{Average} &  & \textbf{4.55} & \textbf{0.17} & \textbf{0.18} & \textbf{4.55} & \textbf{0.19} & \textbf{12.70} & \textbf{4.55} & \textbf{0.17} & \textbf{0.26} \\ 
\midrule
50-5-1a & 90111 & 1.63 & 0.26 & 0.37 & 1.63 & 0.27 & 15.28 & 1.63 & 0.26 & 0.92 \\ 
50-5-1b & 63242 & 1.18 & 0.28 & 0.45 & 1.18 & 0.28 & 10.28 & 1.18 & 0.27 & 6.15 \\ 
50-5-2a & 88293 & 2.74 & 0.32 & 0.45 & 2.82 & 0.32 & 15.33 & 2.68 & 0.32 & 2.02 \\ 
50-5-2b & 67308 & 5.80 & 0.28 & 0.38 & 5.80 & 0.28 & 15.28 & 5.80 & 0.28 & 1.44 \\ 
50-5-2bBIS & 51822 & 4.86 & 0.46 & 0.53 & 4.86 & 0.47 & 15.47 & 4.86 & 0.47 & 40.41 \\ 
50-5-2BIS & 84055 & 1.62 & 0.25 & 0.45 & 2.83 & 0.26 & 15.26 & 1.36 & 0.25 & 0.79 \\ 
50-5-3a & 86203 & 0.96 & 0.37 & 0.59 & 0.61 & 0.38 & 10.38 & 0.50 & 0.37 & 2.61 \\ 
50-5-3b & 61830 & 4.32 & 0.40 & 0.60 & 4.45 & 0.39 & 10.40 & 4.32 & 0.39 & 1.98 \\ 
\midrule
\textbf{Average} &  & \textbf{2.89} & \textbf{0.33} & \textbf{0.48} & \textbf{3.02} & \textbf{0.33} & \textbf{13.46} & \textbf{2.79} & \textbf{0.33} & \textbf{7.04} \\ 
\midrule
100-5-1a & 274814 & 1.49 & 1.40 & 2.39 & 1.90 & 1.36 & 16.37 & 0.99 & 1.36 & 4.51 \\ 
100-5-1b & 213568 & 1.25 & 1.36 & 2.25 & 1.98 & 1.33 & 16.34 & 1.22 & 1.33 & 6.95 \\ 
100-5-2a & 193671 & 0.28 & 0.45 & 2.22 & 1.64 & 0.43 & 10.43 & -0.05 & 0.42 & 6.81 \\ 
100-5-2b & 157095 & 0.79 & 0.51 & 2.18 & 1.76 & 0.44 & 10.45 & 0.73 & 0.44 & 39.22 \\ 
100-5-3a & 200079 & 1.44 & 0.43 & 1.94 & 2.12 & 0.42 & 10.42 & 1.05 & 0.40 & 4.69 \\ 
100-5-3b & 152441 & 3.09 & 0.45 & 1.81 & 3.48 & 0.43 & 10.44 & 2.65 & 0.44 & 531.65 \\ 
100-10-1a & 287661 & 1.52 & 3.78 & 4.60 & 1.81 & 3.43 & 18.44 & 1.52 & 3.52 & 16.88 \\ 
100-10-1b & 230989 & 2.35 & 1.80 & 2.38 & 2.39 & 1.76 & 16.77 & 2.35 & 1.79 & 8.66 \\ 
100-10-2a & 243590 & 1.49 & 2.95 & 3.75 & 2.14 & 2.66 & 17.66 & 1.45 & 2.66 & 12.99 \\ 
100-10-2b & 203988 & 2.54 & 2.83 & 3.56 & 2.65 & 2.55 & 17.55 & 2.54 & 2.55 & 4.68 \\ 
100-10-3a & 250882 & 1.79 & 3.25 & 4.18 & 1.95 & 2.94 & 17.95 & 1.78 & 2.95 & 15.81 \\ 
100-10-3b & 203114 & 3.08 & 25.65 & 26.52 & 3.36 & 22.40 & 37.40 & 3.08 & 22.76 & 27.42 \\ 
\midrule
\textbf{Average} &  & \textbf{1.76} & \textbf{3.74} & \textbf{4.81} & \textbf{2.26} & \textbf{3.35} & \textbf{16.68} & \textbf{1.61} & \textbf{3.39} & \textbf{56.69} \\ 
\midrule
200-10-1a & 474702 & 1.00 & 5.89 & 10.85 & 2.25 & 5.03 & 20.03 & 0.67 & 5.13 & 130.54 \\ 
200-10-1b & 375177 & 0.74 & 60.58 & 64.80 & 2.16 & 60.40 & 75.41 & 0.70 & 60.35 & 123.47 \\ 
200-10-2a & 448077 & 0.65 & 2.28 & 6.40 & 1.50 & 2.14 & 17.15 & 0.50 & 2.20 & 50.80 \\ 
200-10-2b & 373696 & 0.71 & 3.02 & 6.65 & 1.39 & 2.82 & 17.83 & 0.69 & 2.89 & 25.64 \\ 
200-10-3a & 469433 & 1.67 & 5.58 & 10.51 & 2.40 & 5.06 & 20.07 & 1.40 & 5.15 & 176.21 \\ 
200-10-3b & 362320 & 1.30 & 2.96 & 7.33 & 2.68 & 2.60 & 17.60 & 1.28 & 2.67 & 56.56 \\ 
\midrule
\textbf{Average} &  & \textbf{1.01} & \textbf{13.38} & \textbf{17.76} & \textbf{2.06} & \textbf{13.01} & \textbf{28.01} & \textbf{0.87} & \textbf{13.07} & \textbf{93.87} \\ 
\bottomrule
\end{tabular}%
}
\end{table}

\fi

\fi
\subsubsection{Effect of Neural Architecture}
\label{sec:ablation_architecture}

The neural surrogate is a modular component of NEO-LRP, and a natural question is how solution quality depends on the choice of architecture. 
We compare the Deep Sets (DS) architecture used in our main experiments against Graph Transformers (GT) as an alternative neural surrogate. As discussed in Section~\ref{sec:learning}, DS computes independent node-level embeddings $h_u = \phi(z_u)$ and aggregates them via summation. This provides a simple permutation-invariant model that does not exploit the graph structure. In contrast, GT leverages attention-based aggregation over pairwise node interactions and can potentially capture richer spatial relationships between customers. We investigate whether this additional expressiveness translates to improved optimization performance. For GT, we use Weights \& Biases (W\&B)~\citep{wandb} for hyperparameter search and allocate approximately 24 hours for NEO-GT$^{\text{S}}$ trained with scaled labels ($N_2$ normalization). Final GT hyperparameters are provided in
\ifsupplementary
\else
Appendix~\ref{app:hyperparameters}.
\fi

Table~\ref{tab:neods_neogt_summary} presents the results averaged by instance size. NEO-DS outperforms NEO-GT across all instance sizes in both solution quality and computational time. NEO-DS achieves gaps ranging from 1.01\% to 4.55\%, while NEO-GT gaps range from 7.97\% to 15.50\%. As we mentioned earlier, in principle, GTs can capture richer spatial relationships through attention-based aggregation. However, note that this requires more sophisticated training procedures. Prior work on graph neural networks for routing problems has shown that effective training typically requires large datasets on the order of a million samples and careful hyperparameter tuning~\citep{joshi2019efficient}, whereas our training set contains on the order of a hundred thousand samples. In our experiments, we allocated a 12-hour budget for NEO-GT, which allowed evaluation of only 3 to 4 hyperparameter configurations. In contrast, the simpler Deep Sets architecture enabled exploration of over 25 configurations within a 6-hour budget. These results suggest that simpler architectures, which are faster to train, provide a more practical choice in our neural embedded optimization approach when training resources are limited. We leave further exploration of training procedures for GTs to future work. Detailed instance-level results are provided in
\ifsupplementary
\else
Table~\ref{tab:optimization_neods_neogt} in Appendix~\ref{app:detailed_results}.
\fi

\begin{table}[htbp]
\centering
\small
\setlength{\tabcolsep}{4pt}
\caption{Summary of NEO-DS and NEO-GT comparison on the $\mathbb{P}$ benchmark of \citet{prins2004nouveaux}, averaged by instance size.}
\label{tab:neods_neogt_summary}
\begin{tabular}{l ccc ccc}
\toprule
& \multicolumn{3}{c}{NEO-DS} & \multicolumn{3}{c}{NEO-GT} \\
\cmidrule(lr){2-4} \cmidrule(lr){5-7}
$N_C$ & $E^{\text{gap}}_{\text{BKS}}$ & $T^{\text{LA}}$ & $T_{\text{total}}$ & $E^{\text{gap}}_{\text{BKS}}$ & $T^{\text{LA}}$ & $T_{\text{total}}$ \\
\midrule
20  & 4.55 & 0.17 & 0.18 & 12.03 & 0.35 & 0.37 \\
50  & 2.89 & 0.33 & 0.48 & 15.50 & 2.15 & 2.32 \\
100 & 1.76 & 3.74 & 4.81 & 7.97 & 12.32 & 13.28 \\
200 & 1.01 & 13.38 & 17.76 & 10.18 & 40.34 & 44.44 \\
\bottomrule
\end{tabular}
\end{table}

\ifsupplementary
\else
\ifinsupplementary
\section{Effect of Target Normalization}
\label{app:target_normalization}
\else
\subsubsection{Effect of Target Normalization}
\label{sec:target_normalization}
\fi
\ifinsupplementary
To evaluate the impact of target normalization applied to the routing cost labels, we compare alternative normalization schemes.
\else
To evaluate the impact of target normalization applied to the routing cost labels, we compare alternative normalization schemes in addition to the default scheme defined in  Section~\ref{sec:feature_normalization}.
\fi
\ifinsupplementary
Let $\mathcal D = \{ (\mathcal G'^{(m)}, R(\mathcal G'^{(m)})) : m = 1,\dots,M \}$ denote the training dataset of CVRP instances, where $R^{(m)} = R(\mathcal G'^{(m)})$ denotes the CVRP cost of instance $\mathcal G'^{(m)}$, and $P^{(m)} = P(\mathcal G'^{(m)})$ is the scale factor.
\else
Recall that $\mathcal D = \{ (\mathcal G'^{(m)}, R(\mathcal G'^{(m)})) : m = 1,\dots,M \}$ denotes the training dataset of CVRP instances, where $R^{(m)} = R(\mathcal G'^{(m)})$ denotes the CVRP cost of instance $\mathcal G'^{(m)}$, and $P^{(m)} = P(\mathcal G'^{(m)})$ is the scale factor defined in Section~\ref{sec:feature_normalization}.
\fi
Let also $R_{\min}$ and $R_{\max}$ denote the minimum and maximum routing costs across the training dataset.
The four normalization schemes can then be defined as follows.
\[
\bar R^{(m)} =
\begin{cases}
R^{(m)} & (N_1)\;\text{raw costs}\\ %
\frac{R^{(m)}}{P^{(m)}} & (N_2)\;\text{cost-over-}P\\ %
\frac{R^{(m)} - R_{\min}}{R_{\max} - R_{\min}} & (N_3)\;\text{min-max}\\ %
\dfrac{\frac{R^{(m)}}{P^{(m)}} - \min_j \frac{R^{(j)}}{P^{(j)}}}
{\max_j \frac{R^{(j)}}{P^{(j)}} - \min_j \frac{R^{(j)}}{P^{(j)}}} & (N_4)\;\text{cost-over-}P\;\text{min-max}
\end{cases}
\]
We note that $N_2$ is the default normalization scheme used in our main experiments.

\ifinsupplementary
To evaluate these normalization schemes, we train separate NEO-DS surrogate models on 110,000 independently sampled VRP instances using each scheme.
\else
To evaluate these normalization schemes, we train separate NEO-DS surrogate models on 110,000 independently sampled VRP instances (Section~\ref{sec:data_generation}) using each scheme.
\fi
Each trained surrogate is then embedded into the location-allocation MIP and evaluated on the $\mathbb{P}$ benchmark set. Thus, all reported results reflect the final downstream impact of normalization.

Table~\ref{tab:neos_norm_comparison_summary} presents the results averaged by instance size. For small instances (20 customers), NEO-DS achieves similar gaps across all normalization schemes (4.20--4.59\%), with $N_4$ performing best at 4.20\% and $N_1$, $N_2$ showing 4.59\% and 4.55\% respectively.

\begin{table}[htbp]
\centering
\caption{Summary of normalization scheme comparison on the $\mathbb{P}$ benchmark of \citet{prins2004nouveaux}, averaged by instance size.}
\label{tab:neos_norm_comparison_summary}
\small
\setlength{\tabcolsep}{4pt}
\begin{tabular}{l cc cc cc cc}
\toprule
 & \multicolumn{8}{c}{NEO-DS} \\
\cmidrule(lr){2-9}
$N_C$
& \multicolumn{2}{c}{$N_1$} & \multicolumn{2}{c}{$N_2$} & \multicolumn{2}{c}{$N_3$} & \multicolumn{2}{c}{$N_4$} \\
\cmidrule(lr){2-3}\cmidrule(lr){4-5}\cmidrule(lr){6-7}\cmidrule(lr){8-9}
 & $E^{\text{gap}}_{\text{BKS}}$ & $T_{\text{total}}$
 & $E^{\text{gap}}_{\text{BKS}}$ & $T_{\text{total}}$
 & $E^{\text{gap}}_{\text{BKS}}$ & $T_{\text{total}}$
 & $E^{\text{gap}}_{\text{BKS}}$ & $T_{\text{total}}$ \\
\midrule
20  & 4.59 & 0.51 & 4.55 & 0.18 & 4.39 & 0.22 & 4.20 & 0.16 \\
50  & 4.24 & 1.29 & 2.89 & 0.48 & 4.59 & 0.62 & 3.18 & 0.35 \\
100 & 1.65 & 15.82 & 1.76 & 4.81 & 2.47 & 4.16 & 2.51 & 3.73 \\
200 & 0.96 & 35.37 & 1.01 & 17.76 & 1.53 & 15.65 & 1.51 & 8.21 \\
\midrule
Overall & 2.60 & 13.82 & 2.28 & 5.63 & 3.11 & 4.99 & 2.71 & 3.25 \\
\bottomrule
\end{tabular}
\end{table}

For medium instances (50 customers), $N_2$ achieves 2.89\% gap in 0.48s, outperforming $N_1$ (4.24\%, 1.29s) and $N_3$ (4.59\%, 0.62s). At 100 customers, $N_1$ and $N_2$ show comparable gaps (1.65\%, 1.76\%) but $N_2$ is significantly faster (4.81s vs 15.82s).

For large instances (200 customers), $N_2$ achieves the best gap of 1.01\% in 17.76s, substantially outperforming all other schemes. Averaging across all sizes, NEO-DS with $N_2$ achieves 2.28\% gap in 5.63s, compared to $N_1$ (2.60\%, 13.82s), $N_3$ (3.11\%, 4.99s), and $N_4$ (2.71\%, 3.25s). While $N_4$ has the fastest runtime, $N_2$ provides the best solution quality, with gaps improving from 4.55\% to 1.01\% as instance size increases from 20 to 200 customers.

These results demonstrate that the choice of target normalization during surrogate training has an impact on the final optimization performance of the neural embedded approach. The normalization scheme affects how well the surrogate learns to approximate routing costs, which in turn influences the quality of location-allocation decisions made by the embedded MIP.
\ifinsupplementary
Detailed instance-level results are provided in Table~\ref{tab:detailed_normalization}.

\begin{table}[H]
\centering
\caption{Detailed normalization scheme comparison results on the $\mathbb{P}$ benchmark of \citet{prins2004nouveaux}. BKS values are taken from \citet{loffler2023conceptually} and \citet{he2025hybrid}.}
\label{tab:detailed_normalization}
\scriptsize
\setlength{\tabcolsep}{3pt}
\renewcommand{\arraystretch}{0.9}%
\resizebox{\ifdim\width>\linewidth \linewidth\else \width\fi}{!}{%
\begin{tabular}{l r cc cc cc cc}
\toprule
 & & \multicolumn{8}{c}{\textbf{NEO-DS}} \\
\cmidrule(lr){3-10}
Instance & BKS
& \multicolumn{2}{c}{$N_1$} & \multicolumn{2}{c}{$N_2$} & \multicolumn{2}{c}{$N_3$} & \multicolumn{2}{c}{$N_4$} \\
\cmidrule(lr){3-4}\cmidrule(lr){5-6}\cmidrule(lr){7-8}\cmidrule(lr){9-10}
 &  & $E^{\text{gap}}_{\text{BKS}}$ & $T_{\text{total}}$
 & $E^{\text{gap}}_{\text{BKS}}$ & $T_{\text{total}}$
 & $E^{\text{gap}}_{\text{BKS}}$ & $T_{\text{total}}$
 & $E^{\text{gap}}_{\text{BKS}}$ & $T_{\text{total}}$ \\
\midrule
20-5-1a & 54793 & 3.24 & 0.51 & 3.24 & 0.20 & 3.24 & 0.18 & 3.24 & 0.20 \\
20-5-1b & 39104 & 6.30 & 0.51 & 3.65 & 0.19 & 4.04 & 0.27 & 6.30 & 0.11 \\
20-5-2a & 48908 & 3.32 & 0.50 & 5.82 & 0.16 & 6.25 & 0.21 & 3.22 & 0.15 \\
20-5-2b & 37542 & 5.49 & 0.52 & 5.49 & 0.19 & 4.04 & 0.23 & 4.04 & 0.20 \\
\midrule
\textbf{Average} &  & \textbf{4.59} & \textbf{0.51} & \textbf{4.55} & \textbf{0.18} & \textbf{4.39} & \textbf{0.22} & \textbf{4.20} & \textbf{0.16} \\
\midrule
50-5-1a & 90111 & 1.63 & 0.71 & 1.63 & 0.37 & 3.13 & 0.57 & 2.13 & 0.33 \\
50-5-1b & 63242 & 1.18 & 0.87 & 1.18 & 0.45 & 4.41 & 0.54 & 1.36 & 0.30 \\
50-5-2Bis & 84055 & 7.87 & 3.64 & 1.62 & 0.45 & 2.53 & 0.67 & 1.66 & 0.41 \\
50-5-2a & 88293 & 1.14 & 1.10 & 2.74 & 0.45 & 2.81 & 0.61 & 2.79 & 0.39 \\
50-5-2b & 67308 & 5.76 & 1.03 & 5.80 & 0.38 & 8.84 & 0.49 & 7.13 & 0.28 \\
50-5-2bBis & 51822 & 3.25 & 0.75 & 4.86 & 0.53 & 4.86 & 0.52 & 4.89 & 0.32 \\
50-5-3a & 86203 & 8.96 & 1.20 & 0.96 & 0.59 & 5.80 & 0.86 & 2.82 & 0.40 \\
50-5-3b & 61830 & 4.14 & 1.06 & 4.32 & 0.60 & 4.32 & 0.70 & 2.61 & 0.38 \\
\midrule
\textbf{Average} &  & \textbf{4.24} & \textbf{1.29} & \textbf{2.89} & \textbf{0.48} & \textbf{4.59} & \textbf{0.62} & \textbf{3.18} & \textbf{0.35} \\
\midrule
100-10-1a & 287661 & 2.48 & 52.89 & 1.52 & 4.60 & 3.48 & 6.00 & 1.87 & 1.29 \\
100-10-1b & 230989 & 1.76 & 32.99 & 2.35 & 2.38 & 5.90 & 2.17 & 4.35 & 1.35 \\
100-10-2a & 243590 & 1.49 & 15.97 & 1.49 & 3.75 & 1.33 & 15.14 & 1.75 & 1.60 \\
100-10-2b & 203988 & 1.99 & 15.77 & 2.54 & 3.56 & 2.90 & 2.99 & 2.73 & 1.56 \\
100-10-3a & 250882 & 3.23 & 34.80 & 1.79 & 4.18 & 3.92 & 3.24 & 4.12 & 15.29 \\
100-10-3b & 203114 & 3.14 & 22.84 & 3.08 & 26.52 & 2.60 & 2.64 & 8.41 & 13.05 \\
100-5-1a & 274814 & 1.42 & 2.45 & 1.49 & 2.39 & 1.36 & 6.60 & 0.82 & 1.67 \\
100-5-1b & 213568 & 1.50 & 2.40 & 1.25 & 2.25 & 1.82 & 2.30 & 1.63 & 2.05 \\
100-5-2a & 193671 & 0.28 & 2.58 & 0.28 & 2.22 & 0.28 & 2.40 & 0.67 & 1.91 \\
100-5-2b & 157095 & 0.79 & 2.59 & 0.79 & 2.18 & 2.28 & 2.37 & 0.45 & 1.80 \\
100-5-3a & 200079 & 0.76 & 2.43 & 1.44 & 1.94 & 1.89 & 2.33 & 2.00 & 1.79 \\
100-5-3b & 152441 & 1.01 & 2.13 & 3.09 & 1.81 & 1.94 & 1.76 & 1.37 & 1.42 \\
\midrule
\textbf{Average} &  & \textbf{1.65} & \textbf{15.82} & \textbf{1.76} & \textbf{4.81} & \textbf{2.47} & \textbf{4.16} & \textbf{2.51} & \textbf{3.73} \\
\midrule
200-10-1a & 474702 & 1.58 & 31.81 & 1.00 & 10.85 & 1.60 & 10.06 & 0.94 & 6.13 \\
200-10-1b & 375177 & 0.75 & 27.36 & 0.74 & 64.80 & 0.99 & 21.46 & 1.07 & 21.15 \\
200-10-2a & 448077 & 0.44 & 15.14 & 0.65 & 6.40 & 1.79 & 13.07 & 0.56 & 4.92 \\
200-10-2b & 373696 & 0.54 & 10.22 & 0.71 & 6.65 & 0.38 & 13.28 & 1.46 & 5.71 \\
200-10-3a & 469433 & 1.30 & 65.47 & 1.67 & 10.51 & 1.28 & 18.38 & 1.04 & 5.48 \\
200-10-3b & 362320 & 1.18 & 62.22 & 1.30 & 7.33 & 3.15 & 17.62 & 4.01 & 5.86 \\
\midrule
\textbf{Average} &  & \textbf{0.96} & \textbf{35.37} & \textbf{1.01} & \textbf{17.76} & \textbf{1.53} & \textbf{15.65} & \textbf{1.51} & \textbf{8.21} \\
\bottomrule
\end{tabular}%
}
\end{table}

\else
Detailed instance-level results are provided in
\ifsupplementary
the Supplementary Material.
\else
Table~\ref{tab:detailed_normalization} in Appendix~\ref{app:detailed_results}.
\fi
\fi

\fi

\makeatletter
\@ifundefined{ifanon}{\newif\ifanon}{}
\makeatother

\section{Conclusion}
\label{sec:conclusions}

Our paper offers a fresh approach to combine machine learning with integer programming to solve location-routing problems. Our proposed NEO-LRP approach first approximates the vehicle routing cost associated with each open facility using a sparse neural surrogate and then embeds this surrogate into an easy-to-solve MIP model for location-allocation decisions. A key benefit of our framework is that it avoids potentially costly, time-consuming, and complex implementations that are typical of state-of-the-art heuristics while also being able to generalize to very large problems with other types of side constraints. Our experiments across four benchmark sets showed that NEO-LRP can achieve high-quality solutions quickly and scales well to larger instances with up to 600 customers.%

In terms of future work, our methodology can be extended in several directions. While the current neural network training focuses on minimizing prediction error, incorporating a ranking-preserving loss function could help identify promising solutions that might be overestimated and prematurely discarded during the neural-embedded optimization. This approach would ensure the correct ordering of solution costs even when absolute predictions are not highly accurate. Another direction is to investigate alternative neural architectures for the graph transformer regressor to improve solution quality while maintaining computational efficiency. Finally, the modular design of NEO-LRP allows for extension to other location-routing variants such as the Multi-Depot VRP and Two-Echelon CVRP and the integration of additional constraints, including time windows, fleet heterogeneity, or environmental considerations, making it adaptable to a wider range of real-world logistics settings.

\ifanon\else
\section*{Acknowledgments}
We acknowledge support from the United States Department of Energy, award number DE-SC0023361, and the National Research Foundation of Korea, award number RS-2023-00259550. Computations for this research were performed on the Pennsylvania State University's Institute for
Computational and Data Sciences' Roar supercomputer.
\fi

\bibliographystyle{plainnat}%
\bibliography{references}

@article{larsen2024fast,
	author = {Larsen, Eric and Frejinger, Emma and Gendron, Bernard and Lodi, Andrea},
	journal = {INFORMS Journal on Computing},
	number = {1},
	pages = {203--223},
	publisher = {Informs},
	title = {Fast continuous and integer {L}-shaped heuristics through supervised learning},
	volume = {36},
	year = {2024}}

@article{bogyrbayeva2024machine,
	author = {Bogyrbayeva, Aigerim and Meraliyev, Meraryslan and Mustakhov, Taukekhan and Dauletbayev, Bissenbay},
	journal = {IEEE Transactions on Intelligent Transportation Systems},
	number = {6},
	pages = {4754--4772},
	publisher = {IEEE},
	title = {Machine learning to solve vehicle routing problems: A survey},
	volume = {25},
	year = {2024}}

@article{ceccon2022omlt,
	author = {Ceccon, Francesco and Jalving, Jordan and Haddad, Joshua and Thebelt, Alexander and Tsay, Calvin and Laird, Carl D and Misener, Ruth},
	journal = {The Journal of Machine Learning Research},
	number = {1},
	pages = {15829--15836},
	publisher = {JMLRORG},
	title = {{OMLT}: Optimization \& machine learning toolkit},
	volume = {23},
	year = {2022}}

@article{grimstad2019relu,
	author = {Grimstad, Bjarne and Andersson, Henrik},
	journal = {Computers \& Chemical Engineering},
	pages = {106580},
	publisher = {Elsevier},
	title = {{ReLU} networks as surrogate models in mixed-integer linear programs},
	volume = {131},
	year = {2019}}

@article{tjeng2017evaluating,
	author = {Tjeng, Vincent and Xiao, Kai and Tedrake, Russ},
	journal = {arXiv preprint arXiv:1711.07356},
	title = {Evaluating robustness of neural networks with mixed integer programming},
	year = {2017}}

@inproceedings{prins2004nouveaux,
	author = {Prins, Christian and Prodhon, Caroline and Calvo, Roberto Wolfler},
	booktitle = {MOSIM'04 (4{\`e}me Conf. Francophone de Mod{\'e}lisation et Simulation)},
	title = {Nouveaux algorithmes pour le probl{\`e}me de localisation et routage avec contraintes de capacit{\'e}},
	year = {2004}}

@article{mara2021location,
	author = {Mara, Setyo Tri Windras and Kuo, RJ and Asih, Anna Maria Sri},
	journal = {International Transactions in Operational Research},
	number = {6},
	pages = {2941--2983},
	publisher = {Wiley Online Library},
	title = {Location-routing problem: a classification of recent research},
	volume = {28},
	year = {2021}}

@article{drexl2015survey,
	author = {Drexl, Michael and Schneider, Michael},
	journal = {European Journal of Operational Research},
	number = {2},
	pages = {283--308},
	publisher = {Elsevier},
	title = {A survey of variants and extensions of the location-routing problem},
	volume = {241},
	year = {2015}}

@incollection{albaredasambola2020location,
	address = {Cham, Switzerland},
	author = {Albareda-Sambola, Maria and Rodr{\'\i}guez-Pereira, Jessica},
	booktitle = {{Location Science}},
	edition = {2},
	editor = {Laporte, Gilbert and Nickel, Stefan and Saldanha-da-Gama, Francisco},
	pages = {431-451},
	publisher = {Springer},
	title = {{Location-Routing and Location-Arc Routing}},
	year = {2020}}

@inproceedings{zaheer2017deep,
	author = {Zaheer, Manzil and Kottur, Satwik and Ravanbhakhsh, Siamak and P{\'o}czos, Barnab{\'a}s and Salakhutdinov, Ruslan and Smola, Alexander J},
	booktitle = {Proceedings of the 31st International Conference on Neural Information Processing Systems},
	pages = {3394--3404},
	title = {{Deep Sets}},
	year = {2017}}

@article{loffler2023conceptually,
	author = {L{\"o}ffler, Maximilian and Bartolini, Enrico and Schneider, Michael},
	journal = {EURO Journal on Computational Optimization},
	pages = {100063},
	publisher = {Elsevier},
	title = {A conceptually simple algorithm for the capacitated location-routing problem},
	volume = {11},
	year = {2023}}

@article{schneider2017survey,
	author = {Schneider, Michael and Drexl, Michael},
	journal = {Annals of Operations Research},
	pages = {389--414},
	publisher = {Springer},
	title = {A survey of the standard location-routing problem},
	volume = {259},
	year = {2017}}

@book{laporte2019introduction,
	address = {Cham, Switzerland},
	author = {Laporte, Gilbert and Nickel, Stefan and Saldanha-da-Gama, Francisco},
	edition = {2},
	publisher = {Springer},
	title = {{Location Science}},
	year = {2020}}

@article{berger2007location,
	author = {Berger, Rosemary T and Coullard, Collette R and Daskin, Mark S},
	journal = {Transportation Science},
	number = {1},
	pages = {29--43},
	publisher = {INFORMS},
	title = {Location-routing problems with distance constraints},
	volume = {41},
	year = {2007}}

@inproceedings{akca2009branch,
	author = {Akca, Z and Berger, RT and Ralphs, TK},
	booktitle = {Operations Research and Cyber-Infrastructure},
	organization = {Springer},
	pages = {309--330},
	title = {A branch-and-price algorithm for combined location and routing problems under capacity restrictions},
	year = {2009}}

@article{belenguer2011branch,
	author = {Belenguer, Jos{\'e}-Manuel and Benavent, Enrique and Prins, Christian and Prodhon, Caroline and Calvo, Roberto Wolfler},
	journal = {Computers \& Operations Research},
	number = {6},
	pages = {931--941},
	publisher = {Elsevier},
	title = {A branch-and-cut method for the capacitated location-routing problem},
	volume = {38},
	year = {2011}}

@article{baldacci2011exact,
	author = {Baldacci, Roberto and Mingozzi, Aristide and Wolfler Calvo, Roberto},
	journal = {Operations Research},
	number = {5},
	pages = {1284--1296},
	publisher = {INFORMS},
	title = {An exact method for the capacitated location-routing problem},
	volume = {59},
	year = {2011}}

@article{contardo2013computational,
	author = {Contardo, Claudio and Cordeau, Jean-Fran{\c{c}}ois and Gendron, Bernard},
	journal = {Discrete Optimization},
	number = {4},
	pages = {263--295},
	publisher = {Elsevier},
	title = {A computational comparison of flow formulations for the capacitated location-routing problem},
	volume = {10},
	year = {2013}}

@article{contardo2014exact,
	author = {Contardo, Claudio and Cordeau, Jean-Fran{\c{c}}ois and Gendron, Bernard},
	journal = {INFORMS Journal on Computing},
	number = {1},
	pages = {88--102},
	publisher = {INFORMS},
	title = {An exact algorithm based on cut-and-column generation for the capacitated location-routing problem},
	volume = {26},
	year = {2014}}

@book{farahani2009facility,
	author = {Farahani, Reza Zanjirani and Hekmatfar, Masoud},
	publisher = {Springer Science \& Business Media},
	title = {Facility location: concepts, models, algorithms and case studies},
	year = {2009}}

@article{prins2007solving,
	author = {Prins, Christian and Prodhon, Caroline and Ruiz, Angel and Soriano, Patrick and Wolfler Calvo, Roberto},
	journal = {Transportation Science},
	number = {4},
	pages = {470--483},
	publisher = {INFORMS},
	title = {Solving the capacitated location-routing problem by a cooperative Lagrangean relaxation-granular tabu search heuristic},
	volume = {41},
	year = {2007}}

@article{ozyurt2007solving,
	author = {{\"O}zyurt, Zeynep and Aksen, Deniz},
	journal = {Extending the horizons: Advances in computing, optimization, and decision technologies},
	pages = {125--144},
	publisher = {Springer},
	title = {Solving the multi-depot location-routing problem with lagrangian relaxation},
	year = {2007}}

@article{escobar2013two,
	author = {Escobar, John Willmer and Linfati, Rodrigo and Toth, Paolo},
	journal = {Computers \& Operations Research},
	number = {1},
	pages = {70--79},
	publisher = {Elsevier},
	title = {A two-phase hybrid heuristic algorithm for the capacitated location-routing problem},
	volume = {40},
	year = {2013}}

@article{alvim2013popmusic,
	author = {Alvim, Adriana CF and Taillard, {\'E}ric D},
	journal = {EURO Journal on Transportation and Logistics},
	number = {3},
	pages = {231--254},
	publisher = {Elsevier},
	title = {{POPMUSIC} for the world location-routing problem},
	volume = {2},
	year = {2013}}

@article{chan2005multiple,
	author = {Chan, Yupo and Baker, Steven F},
	journal = {Mathematical and Computer Modelling},
	number = {8-9},
	pages = {1035--1053},
	publisher = {Elsevier},
	title = {The multiple depot, multiple traveling salesmen facility-location problem: Vehicle range, service frequency, and heuristic implementations},
	volume = {41},
	year = {2005}}

@inproceedings{bouhafs2006combination,
	author = {Bouhafs, Lyamine and Hajjam, Amir and Koukam, Abder},
	booktitle = {Knowledge-Based Intelligent Information and Engineering Systems: 10th International Conference, KES 2006, Bournemouth, UK, October 9-11, 2006. Proceedings, Part I 10},
	organization = {Springer},
	pages = {409--416},
	title = {A combination of simulated annealing and ant colony system for the capacitated location-routing problem},
	year = {2006}}

@inproceedings{sahraeian2009using,
	author = {Sahraeian, Rashed and Nadizadeh, Ali},
	booktitle = {XIII Congreso de Ingenier{\'\i}a de Organizaci{\'o}n},
	pages = {1721--1729},
	title = {Using greedy clustering method to solve capacitated location-routing problem},
	year = {2009}}

@article{ting2013multiple,
	author = {Ting, Ching-Jung and Chen, Chia-Ho},
	journal = {International Journal of Production Economics},
	number = {1},
	pages = {34--44},
	publisher = {Elsevier},
	title = {A multiple ant colony optimization algorithm for the capacitated location routing problem},
	volume = {141},
	year = {2013}}

@article{prins2006solving,
	author = {Prins, Christian and Prodhon, Caroline and Calvo, Roberto Wolfler},
	journal = {{4OR}},
	pages = {221--238},
	publisher = {Springer},
	title = {Solving the capacitated location-routing problem by a GRASP complemented by a learning process and a path relinking},
	volume = {4},
	year = {2006}}

@article{jabal2011variable,
	author = {Jabal-Ameli, MS and Aryanezhad, MB and Ghaffari-Nasab, N},
	journal = {International Journal of Industrial Engineering Computations},
	number = {1},
	pages = {141--154},
	title = {A variable neighborhood descent based heuristic to solve the capacitated location-routing problem},
	volume = {2},
	year = {2011}}

@article{jarboui2013variable,
	author = {Jarboui, Bassem and Derbel, Houda and Hanafi, Sa{\"\i}d and Mladenovi{\'c}, Nenad},
	journal = {Computers \& Operations Research},
	number = {1},
	pages = {47--57},
	publisher = {Elsevier},
	title = {Variable neighborhood search for location routing},
	volume = {40},
	year = {2013}}

@inproceedings{prins2006memetic,
	author = {Prins, Christian and Prodhon, Caroline and Calvo, Roberto Wolfler},
	booktitle = {European Conference on Evolutionary Computation in Combinatorial Optimization},
	organization = {Springer},
	pages = {183--194},
	title = {A memetic algorithm with population management ({MA|PM}) for the capacitated location-routing problem},
	year = {2006}}

@inproceedings{duhamel2008memetic,
	author = {Duhamel, Christophe and Lacomme, Philippe and Prins, Christian and Prodhon, Caroline},
	booktitle = {Proceedings of the 9th EU/Meeting on Metaheuristics for Logistics and Vehicle Routing, Troyes, France},
	pages = {39},
	title = {A memetic approach for the capacitated location routing problem},
	volume = {38},
	year = {2008}}

@article{duhamel2010grasp,
	author = {Duhamel, Christophe and Lacomme, Philippe and Prins, Christian and Prodhon, Caroline},
	journal = {Computers \& Operations Research},
	number = {11},
	pages = {1912--1923},
	publisher = {Elsevier},
	title = {A {GRASP$\times$ELS} approach for the capacitated location-routing problem},
	volume = {37},
	year = {2010}}

@article{derbel2012genetic,
	author = {Derbel, Houda and Jarboui, Bassem and Hanafi, Sa{\"\i}d and Chabchoub, Habib},
	journal = {Expert Systems with Applications},
	number = {3},
	pages = {2865--2871},
	publisher = {Elsevier},
	title = {Genetic algorithm with iterated local search for solving a location-routing problem},
	volume = {39},
	year = {2012}}

@book{toth2014vehicle,
	author = {Toth, Paolo and Vigo, Daniele},
	publisher = {SIAM},
	title = {Vehicle routing: problems, methods, and applications},
	year = {2014}}

@article{errami2024vrpsolvereasy,
	author = {Errami, Najib and Queiroga, Eduardo and Sadykov, Ruslan and Uchoa, Eduardo},
	journal = {INFORMS Journal on Computing},
	number = {4},
	pages = {956--965},
	publisher = {INFORMS},
	title = {{VRPSolverEasy: a Python library for the exact solution of a rich vehicle routing problem}},
	volume = {36},
	year = {2024}}

@misc{Gurobi2021,
	author = {{Gurobi Optimization, LLC}},
	howpublished = {\url{https://www.gurobi.com}},
	title = {{Gurobi Optimizer Reference Manual}},
	year = {2021}}

@article{paszke2019pytorch,
	author = {Paszke, Adam and Gross, Sam and Massa, Francisco and Lerer, Adam and Bradbury, James and Chanan, Gregory and Killeen, Trevor and Lin, Zeming and Gimelshein, Natalia and Antiga, Luca and others},
	journal = {Advances in Neural Information Processing Systems},
	title = {{PyTorch}: An imperative style, high-performance deep learning library},
	volume = {32},
	year = {2019}}

@article{arnold2021progressive,
	author = {Arnold, Florian and S{\"o}rensen, Kenneth},
	journal = {Computers \& Operations Research},
	pages = {105166},
	publisher = {Elsevier},
	title = {A progressive filtering heuristic for the location-routing problem and variants},
	volume = {129},
	year = {2021}}

@article{salhi1989effect,
	author = {Salhi, Said and Rand, Graham K},
	journal = {European Journal of Operational Research},
	number = {2},
	pages = {150--156},
	publisher = {Elsevier},
	title = {The effect of ignoring routes when locating depots},
	volume = {39},
	year = {1989}}

@article{min1998combined,
	author = {Min, Hokey and Jayaraman, Vaidyanathan and Srivastava, Rajesh},
	journal = {European Journal of Operational Research},
	number = {1},
	pages = {1--15},
	publisher = {Elsevier},
	title = {Combined location-routing problems: A synthesis and future research directions},
	volume = {108},
	year = {1998}}

@misc{ortools,
	author = {Laurent Perron and Vincent Furnon},
	howpublished = {\url{https://developers.google.com/optimization/}},
	organization = {Google},
	title = {{OR-Tools}},
	year = {2023}}

@article{fischetti2018deep,
	author = {Fischetti, Matteo and Jo, Jason},
	journal = {Constraints},
	number = {3},
	pages = {296--309},
	publisher = {Springer},
	title = {Deep neural networks and mixed integer linear optimization},
	volume = {23},
	year = {2018}}

@article{cappart2023combinatorial,
	author = {Cappart, Quentin and Ch{\'e}telat, Didier and Khalil, Elias B and Lodi, Andrea and Morris, Christopher and Veli{\v{c}}kovi{\'c}, Petar},
	journal = {Journal of Machine Learning Research},
	number = {130},
	pages = {1--61},
	title = {Combinatorial optimization and reasoning with graph neural networks},
	volume = {24},
	year = {2023}}

@inproceedings{NEURIPS2022_9793671e,
	author = {Patel, Rahul Mihir and Dumouchelle, Justin and Khalil, Elias and Bodur, Merve},
	booktitle = {Advances in Neural Information Processing Systems},
	pages = {23992--24005},
	publisher = {Curran Associates, Inc.},
	title = {{Neur2SP}: Neural Two-Stage Stochastic Programming},
	volume = {35},
	year = {2022}}

@article{bengio2021machine,
	author = {Bengio, Yoshua and Lodi, Andrea and Prouvost, Antoine},
	journal = {European Journal of Operational Research},
	number = {2},
	pages = {405--421},
	publisher = {Elsevier},
	title = {Machine learning for combinatorial optimization: a methodological tour d'horizon},
	volume = {290},
	year = {2021}}

@inproceedings{NIPS2017_d9896106,
	author = {Khalil, Elias and Dai, Hanjun and Zhang, Yuyu and Dilkina, Bistra and Song, Le},
	booktitle = {Advances in Neural Information Processing Systems},
	publisher = {Curran Associates, Inc.},
	title = {Learning Combinatorial Optimization Algorithms over Graphs},
	volume = {30},
	year = {2017}}

@article{uchoa2017new,
	author = {Uchoa, Eduardo and Pecin, Diego and Pessoa, Artur and Poggi, Marcus and Vidal, Thibaut and Subramanian, Anand},
	journal = {European Journal of Operational Research},
	number = {3},
	pages = {845--858},
	publisher = {Elsevier},
	title = {New benchmark instances for the capacitated vehicle routing problem},
	volume = {257},
	year = {2017}}

@article{sobhanan2024genetic,
	author = {Sobhanan, Abhay and Park, Junyoung and Park, Jinkyoo and Kwon, Changhyun},
	journal = {Transportation Science},
	number = {2},
	pages = {322--339},
	publisher = {INFORMS},
	title = {Genetic algorithms with neural cost predictor for solving hierarchical vehicle routing problems},
	volume = {59},
	year = {2025}}

@article{varol2024neural,
	author = {Varol, Taha and {\"O}zener, Okan {\"O}rsan and Albey, Erin{\c{c}}},
	journal = {Transportation Science},
	number = {1},
	pages = {45--66},
	publisher = {Informs},
	title = {Neural network estimators for optimal tour lengths of traveling salesperson problem instances with arbitrary node distributions},
	volume = {58},
	year = {2024}}

@inproceedings{kaleem2024neural,
	author = {Kaleem, Waquar and Ayala, Harshita and Subramanyam, Anirudh},
	booktitle = {IISE Annual Conference. Proceedings},
	organization = {Institute of Industrial and Systems Engineers (IISE)},
	pages = {1--6},
	title = {Neural Embedded Optimization for Integrated Location and Routing Problems},
	year = {2024}}

@article{prodhon2014survey,
	author = {Prodhon, Caroline and Prins, Christian},
	journal = {European Journal of Operational Research},
	number = {1},
	pages = {1--17},
	publisher = {Elsevier},
	title = {A survey of recent research on location-routing problems},
	volume = {238},
	year = {2014}}

@article{cuda2015survey,
	author = {Cuda, Rosario and Guastaroba, Gianfranco and Speranza, Maria Grazia},
	journal = {Computers \& Operations Research},
	pages = {185--199},
	publisher = {Elsevier},
	title = {A survey on two-echelon routing problems},
	volume = {55},
	year = {2015}}

@inproceedings{queiroga202110,
	author = {Queiroga, Eduardo and Sadykov, Ruslan and Uchoa, Eduardo and Vidal, Thibaut},
	booktitle = {AAAI-22 workshop on machine learning for operations research (ML4OR)},
	title = {10,000 optimal {CVRP} solutions for testing machine learning based heuristics},
	year = {2021}}

@misc{gurobi_ml,
	author = {{Gurobi Optimization LLC}},
	howpublished = {\url{https://github.com/Gurobi/gurobi-machinelearning}},
	title = {{gurobi-machinelearning}: A {Python} Package for Mixed-Integer Programming Formulations of Trained Machine Learning Models},
	year = {2024}}

@article{hemmelmayr2012adaptive,
	author = {Hemmelmayr, Vera C and Cordeau, Jean-Fran{\c{c}}ois and Crainic, Teodor Gabriel},
	journal = {Computers \& Operations Research},
	number = {12},
	pages = {3215--3228},
	publisher = {Elsevier},
	title = {An adaptive large neighborhood search heuristic for two-echelon vehicle routing problems arising in city logistics},
	volume = {39},
	year = {2012}}

@article{tuzun1999two,
	author = {Tuzun, Dilek and Burke, Laura I},
	journal = {European Journal of Operational Research},
	number = {1},
	pages = {87--99},
	publisher = {Elsevier},
	title = {A two-phase tabu search approach to the location routing problem},
	volume = {116},
	year = {1999}}

@manual{vroom_v1.14,
	address = {Besan{\c c}on, France},
	author = {Coupey, Julien and Nicod, Jean-Marc and Varnier, Christophe},
	note = {\url{http://vroom-project.org/}},
	organization = {Verso (\url{https://verso-optim.com/})},
	title = {{VROOM v1.14, Vehicle Routing Open-source Optimization Machine}},
	year = 2024}

@inproceedings{joshi2019efficient,
	author = {Joshi, Chaitanya K and Laurent, Thomas and Bresson, Xavier},
	booktitle = {INFORMS Annual Meeting},
	note = {arXiv preprint arXiv:1906.01227},
	title = {An Efficient Graph Convolutional Network Technique for the Travelling Salesman Problem},
	url = {https://arxiv.org/abs/1906.01227},
	year = {2019}}

@article{velivckovic2023everything,
	author = {Veli{\v{c}}kovi{\'c}, Petar},
	journal = {Current Opinion in Structural Biology},
	pages = {102538},
	publisher = {Elsevier},
	title = {Everything is connected: Graph neural networks},
	volume = {79},
	year = {2023}}

@article{bronstein2021geometric,
	author = {Bronstein, Michael M and Bruna, Joan and Cohen, Taco and Veli{\v{c}}kovi{\'c}, Petar},
	journal = {arXiv preprint arXiv:2104.13478},
	title = {Geometric deep learning: Grids, groups, graphs, geodesics, and gauges},
	year = {2021}}

@article{yun2019graph,
	author = {Yun, Seongjun and Jeong, Minbyul and Kim, Raehyun and Kang, Jaewoo and Kim, Hyunwoo J},
	journal = {Advances in neural information processing systems},
	title = {Graph transformer networks},
	volume = {32},
	year = {2019}}

@article{liguori2023nonrobust,
	author = {Liguori, Pedro Henrique and Mahjoub, A Ridha and Marques, Guillaume and Sadykov, Ruslan and Uchoa, Eduardo},
	journal = {Operations Research},
	number = {5},
	pages = {1577--1595},
	publisher = {INFORMS},
	title = {Nonrobust strong knapsack cuts for capacitated location routing and related problems},
	volume = {71},
	year = {2023}}

@article{lopes2016simple,
	author = {Lopes, Rui Borges and Ferreira, Carlos and Santos, Beatriz Sousa},
	journal = {Computers \& Operations Research},
	pages = {155--162},
	publisher = {Elsevier},
	title = {A simple and effective evolutionary algorithm for the capacitated location--routing problem},
	volume = {70},
	year = {2016}}

@misc{wandb,
	author = {Biewald, Lukas},
	howpublished = {\url{https://www.wandb.com/}},
	title = {Experiment Tracking with {Weights and Biases}},
	year = {2020}}

@article{Egele2025,
	author = {Romain Egele and Prasanna Balaprakash and Gavin M. Wiggins and Brett Eiffert},
	doi = {10.21105/joss.07975},
	journal = {Journal of Open Source Software},
	number = {109},
	pages = {7975},
	publisher = {The Open Journal},
	title = {{DeepHyper}: A {Python} Package for Massively Parallel Hyperparameter Optimization in Machine Learning},
	url = {https://doi.org/10.21105/joss.07975},
	volume = {10},
	year = {2025}}

@phdthesis{barreto_2004,
	advisor = {Carlos, Ferreira, and Pinto, Paix\~{a}o, Jos\'{e}},
	author = {Barreto, S\'{e}rgio dos Santos},
	isbn = {9798515247706},
	note = {AAI28477693},
	school = {Universidade de Aveiro (Portugal)},
	title = {An\'{a}lise e modeliza\c{c}\~{a}o De Problemas De localiza\c{c}\~{a}o-distribui\c{c}\~{a}o},
	year = {2004}}

@article{schneider2019large,
	author = {Schneider, Michael and L{\"o}ffler, Maximilian},
	journal = {Transportation Science},
	number = {1},
	pages = {301--318},
	publisher = {INFORMS},
	title = {Large composite neighborhoods for the capacitated location-routing problem},
	volume = {53},
	year = {2019}}

@article{he2025hybrid,
	author = {He, Pengfei and Hao, Jin-Kao and Wu, Qinghua},
	journal = {INFORMS Journal on Computing},
	publisher = {INFORMS},
	title = {A hybrid genetic algorithm with multi-population for capacitated location routing},
	year = {2025}}

@article{voigt2022hybrid,
  title={Hybrid adaptive large neighborhood search for vehicle routing problems with depot location decisions},
  author={Voigt, Stefan and Frank, Markus and Fontaine, Pirmin and Kuhn, Heinrich},
  journal={Computers \& Operations Research},
  volume={146},
  pages={105856},
  year={2022},
  publisher={Elsevier}
}

\newpage

\begin{appendices}

\newpage

\section{Hyperparameters}
\label{app:hyperparameters}

\begin{table}[h]
\caption{Hyperparameters, search ranges, and selected values for NEO-DS$^{\text{S}}$ and NEO-DS$^{\text{U}}$.}\label{tab:ds_hyperparams}%
\resizebox{\ifdim\width>\linewidth \linewidth\else \width\fi}{!}{%
\begin{tabular}{@{}llcc@{}}
\toprule
\textbf{Hyperparameter} & \textbf{Range} & \textbf{NEO-DS$^{\text{S}}$} & \textbf{NEO-DS$^{\text{U}}$}\\
\midrule
Latent dimension & \{4, 6, 8\} & 6 & 6\\
\# hidden layers in $\phi$ & \{2, 3, 4, 5, 6\} & 5 & 3\\
\# hidden layers in $\rho$ & \{1\} & 1 & 1\\
Hidden units per layer in $\phi$ & \{32, 64, 128, 256, 512, 1024, 2048\} & 32 & 1024\\
Hidden units in $\rho$ & \{4, 6, 8\} & 6 & 6\\
Batch size & \{32\} & 32 & 32\\
Learning rate & \{0.001\} & 0.001 & 0.001\\
Optimizer & Adam & Adam & Adam\\
Loss function & MSE & MSE & MSE\\
Activation function & ReLU & ReLU & ReLU\\
Early stopping patience & \{15, 20\} & 20 & 15\\
Epochs & \{50, 100, 200, 400, 600, 800, 1000\} & 200 & 600\\
\bottomrule
\end{tabular}%
}
\end{table}

\begin{table}[h]
\caption{Hyperparameters, search ranges, and selected values for NEO-GT$^{\text{S}}$ and NEO-GT$^{\text{U}}$.}\label{tab:gt_hyperparams}%
\resizebox{\ifdim\width>\linewidth \linewidth\else \width\fi}{!}{%
\begin{tabular}{@{}llcc@{}}
\toprule
\textbf{Hyperparameter} & \textbf{Range} & \textbf{NEO-GT$^{\text{S}}$} & \textbf{NEO-GT$^{\text{U}}$}\\
\midrule
Encoding dimension & \{8, 16, 32, 64, 128\} & 16 & 128\\
Latent dimension & \{6\} & 6 & 6\\
\# GT layers & \{2, 3, 4, 5\} & 5 & 4\\
\# hidden layers in $\rho$ & \{1\} & 1 & 1\\
Hidden units in $\rho$ & \{8\} & 8 & 8\\
\# attention heads & \{4, 8\} & 4 & 4\\
Dropout rate & \{0.0, 0.1, 0.2, 0.3, 0.4, 0.5\} & 0.3 & 0.3\\
Normalization & \{Graph Norm, Batch Norm, Layer Norm\} & Graph Norm & Graph Norm\\
Activation function & \{ELU, ReLU, Leaky ReLU\} & ELU & ELU\\
Beta & \{True, False\} & False & True\\
Batch size & \{8, 16, 32, 64, 128\} & 8 & 8\\
Initial learning rate & \{0.1, 0.01, 0.001, 0.0001\} & 0.001 & 0.001\\
Optimizer & AdamW (schedule-free) & AdamW & AdamW\\
Loss function & \{MSE, Huber, Smooth L1\} & Huber & Smooth L1\\
Decode method & pool & pool & pool\\
Epochs & 150 (fixed) & 50 & 50\\
\bottomrule
\end{tabular}%
}
\end{table}

\section{FLP-VRP vs NEO-DS Visualizations}
\label{app:flpvrp_figures}
This appendix provides visualizations comparing FLP-VRP and NEO-DS solutions on representative instances from the $\mathbb{P}$ benchmark set of \citet{prins2004nouveaux}.

\begin{figure}[H]
\centering
\includegraphics[width=1\linewidth]{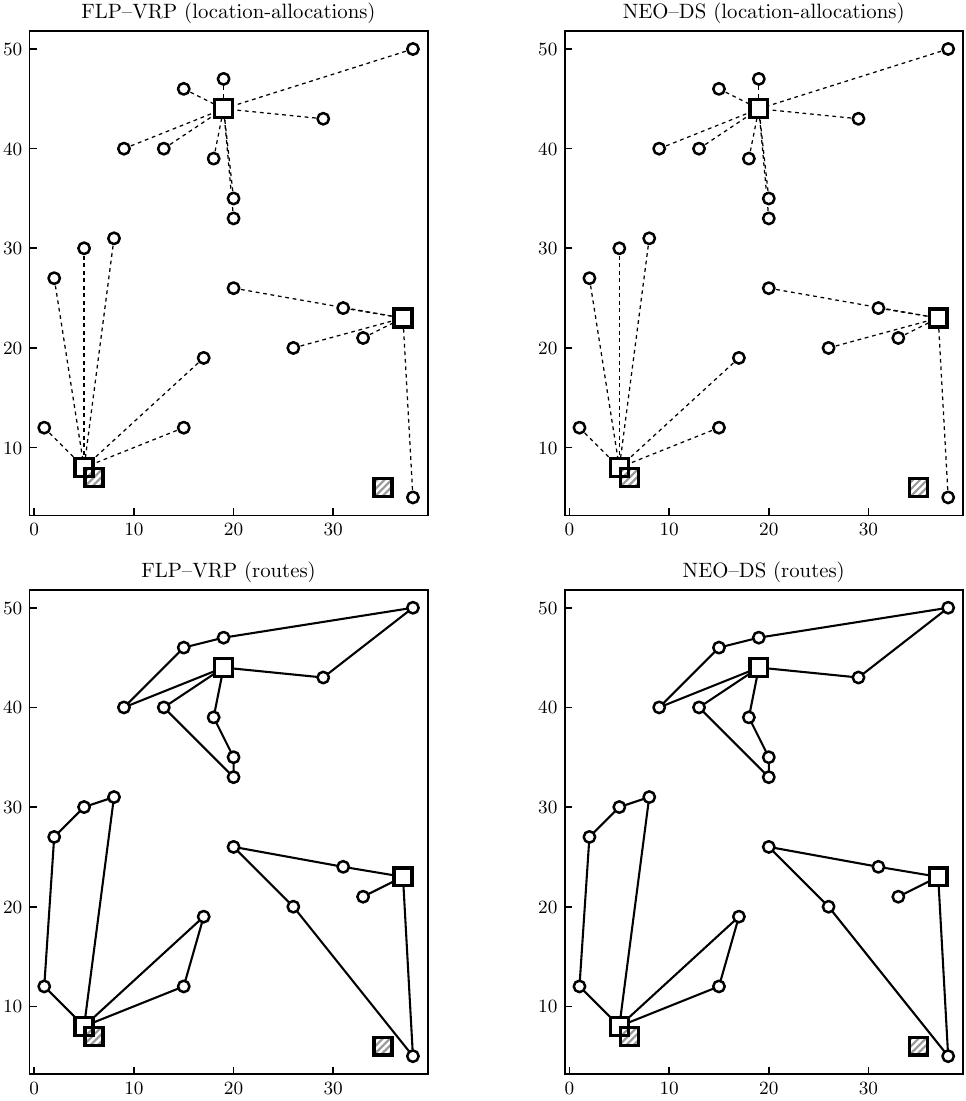}
\caption{
FLP-VRP and NEO-DS results for instance \texttt{20-5-1} of $\mathbb{P}$ benchmark.  
Both methods open the same depot set ($\|\Delta \mathcal{D}\|_{\text{norm}} = 0$) and produce identical customer assignments 
($\|\Delta A\|_{\text{norm}} = 0$)}
\label{fig:inst_20_5_1}
\end{figure}

\begin{figure}[H]
\centering
\includegraphics[width=1\linewidth]{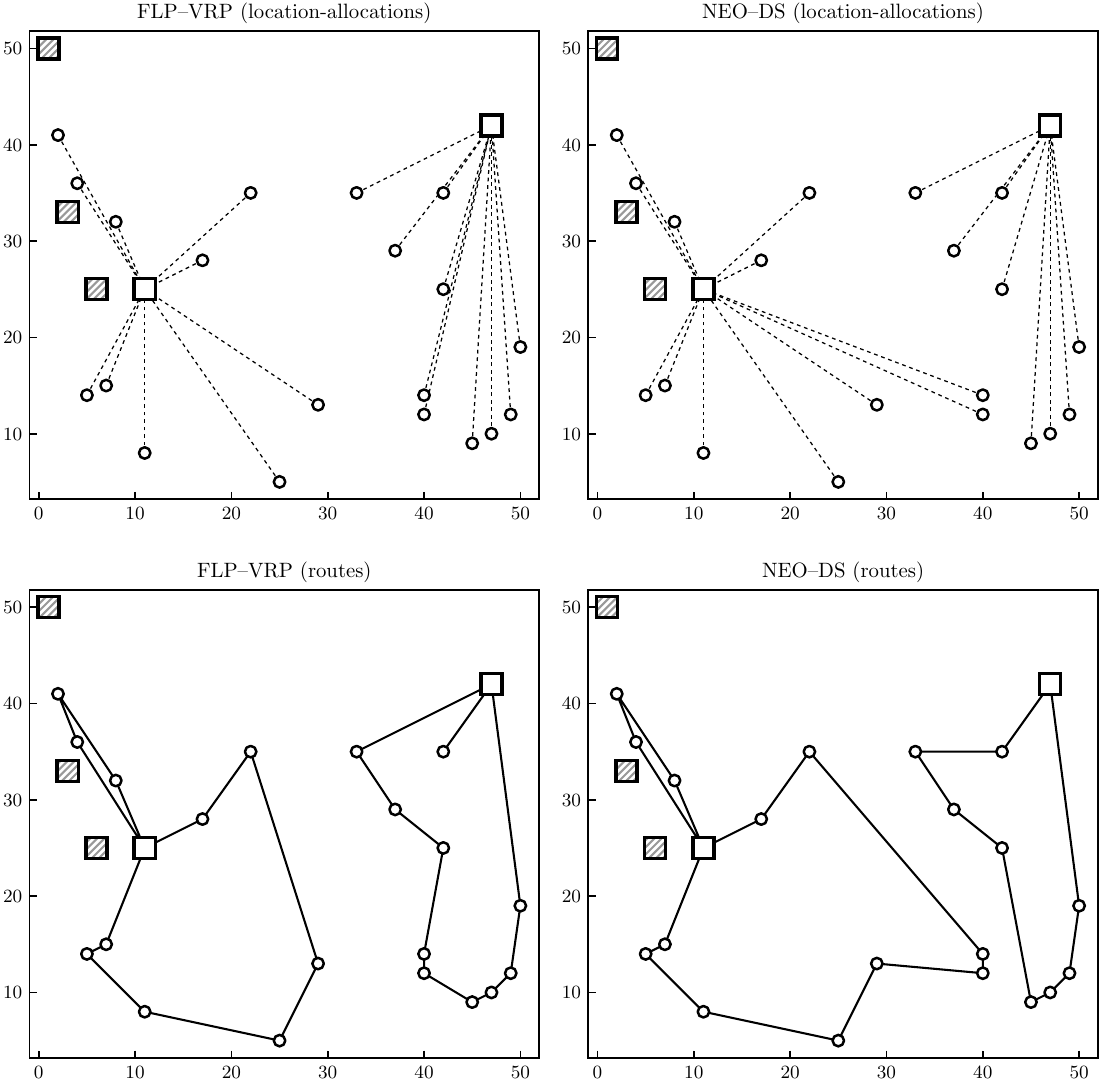}
\caption{
FLP-VRP and NEO-DS results for instance \texttt{20-5-1b} of $\mathbb{P}$ benchmark.  
The depot set is identical across methods ($\|\Delta \mathcal{D}\|_{\text{norm}} = 0$),  
but the customer-to-depot assignments differ ($\|\Delta A\|_{\text{norm}} = 0.10$).  
NEO-DS uses 3 routes compared to 4 for FLP-VRP, achieving a lower routing cost ($\Delta R = -1038$). This illustrates how the surrogate informs assignment decisions to reduce the number of routes.}
\label{fig:inst_20_5_1b}
\end{figure}

\begin{figure}[H]
\centering
\includegraphics[width=1\linewidth]{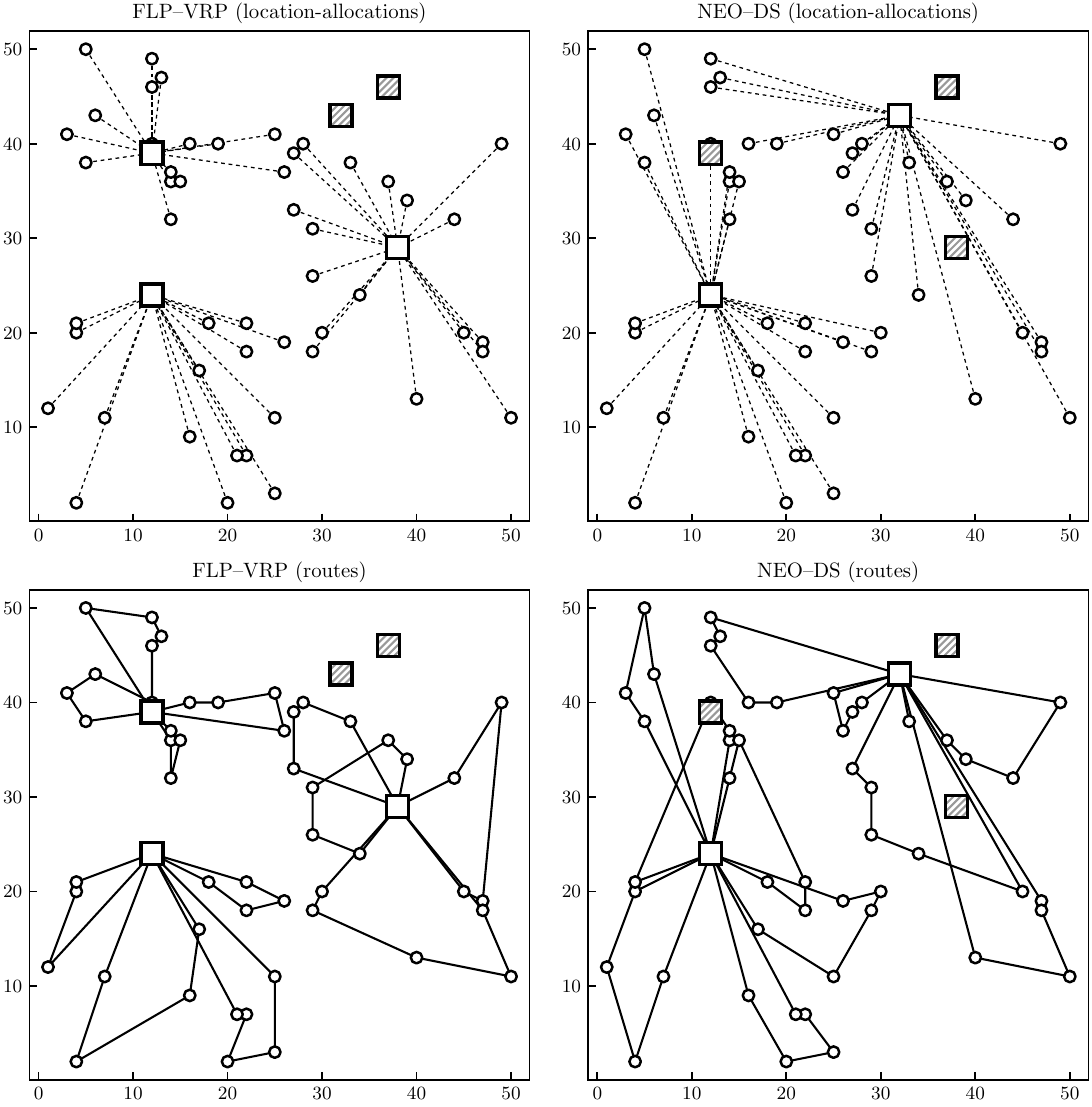}
\caption{
FLP-VRP and NEO-DS results for instance \texttt{50-5-3} of $\mathbb{P}$ benchmark.  
NEO-DS opens one fewer depot ($\Delta n^{\text{open}} = -1$), with $\Delta F = -18993$ and $\Delta R = +10655$, 
yielding a net cost reduction. This illustrates how the surrogate balances depot opening against routing costs.}
\label{fig:inst_50_5_3}
\end{figure}

\begin{figure}[H]
\centering
\includegraphics[width=1\linewidth]{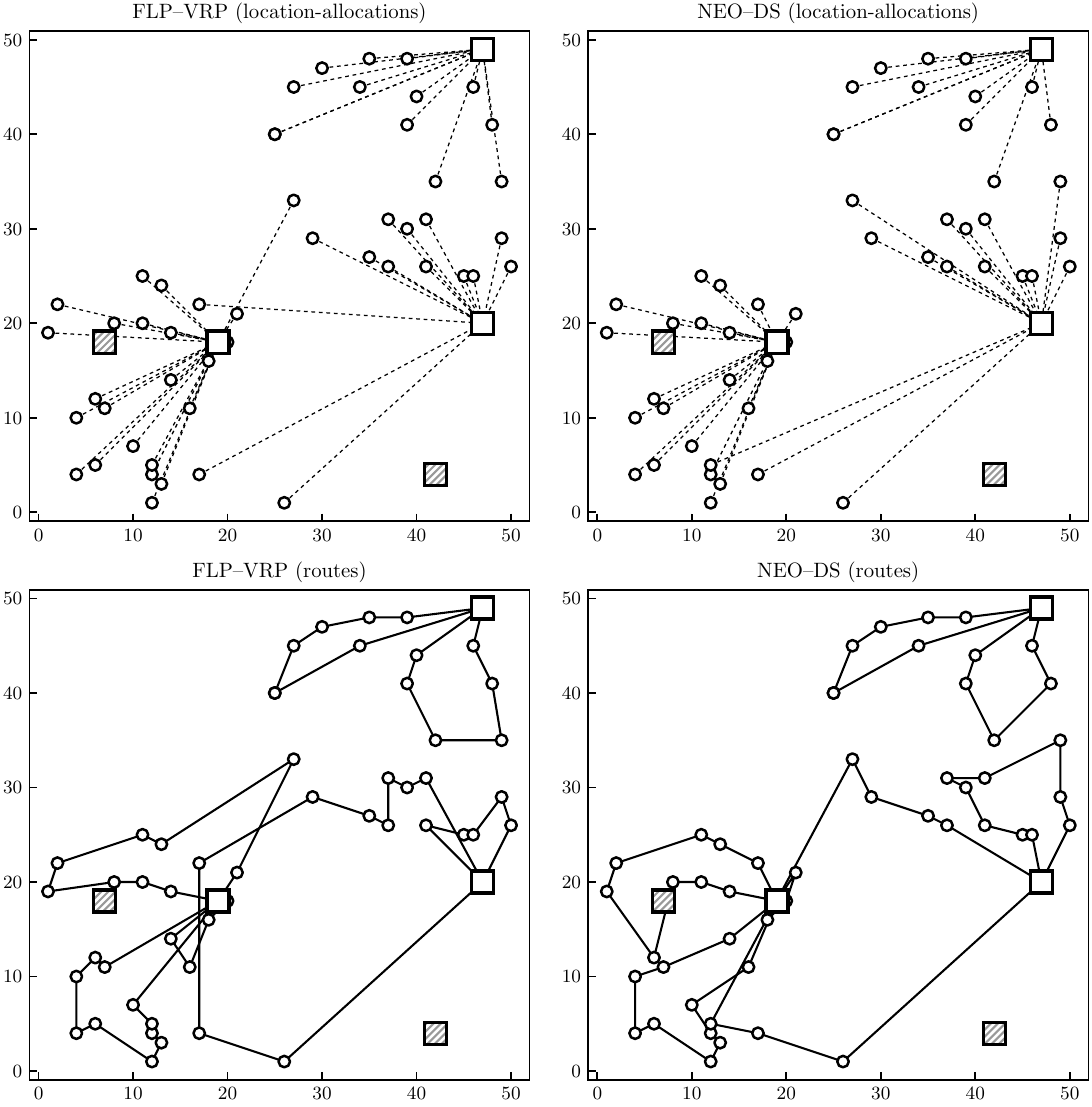}
\caption{
FLP-VRP and NEO-DS results for instance \texttt{50-5-2b} of $\mathbb{P}$ benchmark.  
Both methods open the same depot set ($\|\Delta \mathcal{D}\|_{\text{norm}} = 0$),  
yet their customer-to-depot assignments differ ($\|\Delta A\|_{\text{norm}} = 0.08$).
NEO-DS achieves lower routing cost ($\Delta R = -1113$), demonstrating that the neural surrogate guides allocation decisions toward more cost efficient routes.}
\label{fig:inst_50_5_2b}
\end{figure}
\FloatBarrier

\section{Detailed Results}
\label{app:detailed_results}
This appendix provides instance-level results for the experiments presented in the main paper.

\subsection{Baseline Comparisons}
Tables~\ref{tab:detailed_prodhon}--\ref{tab:detailed_schneider_500_600} report detailed results for Section~\ref{comparison-sota}.

\begin{sidewaystable}
\caption{Detailed results on the $\mathbb{P}$ benchmark of \citet{prins2004nouveaux}.
BKS values are taken from \citet{loffler2023conceptually} and \citet{he2025hybrid}.}
\label{tab:detailed_prodhon}
\label{tab:optimization_main_full_extended}
\scriptsize
\setlength{\tabcolsep}{2pt}
\begin{tabular}{@{}l r rr rr rr rr rrr@{}}
\toprule
& & \multicolumn{2}{c}{HCC-500K}
  & \multicolumn{2}{c}{TSBA$_\text{speed}$}
  & \multicolumn{2}{c}{GRASP/VNS}
  & \multicolumn{2}{c}{HALNS}
  & \multicolumn{3}{c}{NEO-DS} \\
\cmidrule(r){3-4}\cmidrule(r){5-6}\cmidrule(r){7-8}\cmidrule(r){9-10}\cmidrule(l){11-13}
Instance & BKS
& $E^{\text{gap}}_{\text{BKS}}$ (\%) & $T_{\text{total}}$ (s)
& $E^{\text{gap}}_{\text{BKS}}$ (\%) & $T_{\text{total}}$ (s)
& $E^{\text{gap}}_{\text{BKS}}$ (\%) & $T_{\text{total}}$ (s)
& $E^{\text{gap}}_{\text{BKS}}$ (\%) & $T_{\text{total}}$ (s)
& $E^{\text{gap}}_{\text{BKS}}$ (\%) & $T^{\text{LA}}$ (s) & $T_{\text{total}}$ (s) \\
\midrule
20-5-1a & 54,793 & 0.00 & 39 & 0.00 & 0.80 & 0.08 & 0.78 & 0.00 & 26 &  3.24  &  0.14  &  0.16 \\
20-5-1b & 39,104 & 0.00 & 54 & 0.00 & 0.53 & 0.00 & 0.67 & 0.00 & 21 &  3.65  &  0.14  &  0.16 \\
20-5-2a & 48,908 & 0.00 & 38 & 0.00 & 0.74 & 0.00 & 0.76 & 0.00 & 23 &  5.82  &  0.11  &  0.13 \\
20-5-2b & 37,542 & 0.00 & 67 & 0.00 & 0.51 & 0.00 & 0.65 & 0.00 & 23 &  5.49  &  0.13  &  0.16 \\
\midrule
\textbf{Average} &  & \textbf{0.00} & \textbf{49.50} & \textbf{0.00} & \textbf{0.65} & \textbf{0.02} & \textbf{0.71} & \textbf{0.00} & \textbf{23.15} & \textbf{ 4.55 } & \textbf{ 0.13 } & \textbf{ 0.15 } \\
\midrule
50-5-1a & 90,111 & 0.00 & 101 & 0.00 & 2.48 & 0.00 & 7.95 & 0.00 & 61 &  1.63  &  0.24  &  0.37 \\
50-5-1b & 63,242 & 0.00 & 65 & 0.00 & 2.35 & 0.00 & 8.59 & 0.00 & 54 &  1.18  &  0.25  &  0.46 \\
50-5-2a & 88,293 & 0.32 & 99 & 0.06 & 3.32 & 0.35 & 8.52 & 0.01 & 59 &  3.56  &  0.30  &  0.47 \\
50-5-2b & 67,308 & 0.21 & 200 & 0.14 & 3.07 & 0.54 & 9.18 & 0.00 & 59 &  5.80  &  0.28  &  0.41 \\
50-5-2bBIS & 51,822 & 0.03 & 98 & 0.08 & 2.70 & 0.02 & 8.98 & 0.00 & 53 &  3.18  &  0.27  &  0.34 \\
50-5-2BIS & 84,055 & 0.08 & 107 & 0.00 & 3.40 & 0.00 & 7.90 & 0.00 & 61 &  1.62  &  0.23  &  0.46 \\
50-5-3a & 86,203 & 0.07 & 101 & 0.19 & 3.34 & 0.19 & 7.78 & 0.00 & 67 &  2.37  &  0.23  &  0.51 \\
50-5-3b & 61,830 & 0.00 & 137 & 0.01 & 2.35 & 0.00 & 7.59 & 0.00 & 55 &  4.32  &  0.38  &  0.62 \\
\midrule
\textbf{Average} &  & \textbf{0.09} & \textbf{113.50} & \textbf{0.06} & \textbf{2.88} & \textbf{0.14} & \textbf{8.31} & \textbf{0.00} & \textbf{58.71} & \textbf{ 2.96 } & \textbf{ 0.27 } & \textbf{ 0.45 } \\
\midrule
100-5-1a & 274,814 & 0.56 & 520 & 0.37 & 15.14 & 0.44 & 70.15 & 0.11 & 385 &  1.49  &  1.22  &  2.39 \\
100-5-1b & 213,568 & 0.69 & 1190 & 0.50 & 11.68 & 0.38 & 70.81 & 0.01 & 280 &  1.25  &  1.44  &  2.47 \\
100-5-2a & 193,671 & 0.12 & 463 & 0.07 & 11.86 & 0.23 & 82.00 & 0.00 & 210 &  0.28  &  0.42  &  2.31 \\
100-5-2b & 157,095 & 0.04 & 859 & 0.05 & 8.11 & 0.07 & 61.93 & 0.00 & 186 &  0.45  &  0.46  &  2.11 \\
100-5-3a & 200,079 & 0.21 & 454 & 0.21 & 14.05 & 0.24 & 64.37 & 0.00 & 247 &  1.44  &  0.41  &  2.05 \\
100-5-3b & 152,441 & 0.30 & 684 & 0.03 & 8.39 & 1.03 & 57.29 & 0.00 & 139 &  3.09  &  0.41  &  1.74 \\
100-10-1a & 287,661 & 4.28 & 210 & 0.24 & 25.54 & 0.61 & 78.81 & 0.07 & 395 &  1.91  &  1.90  &  2.74 \\
100-10-1b & 230,989 & 4.03 & 188 & 0.47 & 16.57 & 1.19 & 87.95 & 0.00 & 223 &  2.35  &  1.95  &  2.62 \\
100-10-2a & 243,590 & 0.80 & 136 & 0.05 & 21.16 & 2.06 & 75.65 & 0.01 & 255 &  1.49  &  3.11  &  4.02 \\
100-10-2b & 203,988 & 0.25 & 261 & 0.00 & 10.93 & 1.23 & 67.50 & 0.00 & 147 &  2.54  &  2.83  &  3.59 \\
100-10-3a & 250,882 & 1.59 & 202 & 0.93 & 22.60 & 3.83 & 71.87 & 0.03 & 332 &  1.38  &  23.15  &  24.11 \\
100-10-3b & 203,114 & 1.51 & 224 & 0.29 & 14.88 & 5.53 & 79.76 & 0.14 & 207 &  5.37  &  39.31  &  40.17 \\
\midrule
\textbf{Average} &  & \textbf{1.20} & \textbf{449.25} & \textbf{0.27} & \textbf{15.08} & \textbf{1.40} & \textbf{72.34} & \textbf{0.03} & \textbf{250.53} & \textbf{ 1.92 } & \textbf{ 6.38 } & \textbf{ 7.53 } \\
\midrule
200-10-1a & 474,702 & 1.79 & 752 & 0.65 & 179.62 & 3.19 & 752.03 & 0.52 & 1479 &  1.00  &  6.07  &  11.50 \\
200-10-1b & 375,177 & 1.43 & 1346 & 0.42 & 115.72 & 2.74 & 735.75 & 0.09 & 1032 &  0.74  &  38.54  &  43.26 \\
200-10-2a & 448,077 & 0.82 & 1201 & 0.35 & 147.04 & 0.38 & 642.16 & 0.24 & 1283 &  0.63  &  2.59  &  7.26 \\
200-10-2b & 373,696 & 0.65 & 1349 & 0.14 & 69.52 & 0.23 & 683.19 & 0.00 & 951 &  0.55  &  3.71  &  7.74 \\
200-10-3a & 469,433 & 2.12 & 1251 & 0.49 & 176.25 & 0.48 & 661.82 & 0.37 & 1197 &  1.62  &  13.50  &  18.85 \\
200-10-3b & 362,320 & 2.01 & 1137 & 0.14 & 67.58 & 0.45 & 818.25 & 0.12 & 1092 &  1.54  &  3.09  &  7.80 \\
\midrule
\textbf{Average} &  & \textbf{1.47} & \textbf{1172.67} & \textbf{0.37} & \textbf{125.96} & \textbf{1.25} & \textbf{715.53} & \textbf{0.22} & \textbf{1172.16} & \textbf{ 1.01 } & \textbf{ 11.25 } & \textbf{ 16.07 } \\
\midrule
\multicolumn{2}{l}{Processor}
  & \multicolumn{2}{r}{Opteron 275}
  & \multicolumn{2}{r}{Xeon E5-2670}
  & \multicolumn{2}{r}{Xeon E5-2430v2}
  & \multicolumn{2}{r}{Ryzen 9 3900X}
  & \multicolumn{3}{r}{Xeon Gold 6248R} \\
\multicolumn{2}{l}{GHz}
  & \multicolumn{2}{r}{2.2}
  & \multicolumn{2}{r}{2.6}
  & \multicolumn{2}{r}{2.5}
  & \multicolumn{2}{r}{3.8}
  & \multicolumn{3}{r}{3.00}\\
\bottomrule
\end{tabular}
\end{sidewaystable}

\begin{sidewaystable}
\caption{Detailed results on the $\mathbb{T}$ benchmark of \citet{tuzun1999two}. BKS values are taken from \citet{sobhanan2024genetic}, \citet{he2025hybrid}, and \citet{loffler2023conceptually}.}
\label{tab:detailed_tuzun}
\label{tab:optimization_tuzun}
\scriptsize
\setlength{\tabcolsep}{4pt}  
\begin{tabular}{@{}l r rr rr rr rr rrr@{}}
\toprule
& & \multicolumn{2}{c}{HCC-500K}
  & \multicolumn{2}{c}{TSBA$_\text{speed}$}
  & \multicolumn{2}{c}{GRASP/VNS}
  & \multicolumn{2}{c}{HALNS}
  & \multicolumn{3}{c}{NEO-DS} \\
\cmidrule(r){3-4}\cmidrule(r){5-6}\cmidrule(r){7-8}\cmidrule(r){9-10}\cmidrule(l){11-13}
Instance & BKS
& $E^{\text{gap}}_{\text{BKS}}$ (\%) & $T_{\text{total}}$ (s)
& $E^{\text{gap}}_{\text{BKS}}$ (\%) & $T_{\text{total}}$ (s)
& $E^{\text{gap}}_{\text{BKS}}$ (\%) & $T_{\text{total}}$ (s)
& $E^{\text{gap}}_{\text{BKS}}$ (\%) & $T_{\text{total}}$ (s)
& $E^{\text{gap}}_{\text{BKS}}$ (\%) & $T^{\text{LA}}$ (s) & $T_{\text{total}}$ (s) \\
\midrule
111112 & 1467.68 & 0.54 & 275 & 0.07 & 23 & 2.29 & 65 & 0.00 & 163 &  4.51  &  3.27  &  4.27 \\
111122 & 1448.37 & 1.13 & 321 & 0.44 & 24 & 2.05 & 81 & 0.01 & 239 &  8.39  &  12.13  &  13.10 \\
111212 & 1394.80 & 0.41 & 244 & 0.12 & 22 & 0.64 & 59 & 0.00 & 163 &  3.26  &  3.68  &  4.48 \\
111222 & 1432.29 & 0.62 & 376 & 0.03 & 20 & 2.54 & 77 & 0.02 & 235 &  4.09  &  15.24  &  16.02 \\
112112 & 1167.16 & 0.50 & 489 & 0.07 & 21 & 0.44 & 69 & 0.00 & 184 &  4.67  &  3.37  &  5.63 \\
112122 & 1102.24 & 0.01 & 373 & 0.03 & 16 & 1.84 & 100 & 0.00 & 152 &  3.43  &  10.30  &  11.80 \\
112212 & 791.66 & 0.02 & 739 & 0.34 & 16 & 0.39 & 97 & 0.00 & 169 &  0.79  &  3.14  &  4.66 \\
112222 & 728.30 & 0.00 & 384 & 0.00 & 13 & 0.18 & 91 & 0.00 & 138 &  0.50  &  10.07  &  11.49 \\
113112 & 1238.24 & 0.17 & 357 & 0.10 & 20 & 1.14 & 74 & 0.02 & 197 &  3.03  &  3.00  &  4.07 \\
113122 & 1245.30 & 0.23 & 445 & 0.11 & 19 & 0.64 & 86 & 0.01 & 151 &  4.57  &  10.79  &  12.38 \\
113212 & 902.26 & 0.00 & 321 & 0.06 & 16 & 0.10 & 70 & 0.00 & 142 &  0.24  &  3.29  &  3.98 \\
113222 & 1018.29 & 0.03 & 386 & 0.47 & 19 & 0.13 & 89 & 0.18 & 155 &  7.24  &  11.45  &  12.21 \\
121112 & 2237.73 & 1.81 & 944 & 1.66 & 125 & 1.73 & 574 & 0.33 & 857 &  2.74  &  8.77  &  14.62 \\
121122 & 2137.45 & 2.58 & 847 & 1.74 & 107 & 5.35 & 681 & 0.13 & 1011 &  4.11  &  23.45  &  26.78 \\
121212 & 2195.17 & 2.40 & 907 & 1.25 & 145 & 1.04 & 558 & 0.20 & 859 &  5.08  &  9.24  &  13.49 \\
121222 & 2214.86 & 2.18 & 860 & 1.08 & 131 & 3.19 & 662 & 0.06 & 1015 &  3.44  &  58.21  &  63.81 \\
122112 & 2070.43 & 1.13 & 1606 & 0.65 & 117 & 0.62 & 704 & 0.08 & 862 &  2.19  &  11.72  &  20.49 \\
122122 & 1685.52 & 2.76 & 941 & 0.60 & 128 & 7.56 & 927 & 0.01 & 1029 &  4.23  &  36.81  &  41.30 \\
122212 & 1449.62 & 0.86 & 1861 & 0.45 & 128 & 1.58 & 697 & 0.01 & 1068 &  2.39  &  8.88  &  17.36 \\
122222 & 1082.46 & 0.33 & 812 & 0.33 & 101 & 0.52 & 656 & 0.01 & 842 &  1.81  &  25.75  &  29.76 \\
123112 & 1942.23 & 1.48 & 968 & 1.24 & 120 & 5.25 & 595 & 0.55 & 1047 &  1.81  &  6.94  &  10.19 \\
123122 & 1910.08 & 2.21 & 740 & 1.08 & 157 & 3.26 & 770 & 0.16 & 1005 &  3.87  &  26.42  &  29.80 \\
123212 & 1760.20 & 0.23 & 2055 & 0.38 & 103 & 0.78 & 690 & 0.05 & 1119 &  5.34  &  8.48  &  12.41 \\
123222 & 1390.74 & 0.33 & 1038 & 0.16 & 91 & 4.01 & 704 & 0.00 & 663 &  0.98  &  29.73  &  31.48 \\
131112 & 1866.75 & 3.90 & 504 & 2.83 & 56 & 3.87 & 230 & 1.45 & 482 &  4.22  &  7.11  &  10.25 \\
131122 & 1819.68 & 2.07 & 635 & 0.59 & 61 & 2.51 & 259 & 0.17 & 473 &  3.40  &  38.26  &  39.74 \\
131212 & 1960.02 & 2.52 & 664 & 0.55 & 49 & 0.58 & 233 & 0.27 & 407 &  6.08  &  7.37  &  9.32 \\
131222 & 1792.77 & 2.55 & 485 & 0.76 & 63 & 0.93 & 273 & 0.14 & 378 &  2.73  &  35.83  &  37.60 \\
132112 & 1443.32 & 0.40 & 1049 & 0.12 & 72 & 1.15 & 260 & 0.00 & 369 &  8.48  &  6.75  &  10.98 \\
132122 & 1429.30 & 1.23 & 805 & 1.31 & 72 & 1.16 & 344 & 0.68 & 442 &  1.78  &  26.59  &  31.09 \\
132212 & 1204.42 & 0.12 & 2197 & 0.26 & 51 & 1.14 & 256 & 0.00 & 428 &  0.91  &  11.80  &  13.98 \\
132222 & 924.68 & 0.91 & 982 & 0.26 & 55 & 1.01 & 318 & 0.00 & 484 &  1.63  &  17.47  &  19.20 \\
133112 & 1694.18 & 0.37 & 1046 & 0.83 & 50 & 1.05 & 281 & 0.26 & 530 &  4.29  &  5.73  &  10.15 \\
133122 & 1392.00 & 0.83 & 925 & 1.34 & 56 & 1.85 & 309 & 0.13 & 474 &  4.93  &  30.22  &  32.46 \\
133212 & 1197.95 & 0.11 & 1375 & 0.14 & 42 & 0.37 & 276 & 0.00 & 479 &  4.07  &  5.42  &  9.54 \\
133222 & 1151.37 & 0.26 & 911 & 0.34 & 58 & 0.57 & 343 & 0.04 & 385 &  2.03  &  20.34  &  22.65 \\
\midrule
\textbf{Average} & & \textbf{1.03} & \textbf{830} & \textbf{0.61} & \textbf{66} & \textbf{1.76} & \textbf{349} & \textbf{0.14} & \textbf{522} & \textbf{3.53} & \textbf{15.47} & \textbf{18.40} \\
\midrule
\multicolumn{2}{l}{Processor}
  & \multicolumn{2}{r}{Opteron 275}
  & \multicolumn{2}{r}{Xeon E5-2670}
  & \multicolumn{2}{r}{Xeon E5-2430v2}
  & \multicolumn{2}{r}{Ryzen 9 3900X}
  & \multicolumn{3}{r}{Xeon Gold 6248R} \\
\multicolumn{2}{l}{GHz}
  & \multicolumn{2}{r}{2.2}
  & \multicolumn{2}{r}{2.6}
  & \multicolumn{2}{r}{2.5}
  & \multicolumn{2}{r}{3.8}
  & \multicolumn{3}{r}{3.00}\\
\bottomrule
\end{tabular}
\end{sidewaystable}

\begin{sidewaystable}
\caption{Detailed results on the $\mathbb{B}$ benchmark of \citet{barreto_2004}. BKS values are taken from \citet{sobhanan2024genetic}, \citet{he2025hybrid}, and \citet{loffler2023conceptually}.}
\label{tab:detailed_barreto}
\scriptsize
\setlength{\tabcolsep}{4pt}  
\begin{tabular}{@{}l r rr rr rr rr rrr@{}}
\toprule
& & \multicolumn{2}{c}{HCC-500K}
  & \multicolumn{2}{c}{TSBA$_\text{speed}$}
  & \multicolumn{2}{c}{GRASP/VNS}
  & \multicolumn{2}{c}{HALNS}
  & \multicolumn{3}{c}{NEO-DS} \\
\cmidrule(r){3-4}\cmidrule(r){5-6}\cmidrule(r){7-8}\cmidrule(r){9-10}\cmidrule(l){11-13}
Instance & BKS
& $E^{\text{gap}}_{\text{BKS}}$ (\%) & $T_{\text{total}}$ (s)
& $E^{\text{gap}}_{\text{BKS}}$ (\%) & $T_{\text{total}}$ (s)
& $E^{\text{gap}}_{\text{BKS}}$ (\%) & $T_{\text{total}}$ (s)
& $E^{\text{gap}}_{\text{BKS}}$ (\%) & $T_{\text{total}}$ (s)
& $E^{\text{gap}}_{\text{BKS}}$ (\%) & $T^{\text{LA}}$ (s) & $T_{\text{total}}$ (s) \\
\midrule
Christofides69-50x5 & 565.6 & 0.00 & 73 & 0.00 & 2 & 2.55 & 7 & 0.00 & 55 &  7.06  &  0.15  &  0.34 \\
Christofides69-75x10 & 844.4 & 1.24 & 207 & 0.72 & 12 & 0.92 & 27 & 0.53 & 121 &  7.07  &  1.64  &  1.99 \\
Christofides69-100x10 & 833.4 & 0.24 & 403 & 0.14 & 26 & 1.31 & 64 & 0.07 & 197 &  11.88  &  2.90  &  4.00 \\
Daskin95-88x8 & 355.8 & 0.01 & 250 & 0.03 & 6 & 0.90 & 46 & -0.01 & 92 &  20.39  &  1.31  &  1.87 \\
Daskin95-150x10 & 43919.9 & 1.31 & 613 & 0.38 & 52 & 0.82 & 345 & 0.00 & 411 &  9.97  &  22.99  &  24.72 \\
Gaskell67-21x5 & 424.9 & 0.00 & 25 & 0.00 & 1 & 0.00 & 1 & 0.00 & 25 &  3.76  &  0.11  &  0.13 \\
Gaskell67-22x5 & 585.1 & 0.00 & 21 & 0.00 & 0 & 0.00 & 1 & 0.00 & 24 &  0.00  &  0.18  &  0.27 \\
Gaskell67-29x5 & 512.1 & 0.00 & 40 & 0.00 & 1 & 0.00 & 2 & 0.00 & 29 &  12.31  &  0.23  &  0.28 \\
Gaskell67-32x5-1 & 562.2 & 0.00 & 58 & 0.00 & 1 & 0.00 & 2 & 0.00 & 34 &  0.39  &  0.10  &  0.36 \\
Gaskell67-32x5-2 & 504.3 & 0.01 & 55 & 0.00 & 1 & 0.00 & 2 & 0.01 & 33 &  0.05  &  0.10  &  0.30 \\
Gaskell67-36x5 & 460.4 & 0.01 & 61 & 0.00 & 1 & 0.00 & 2 & -0.01 & 36 &  -0.01  &  0.15  &  0.31 \\
Min92-27x5 & 3062.0 & 0.00 & 38 & 0.00 & 1 & 0.00 & 2 & 0.00 & 29 &  23.23  &  0.08  &  0.12 \\
Min92-134x8 & 5709.0 & 0.41 & 460 & 0.23 & 29 & 1.11 & 166 & 0.00 & 240 &  28.98  &  1.63  &  2.54 \\
\midrule
\textbf{Average} & & \textbf{0.25} & \textbf{177} & \textbf{0.12} & \textbf{10} & \textbf{0.59} & \textbf{51} & \textbf{0.05} & \textbf{102} & \textbf{9.62} & \textbf{2.43} & \textbf{2.86} \\
\midrule
\multicolumn{2}{l}{Processor}
  & \multicolumn{2}{r}{Opteron 275}
  & \multicolumn{2}{r}{Xeon E5-2670}
  & \multicolumn{2}{r}{Xeon E5-2430v2}
  & \multicolumn{2}{r}{Ryzen 9 3900X}
  & \multicolumn{3}{r}{Xeon Gold 6248R} \\
\multicolumn{2}{l}{GHz}
  & \multicolumn{2}{r}{2.2}
  & \multicolumn{2}{r}{2.6}
  & \multicolumn{2}{r}{2.5}
  & \multicolumn{2}{r}{3.8}
  & \multicolumn{3}{r}{3.00}\\
\bottomrule
\end{tabular}
\end{sidewaystable}

\clearpage
\begingroup
\tiny
\captionsetup{font=normalsize}
\setlength{\tabcolsep}{4pt}
\renewcommand{\arraystretch}{1.05}
\sisetup{
  group-separator = {},
  group-minimum-digits = 4,
  table-number-alignment = center,
  detect-weight,
  mode = match,
}
\begin{longtable}{
    l
    l
    S[table-format=2.2]
    S[table-format=5.0]
    S[table-format=2.2]
    S[table-format=5.0]
    S[table-format=2.2]
    S[table-format=4.2]
    S[table-format=4.2]
}
\caption{Detailed results on the $\mathbb{S}$ benchmark of \citet{schneider2019large} (100--200 customers). BKS values are taken from \citet{he2025hybrid}.}
\label{tab:detailed_schneider_100_200} \\
\toprule
\multicolumn{2}{c}{} &
\multicolumn{2}{c}{\text{TSBA$_{\text{basic}}$}} &
\multicolumn{2}{c}{\text{HGAMP}} &
\multicolumn{3}{c}{\text{NEO-DS}} \\
\cmidrule(lr){3-4} \cmidrule(lr){5-6} \cmidrule(lr){7-9}
Instance & BKS &
\text{$E^{\text{gap}}_{BKS}$ (\%)} & \text{$T_{\text{total}}$ (s)} &
\text{$E^{\text{gap}}_{BKS}$ (\%)} & \text{$T_{\text{total}}$ (s)} &
\text{$E^{\text{gap}}_{BKS}$ (\%)} & \text{$T^{\text{LA}}$ (s)} & \text{$T_{\text{total}}$ (s)} \\
\midrule
\endfirsthead
\toprule
\multicolumn{2}{c}{} &
\multicolumn{2}{c}{\text{TSBA$_{\text{basic}}$}} &
\multicolumn{2}{c}{\text{HGAMP}} &
\multicolumn{3}{c}{\text{NEO-DS}} \\
\cmidrule(lr){3-4} \cmidrule(lr){5-6} \cmidrule(lr){7-9}
Instance & BKS &
\text{$E^{\text{gap}}_{BKS}$ (\%)} & \text{$T_{\text{total}}$ (s)} &
\text{$E^{\text{gap}}_{BKS}$ (\%)} & \text{$T_{\text{total}}$ (s)} &
\text{$E^{\text{gap}}_{BKS}$ (\%)} & \text{$T^{\text{LA}}$ (s)} & \text{$T_{\text{total}}$ (s)} \\
\midrule
\endhead
\midrule
\multicolumn{9}{r}{\textit{Continued on next page}} \\
\endfoot
\endlastfoot
100-5-1c & 134516 & 0.07 & 81 & 0.00 & 1070 & 3.44 & 0.25 & 0.84 \\
100-5-1d & 275749 & 0.02 & 53 & 0.00 & 950 & 1.88 & 1.17 & 2.04 \\
100-5-1e & 292301 & 0.03 & 69 & 0.32 & 1104 & 0.91 & 0.27 & 1.82 \\
100-5-2c & 83855 & 0.45 & 67 & 0.15 & 1070 & 4.70 & 0.21 & 0.97 \\
100-5-2d & 242105 & 0.13 & 52 & 0.00 & 1099 & 1.15 & 1.17 & 2.13 \\
100-5-2e & 253888 & 0.07 & 66 & 0.05 & 2171 & 0.34 & 0.24 & 4.42 \\
100-5-3c & 87555 & 0.06 & 55 & 0.00 & 949 & 5.54 & 0.17 & 0.68 \\
100-5-3d & 226634 & 0.09 & 49 & 0.06 & 892 & 1.02 & 0.61 & 1.39 \\
100-5-3e & 252603 & 0.02 & 79 & 0.01 & 1436 & 0.54 & 0.27 & 4.33 \\
100-5-4a & 255853 & 0.02 & 52 & 0.00 & 1087 & 1.62 & 0.28 & 1.17 \\
100-5-4b & 214425 & 0.00 & 30 & 0.00 & 1068 & 2.50 & 0.30 & 1.12 \\
100-5-4c & 98104 & 0.09 & 53 & 0.01 & 1091 & 5.10 & 0.32 & 1.22 \\
100-5-4d & 250301 & 0.19 & 46 & 0.28 & 1062 & 1.89 & 1.02 & 2.00 \\
100-5-4e & 211113 & 0.05 & 58 & 0.07 & 1651 & 1.28 & 11.51 & 13.28 \\
100-10-1c & 92629 & 0.06 & 98 & 0.00 & 1285 & 8.28 & 0.59 & 0.70 \\
100-10-1d & 363930 & 0.07 & 61 & 0.00 & 1132 & 0.70 & 1.53 & 1.85 \\
100-10-1e & 344322 & 0.08 & 60 & 0.15 & 1114 & 1.18 & 1.33 & 2.30 \\
100-10-2c & 84717 & 0.03 & 99 & 0.00 & 1169 & 16.69 & 1.28 & 1.54 \\
100-10-2d & 343252 & 0.00 & 56 & 0.00 & 1255 & 1.39 & 1.33 & 1.70 \\
100-10-2e & 332900 & 0.08 & 61 & 0.26 & 1000 & 1.44 & 1.28 & 2.14 \\
100-10-3c & 85369 & 0.40 & 97 & 0.00 & 1129 & 10.18 & 1.01 & 1.16 \\
100-10-3d & 329990 & 0.01 & 52 & 0.01 & 1161 & 0.67 & 1.08 & 1.39 \\
100-10-3e & 318109 & 0.05 & 64 & 0.04 & 1185 & 0.63 & 1.34 & 2.15 \\
100-10-4a & 253471 & 0.29 & 48 & 0.01 & 1333 & 0.72 & 19.10 & 20.01 \\
100-10-4b & 211354 & 0.00 & 39 & 0.00 & 1340 & 1.88 & 25.87 & 26.60 \\
100-10-4c & 86215 & 0.00 & 74 & 0.00 & 1320 & 6.90 & 1.03 & 1.44 \\
100-10-4d & 328181 & 0.05 & 63 & 0.07 & 1381 & 0.83 & 1.30 & 1.65 \\
100-10-4e & 308757 & 0.18 & 72 & 0.36 & 1320 & 1.83 & 0.76 & 1.58 \\
200-10-1c & 156029 & 0.48 & 590 & 0.09 & 3218 & 5.17 & 1.80 & 2.80 \\
200-10-1d & 638068 & 0.32 & 418 & 0.19 & 2590 & 0.78 & 3.75 & 5.87 \\
200-10-1e & 599069 & 0.11 & 464 & 0.35 & 2521 & 1.57 & 12.26 & 17.47 \\
200-10-2c & 144046 & 0.23 & 476 & 0.21 & 2756 & 2.53 & 2.64 & 4.90 \\
200-10-2d & 663154 & 0.16 & 343 & 0.14 & 2801 & 0.82 & 1.60 & 3.60 \\
200-10-2e & 618858 & 0.05 & 318 & 0.11 & 2834 & 0.66 & 1.59 & 6.29 \\
200-10-3c & 184783 & 0.61 & 588 & 0.66 & 2651 & 5.90 & 19.42 & 21.88 \\
200-10-3d & 640289 & 0.02 & 385 & 0.18 & 2542 & 2.25 & 1.74 & 4.04 \\
200-10-3e & 604480 & 0.13 & 388 & 0.33 & 2509 & 1.55 & 2.38 & 7.60 \\
200-10-4a & 452430 & 0.15 & 316 & 0.13 & 2997 & 1.18 & 89.60 & 94.31 \\
200-10-4b & 369580 & 0.18 & 215 & 0.10 & 3254 & 0.89 & 17.89 & 22.35 \\
200-10-4c & 144013 & 0.41 & 611 & 0.34 & 3323 & 6.58 & 1.79 & 3.09 \\
200-10-4d & 617932 & 0.18 & 370 & 0.15 & 3105 & 0.62 & 2.97 & 4.94 \\
200-10-4e & 562843 & 0.10 & 403 & 0.19 & 2959 & 1.00 & 2.81 & 8.13 \\
200-15-1a & 460430 & 0.29 & 589 & 0.33 & 2728 & 1.06 & 85.45 & 89.84 \\
200-15-1b & 366359 & 0.38 & 274 & 0.18 & 2781 & 1.75 & 64.98 & 68.68 \\
200-15-1c & 148141 & 0.71 & 692 & 0.67 & 3693 & 7.95 & 3.51 & 4.41 \\
200-15-1d & 813576 & 0.12 & 593 & 0.16 & 2755 & 1.17 & 73.21 & 74.38 \\
200-15-1e & 708585 & 1.58 & 624 & 0.16 & 2834 & 0.58 & 3.50 & 5.56 \\
200-15-2a & 513512 & 0.11 & 627 & 0.15 & 3046 & 0.79 & 164.56 & 169.14 \\
200-15-2b & 406685 & 0.11 & 271 & 0.06 & 3046 & 1.64 & 216.64 & 221.26 \\
200-15-2c & 134779 & 0.54 & 466 & 0.42 & 3599 & 7.06 & 4.51 & 5.41 \\
200-15-2d & 811361 & 0.14 & 544 & 0.22 & 2854 & 1.53 & 24.91 & 25.93 \\
200-15-2e & 712524 & 0.09 & 569 & 0.29 & 3326 & 2.04 & 268.89 & 271.13 \\
200-15-3a & 455351 & 0.16 & 558 & 0.07 & 2871 & 1.26 & 101.05 & 105.69 \\
200-15-3b & 356887 & 0.26 & 218 & 0.11 & 2933 & 0.98 & 61.17 & 65.45 \\
200-15-3c & 140765 & 0.56 & 536 & 0.24 & 3741 & 8.68 & 4.33 & 4.94 \\
200-15-3d & 877543 & 0.10 & 502 & 0.10 & 2551 & 1.44 & 4.34 & 5.39 \\
200-15-3e & 816001 & 0.12 & 524 & 0.05 & 2723 & 0.55 & 30.98 & 32.69 \\
200-15-4a & 432672 & 0.23 & 580 & 0.15 & 3260 & 0.66 & 520.94 & 525.94 \\
200-15-4b & 349088 & 0.20 & 274 & 0.14 & 3224 & 0.96 & 298.16 & 302.90 \\
200-15-4c & 143052 & 0.71 & 558 & 0.48 & 3807 & 9.51 & 4.62 & 5.69 \\
200-15-4d & 826829 & 0.18 & 547 & 0.14 & 2947 & 1.64 & 38.12 & 39.31 \\
200-15-4e & 700013 & 0.15 & 744 & 0.24 & 3050 & 1.14 & 5.01 & 8.07 \\
\midrule
\bfseries Average & \multicolumn{1}{c}{\bfseries } & \bfseries 0.20 & \bfseries 289 & \bfseries 0.15 & \bfseries 2188 & \bfseries 2.73 & \bfseries 35.77 & \bfseries 37.79 \\
\midrule
\multicolumn{2}{l}{Processor}
  & \multicolumn{2}{r}{Xeon E5-2670}
  & \multicolumn{2}{r}{Xeon E-2670}
  & \multicolumn{3}{r}{Xeon Gold 6248R} \\
\multicolumn{2}{l}{GHz}
  & \multicolumn{2}{r}{2.6}
  & \multicolumn{2}{r}{2.5}
  & \multicolumn{3}{r}{3.0}\\
\bottomrule
\end{longtable}
\endgroup

\clearpage
\begingroup
\tiny
\captionsetup{font=normalsize}
\setlength{\tabcolsep}{4pt}
\renewcommand{\arraystretch}{1.05}
\sisetup{
  group-separator = {},
  group-minimum-digits = 4,
  table-number-alignment = center,
  detect-weight,
  mode = match,
}
\begin{longtable}{
    l
    l
    S[table-format=2.2]
    S[table-format=5.0]
    S[table-format=2.2]
    S[table-format=5.0]
    S[table-format=2.2]
    S[table-format=4.2]
    S[table-format=4.2]
}
\caption{Detailed results on the $\mathbb{S}$ benchmark of \citet{schneider2019large} (300--400 customers). BKS values are taken from \citet{he2025hybrid}.}
\label{tab:detailed_schneider_300_400} \\
\toprule
\multicolumn{2}{c}{} &
\multicolumn{2}{c}{\text{TSBA$_{\text{basic}}$}} &
\multicolumn{2}{c}{\text{HGAMP}} &
\multicolumn{3}{c}{\text{NEO-DS}} \\
\cmidrule(lr){3-4} \cmidrule(lr){5-6} \cmidrule(lr){7-9}
Instance & BKS &
\text{$E^{\text{gap}}_{BKS}$ (\%)} & \text{$T_{\text{total}}$ (s)} &
\text{$E^{\text{gap}}_{BKS}$ (\%)} & \text{$T_{\text{total}}$ (s)} &
\text{$E^{\text{gap}}_{BKS}$ (\%)} & \text{$T^{\text{LA}}$ (s)} & \text{$T_{\text{total}}$ (s)} \\
\midrule
\endfirsthead
\toprule
\multicolumn{2}{c}{} &
\multicolumn{2}{c}{\text{TSBA$_{\text{basic}}$}} &
\multicolumn{2}{c}{\text{HGAMP}} &
\multicolumn{3}{c}{\text{NEO-DS}} \\
\cmidrule(lr){3-4} \cmidrule(lr){5-6} \cmidrule(lr){7-9}
Instance & BKS &
\text{$E^{\text{gap}}_{BKS}$ (\%)} & \text{$T_{\text{total}}$ (s)} &
\text{$E^{\text{gap}}_{BKS}$ (\%)} & \text{$T_{\text{total}}$ (s)} &
\text{$E^{\text{gap}}_{BKS}$ (\%)} & \text{$T^{\text{LA}}$ (s)} & \text{$T_{\text{total}}$ (s)} \\
\midrule
\endhead
\midrule
\multicolumn{9}{r}{\textit{Continued on next page}} \\
\endfoot
\endlastfoot
300-15-1a & 854503 & 0.31 & 2085 & 0.38 & 4776 & 2.41 & 46.18 & 57.72 \\
300-15-1b & 621894 & 0.39 & 1387 & 0.23 & 5066 & 1.30 & 158.56 & 169.48 \\
300-15-1c & 364979 & 0.86 & 1613 & 0.18 & 6459 & 5.51 & 5.35 & 7.62 \\
300-15-1d & 1337930 & 0.14 & 1557 & 0.09 & 4300 & 1.17 & 145.99 & 149.15 \\
300-15-1e & 1217690 & 0.39 & 1949 & 0.25 & 4740 & 1.82 & 109.47 & 115.51 \\
300-15-2a & 757931 & 0.31 & 2135 & 0.38 & 5265 & 1.43 & 17.68 & 29.57 \\
300-15-2b & 556948 & 0.40 & 1250 & 0.12 & 4935 & 0.89 & 13.69 & 24.33 \\
300-15-2c & 310061 & 0.83 & 1465 & 0.31 & 6299 & 3.66 & 6.74 & 9.63 \\
300-15-2d & 1301210 & 0.17 & 1579 & 0.12 & 4937 & 1.47 & 77.02 & 79.99 \\
300-15-2e & 1272700 & 0.26 & 1853 & 0.33 & 4915 & 1.73 & 67.18 & 72.39 \\
300-15-3a & 776531 & 0.24 & 1897 & 0.21 & 5246 & 0.96 & 11.01 & 21.18 \\
300-15-3b & 593743 & 0.21 & 1065 & 0.04 & 4857 & 0.43 & 158.94 & 169.42 \\
300-15-3c & 340155 & 0.52 & 1434 & 0.60 & 6186 & 4.01 & 6.13 & 7.99 \\
300-15-3d & 1355890 & 0.30 & 1495 & 0.15 & 4595 & 0.96 & 151.24 & 154.25 \\
300-15-3e & 1286877 & 0.13 & 2209 & 0.20 & 4778 & 1.45 & 7.08 & 15.70 \\
300-15-4a & 746407 & 0.35 & 1976 & 0.46 & 5757 & 1.18 & 66.49 & 79.23 \\
300-15-4b & 558831 & 0.26 & 1291 & 0.08 & 5236 & 1.23 & 171.18 & 182.46 \\
300-15-4c & 302390 & 0.74 & 1843 & 0.42 & 6658 & 3.59 & 6.25 & 10.91 \\
300-15-4d & 1285714 & 0.32 & 1779 & 0.25 & 4785 & 1.42 & 178.40 & 181.73 \\
300-15-4e & 1172200 & 0.26 & 1939 & 0.14 & 5074 & 0.89 & 69.23 & 77.72 \\
300-20-1a & 944798 & 7.04 & 1613 & 0.34 & 5098 & 1.09 & 15.89 & 23.49 \\
300-20-1b & 739410 & 2.14 & 1409 & 0.34 & 5252 & 0.49 & 15.20 & 22.59 \\
300-20-1c & 361735 & 0.85 & 1554 & 0.29 & 6816 & 8.83 & 9.14 & 10.14 \\
300-20-1d & 1572968 & 0.26 & 1625 & 0.38 & 5239 & 1.95 & 10.17 & 12.14 \\
300-20-1e & 1316880 & 5.78 & 1602 & 0.31 & 4772 & 1.77 & 37.92 & 46.46 \\
300-20-2a & 909306 & 0.16 & 2044 & 0.37 & 5634 & 1.20 & 82.57 & 90.36 \\
300-20-2b & 694525 & 0.26 & 1437 & 0.20 & 5731 & 0.96 & 15.97 & 23.22 \\
300-20-2c & 298522 & 0.64 & 1564 & 0.34 & 6632 & 4.45 & 10.19 & 11.97 \\
300-20-2d & 1569042 & 0.10 & 1900 & 0.22 & 5864 & 1.08 & 8.89 & 10.88 \\
300-20-2e & 1282040 & 8.29 & 2240 & 0.21 & 4994 & 0.85 & 8.14 & 20.30 \\
300-20-3a & 927452 & 0.36 & 1779 & 0.31 & 4924 & 1.13 & 439.32 & 446.66 \\
300-20-3b & 750633 & 0.17 & 1425 & 0.23 & 4980 & 0.62 & 279.70 & 286.61 \\
300-20-3c & 304269 & 0.96 & 1177 & 0.42 & 6277 & 7.29 & 10.00 & 12.19 \\
300-20-3d & 1538930 & 0.14 & 1657 & 0.15 & 4736 & 1.16 & 58.72 & 60.69 \\
300-20-3e & 1261660 & 2.30 & 2091 & 0.16 & 4501 & 1.55 & 129.68 & 138.59 \\
300-20-4a & 858649 & 0.17 & 1993 & 0.14 & 5835 & 1.16 & 48.29 & 56.73 \\
300-20-4b & 687179 & 0.33 & 1333 & 0.16 & 5952 & 1.69 & 211.05 & 219.11 \\
300-20-4c & 299157 & 0.66 & 1678 & 0.13 & 7350 & 5.71 & 8.68 & 11.35 \\
300-20-4d & 1539900 & 0.17 & 1751 & 0.18 & 5911 & 0.50 & 50.97 & 53.01 \\
300-20-4e & 1322150 & 1.76 & 2624 & 0.20 & 5384 & 1.35 & 105.33 & 114.62 \\
400-20-1a & 1139820 & 0.27 & 5460 & 0.27 & 7168 & 0.88 & 486.13 & 501.40 \\
400-20-1b & 876157 & 0.57 & 2844 & 0.29 & 7833 & 1.06 & 657.95 & 672.88 \\
400-20-1c & 466857 & 0.73 & 3935 & 0.28 & 10503 & 6.52 & 12.45 & 16.10 \\
400-20-1d & 1956824 & 0.20 & 3780 & 0.49 & 8334 & 1.84 & 16.03 & 20.45 \\
400-20-1e & 1641880 & 6.64 & 4669 & 0.20 & 6697 & 1.93 & 115.67 & 131.14 \\
400-20-2a & 1053445 & 0.15 & 5850 & 0.33 & 8246 & 1.44 & 109.96 & 124.69 \\
400-20-2b & 827077 & 0.36 & 2957 & 0.15 & 7771 & 0.56 & 316.33 & 329.84 \\
400-20-2c & 394544 & 0.39 & 3976 & 0.32 & 10390 & 3.98 & 13.53 & 18.16 \\
400-20-2d & 1875072 & 0.14 & 4029 & 0.14 & 8307 & 1.02 & 15.78 & 20.18 \\
400-20-2e & 1560160 & 3.21 & 6085 & 0.13 & 6865 & 1.28 & 62.15 & 79.74 \\
400-20-3a & 1098670 & 0.20 & 4705 & 0.13 & 7985 & 0.63 & 27.59 & 44.85 \\
400-20-3b & 847311 & 0.40 & 2731 & 0.22 & 7438 & 0.70 & 574.42 & 589.43 \\
400-20-3c & 390047 & 0.82 & 2898 & 0.29 & 9916 & 3.41 & 13.54 & 16.18 \\
400-20-3d & 1927290 & 0.19 & 3356 & 0.12 & 7306 & 1.12 & 14.31 & 18.33 \\
400-20-3e & 1656860 & 7.45 & 5252 & 0.22 & 6896 & 0.95 & 37.23 & 54.28 \\
400-20-4a & 1080890 & 0.23 & 5802 & 0.21 & 8612 & 1.24 & 153.09 & 169.26 \\
400-20-4b & 841464 & 0.44 & 3261 & 0.21 & 9067 & 1.71 & 338.75 & 355.14 \\
400-20-4c & 350273 & 0.55 & 2671 & 0.37 & 10777 & 5.10 & 16.98 & 23.80 \\
400-20-4d & 1834809 & 0.09 & 4318 & 0.18 & 9223 & 0.75 & 63.99 & 68.61 \\
400-20-4e & 1558411 & 4.07 & 4993 & 0.15 & 8343 & 1.39 & 10.36 & 21.20 \\
400-25-1a & 1153960 & 0.24 & 5401 & 0.30 & 7789 & 1.65 & 125.43 & 139.05 \\
400-25-1b & 889193 & 0.34 & 2878 & 0.15 & 7763 & 0.98 & 440.19 & 454.93 \\
400-25-1c & 393581 & 0.78 & 3579 & 0.29 & 12532 & 5.88 & 21.44 & 22.74 \\
400-25-1d & 2336630 & 0.31 & 3282 & 0.05 & 7791 & 0.66 & 680.21 & 683.17 \\
400-25-1e & 1863200 & 10.48 & 4181 & 0.11 & 6579 & 0.42 & 102.58 & 116.28 \\
400-25-2a & 1090120 & 0.19 & 7027 & 0.15 & 8625 & 0.45 & 893.59 & 909.58 \\
400-25-2b & 867716 & 0.26 & 3608 & 0.18 & 8011 & 1.65 & 991.08 & 1003.88 \\
400-25-2c & 359606 & 0.61 & 3497 & 0.24 & 11149 & 4.85 & 23.10 & 26.22 \\
400-25-2d & 2351903 & 0.07 & 3710 & 0.17 & 7520 & 0.99 & 834.05 & 837.48 \\
400-25-2e & 1902730 & 4.26 & 4965 & 0.13 & 7902 & 1.09 & 274.24 & 283.33 \\
400-25-3a & 1105783 & 0.11 & 5606 & 0.30 & 8249 & 0.65 & 97.79 & 111.95 \\
400-25-3b & 861981 & 0.42 & 3396 & 0.22 & 8212 & 1.09 & 669.06 & 682.42 \\
400-25-3c & 392264 & 0.58 & 2998 & 0.22 & 9778 & 4.74 & 24.44 & 26.46 \\
400-25-3d & 2321358 & 0.12 & 4212 & 0.26 & 7728 & 1.00 & 981.10 & 983.90 \\
400-25-3e & 1890950 & 3.16 & 4099 & 0.21 & 6966 & 0.80 & 220.49 & 229.85 \\
400-25-4a & 1014350 & 0.24 & 5540 & 0.25 & 8496 & 0.97 & 587.29 & 603.37 \\
400-25-4b & 801135 & 0.25 & 3118 & 0.13 & 9085 & 1.13 & 221.78 & 238.22 \\
400-25-4c & 378419 & 0.89 & 3195 & 0.16 & 11428 & 5.16 & 20.92 & 27.41 \\
400-25-4d & 2360020 & 0.21 & 4334 & 0.13 & 8912 & 0.80 & 1067.81 & 1070.96 \\
400-25-4e & 1935600 & 2.99 & 4204 & 0.18 & 7785 & 0.69 & 21.77 & 31.28 \\
\midrule
\bfseries Average & \multicolumn{1}{c}{\bfseries } & \bfseries 1.17 & \bfseries 2934 & \bfseries 0.23 & \bfseries 6959 & \bfseries 1.94 & \bfseries 179.80 & \bfseries 188.07 \\
\midrule
\multicolumn{2}{l}{Processor}
  & \multicolumn{2}{r}{Xeon E5-2670}
  & \multicolumn{2}{r}{Xeon E-2670}
  & \multicolumn{3}{r}{Xeon Gold 6248R} \\
\multicolumn{2}{l}{GHz}
  & \multicolumn{2}{r}{2.6}
  & \multicolumn{2}{r}{2.5}
  & \multicolumn{3}{r}{3.0}\\
\bottomrule
\end{longtable}
\endgroup

\clearpage
\begingroup
\tiny
\captionsetup{font=normalsize}
\setlength{\tabcolsep}{4pt}
\renewcommand{\arraystretch}{1.05}
\sisetup{
  group-separator = {},
  group-minimum-digits = 4,
  table-number-alignment = center,
  detect-weight,
  mode = match,
}
\begin{longtable}{
    l
    l
    S[table-format=2.2]
    S[table-format=5.0]
    S[table-format=2.2]
    S[table-format=5.0]
    S[table-format=2.2]
    S[table-format=4.2]
    S[table-format=4.2]
}
\caption{Detailed results on the $\mathbb{S}$ benchmark of \citet{schneider2019large} (500--600 customers). BKS values are taken from \citet{he2025hybrid}.}
\label{tab:detailed_schneider_500_600} \\
\toprule
\multicolumn{2}{c}{} &
\multicolumn{2}{c}{\text{TSBA$_{\text{basic}}$}} &
\multicolumn{2}{c}{\text{HGAMP}} &
\multicolumn{3}{c}{\text{NEO-DS}} \\
\cmidrule(lr){3-4} \cmidrule(lr){5-6} \cmidrule(lr){7-9}
Instance & BKS &
\text{$E^{\text{gap}}_{BKS}$ (\%)} & \text{$T_{\text{total}}$ (s)} &
\text{$E^{\text{gap}}_{BKS}$ (\%)} & \text{$T_{\text{total}}$ (s)} &
\text{$E^{\text{gap}}_{BKS}$ (\%)} & \text{$T^{\text{LA}}$ (s)} & \text{$T_{\text{total}}$ (s)} \\
\midrule
\endfirsthead
\toprule
\multicolumn{2}{c}{} &
\multicolumn{2}{c}{\text{TSBA$_{\text{basic}}$}} &
\multicolumn{2}{c}{\text{HGAMP}} &
\multicolumn{3}{c}{\text{NEO-DS}} \\
\cmidrule(lr){3-4} \cmidrule(lr){5-6} \cmidrule(lr){7-9}
Instance & BKS &
\text{$E^{\text{gap}}_{BKS}$ (\%)} & \text{$T_{\text{total}}$ (s)} &
\text{$E^{\text{gap}}_{BKS}$ (\%)} & \text{$T_{\text{total}}$ (s)} &
\text{$E^{\text{gap}}_{BKS}$ (\%)} & \text{$T^{\text{LA}}$ (s)} & \text{$T_{\text{total}}$ (s)} \\
\midrule
\endhead
\midrule
\multicolumn{9}{r}{\textit{Continued on next page}} \\
\endfoot
\endlastfoot
500-25-1a & 1771040 & 0.22 & 10864 & 0.10 & 10235 & 0.70 & 875.16 & 899.34 \\
500-25-1b & 1327800 & 0.50 & 5803 & 0.28 & 10432 & 0.94 & 514.26 & 539.51 \\
500-25-1c & 670025 & 0.86 & 6001 & 0.35 & 16214 & 4.36 & 32.90 & 36.15 \\
500-25-1d & 3321740 & 0.22 & 6627 & 0.09 & 10120 & 0.52 & 318.73 & 323.70 \\
500-25-1e & 2682740 & 10.91 & 7694 & 0.29 & 9314 & 0.90 & 18.32 & 38.78 \\
500-25-2a & 1619689 & 0.17 & 11998 & 0.21 & 12081 & 1.00 & 192.56 & 217.21 \\
500-25-2b & 1250270 & 0.31 & 5213 & 0.09 & 9587 & 0.58 & 224.27 & 248.45 \\
500-25-2c & 571155 & 0.97 & 5694 & 0.22 & 14168 & 4.31 & 28.89 & 34.03 \\
500-25-2d & 3335790 & 0.13 & 8014 & 0.09 & 10595 & 0.67 & 572.80 & 578.17 \\
500-25-2e & 2725350 & 5.02 & 11992 & 0.13 & 10784 & 0.70 & 71.45 & 91.25 \\
500-25-3a & 1724540 & 0.25 & 13855 & 0.23 & 10792 & 0.79 & 173.57 & 199.45 \\
500-25-3b & 1302820 & 0.32 & 5421 & 0.20 & 11222 & 1.02 & 994.62 & 1018.30 \\
500-25-3c & 578519 & 0.81 & 5099 & 0.26 & 14532 & 3.96 & 26.01 & 29.74 \\
500-25-3d & 3248557 & 0.05 & 6323 & 0.21 & 8733 & 0.97 & 1512.42 & 1517.61 \\
500-25-3e & 2681360 & 4.28 & 11756 & 0.17 & 10272 & 1.35 & 41.26 & 62.76 \\
500-25-4a & 1653680 & 0.30 & 13964 & 0.21 & 12931 & 1.49 & 207.02 & 235.25 \\
500-25-4b & 1259400 & 0.62 & 6661 & 0.22 & 11870 & 1.98 & 416.23 & 441.28 \\
500-25-4c & 664089 & 0.52 & 7012 & 0.45 & 14370 & 4.24 & 34.17 & 43.94 \\
500-25-4d & 3362588 & 0.06 & 8394 & 0.31 & 12719 & 1.24 & 972.71 & 978.20 \\
500-25-4e & 2625360 & 2.80 & 11932 & 0.17 & 11268 & 0.69 & 135.73 & 157.78 \\
500-30-1a & 1984150 & 0.36 & 11301 & 0.87 & 11688 & 0.67 & 71.89 & 90.17 \\
500-30-1b & 1532670 & 0.54 & 6369 & 0.40 & 11862 & 0.82 & 306.68 & 324.53 \\
500-30-1c & 611584 & 1.07 & 6389 & 0.22 & 15736 & 5.80 & 30.83 & 33.02 \\
500-30-1d & 3740220 & 0.15 & 6002 & 0.17 & 12247 & 1.42 & 259.39 & 263.01 \\
500-30-1e & 3236740 & 7.79 & 7836 & 0.19 & 11263 & 1.22 & 216.69 & 229.52 \\
500-30-2a & 1820346 & 0.13 & 15908 & 0.21 & 12927 & 1.12 & 437.05 & 455.37 \\
500-30-2b & 1452028 & 0.68 & 6908 & 0.24 & 12417 & 1.20 & 59.91 & 77.52 \\
500-30-2c & 647932 & 0.65 & 5547 & 0.27 & 15833 & 5.78 & 31.04 & 37.93 \\
500-30-2d & 3815210 & 0.02 & 8543 & 0.10 & 11682 & 1.73 & 21.35 & 25.17 \\
500-30-2e & 3218410 & 3.91 & 9895 & 0.18 & 11132 & 0.45 & 77.41 & 90.07 \\
500-30-3a & 1781220 & 0.20 & 9979 & 0.15 & 11986 & 1.36 & 68.65 & 87.82 \\
500-30-3b & 1417870 & 0.44 & 5206 & 0.20 & 10957 & 1.48 & 70.12 & 86.90 \\
500-30-3c & 567988 & 0.91 & 5936 & 0.29 & 12775 & 4.70 & 35.05 & 38.77 \\
500-30-3d & 3690995 & 0.12 & 8058 & 0.41 & 11606 & 1.01 & 188.18 & 191.83 \\
500-30-3e & 3059470 & 7.39 & 10411 & 2.52 & 10203 & 2.71 & 156.21 & 172.39 \\
500-30-4a & 1716476 & 0.36 & 12071 & 0.37 & 13181 & 1.81 & 395.23 & 413.59 \\
500-30-4b & 1398401 & 4.50 & 6314 & 0.62 & 13323 & 1.70 & 245.50 & 263.59 \\
500-30-4c & 558398 & 1.19 & 9257 & 0.32 & 13568 & 5.94 & 34.44 & 42.96 \\
500-30-4d & 3708479 & 0.08 & 8341 & 0.40 & 14431 & 0.83 & 43.12 & 46.98 \\
500-30-4e & 3057950 & 5.21 & 9918 & 0.92 & 11804 & 0.62 & 26.20 & 41.76 \\
600-30-1a & 2197860 & 0.13 & 20141 & 0.24 & 15026 & 0.70 & 1102.76 & 1129.95 \\
600-30-1b & 1689980 & 0.46 & 10673 & 0.35 & 15707 & 0.71 & 1447.16 & 1474.77 \\
600-30-1c & 742450 & 1.41 & 10822 & 0.44 & 19931 & 5.25 & 47.70 & 51.31 \\
600-30-1d & 4213337 & 0.09 & 14804 & 0.38 & 16852 & 0.58 & 281.83 & 288.23 \\
600-30-1e & 3552440 & 8.19 & 15500 & 0.16 & 14743 & 0.70 & 209.55 & 233.40 \\
600-30-2a & 2015580 & 0.18 & 23639 & 0.15 & 14461 & 0.99 & 277.82 & 307.02 \\
600-30-2b & 1601140 & 0.23 & 10779 & 0.10 & 14258 & 1.20 & 1489.49 & 1515.67 \\
600-30-2c & 633935 & 0.70 & 10719 & 0.21 & 16477 & 2.75 & 45.62 & 51.69 \\
600-30-2d & 4160240 & 0.16 & 15024 & 0.11 & 15083 & 0.73 & 272.19 & 279.11 \\
600-30-2e & 3563860 & 4.22 & 15010 & 0.61 & 14761 & 0.78 & 37.27 & 57.51 \\
600-30-3a & 2082824 & 0.18 & 24898 & 0.48 & 15086 & 0.89 & 881.73 & 909.03 \\
600-30-3b & 1612210 & 0.54 & 10761 & 0.26 & 15231 & 1.03 & 475.45 & 505.33 \\
600-30-3c & 658376 & 0.98 & 9549 & 0.36 & 18726 & 2.98 & 48.56 & 53.80 \\
600-30-3d & 4068474 & 0.07 & 15657 & 0.46 & 14432 & 0.99 & 61.96 & 67.90 \\
600-30-3e & 3327790 & 5.33 & 15180 & 0.23 & 12476 & 1.41 & 31.31 & 56.49 \\
600-30-4a & 1938780 & 0.27 & 25361 & 0.21 & 15596 & 1.36 & 1634.74 & 1660.38 \\
600-30-4b & 1554350 & 0.46 & 11437 & 0.34 & 15512 & 1.54 & 954.02 & 983.29 \\
600-30-4c & 705601 & 0.41 & 14535 & 0.27 & 17902 & 4.13 & 45.32 & 58.49 \\
600-30-4d & 4150775 & 0.10 & 17603 & 0.98 & 18025 & 1.11 & 67.78 & 74.82 \\
600-30-4e & 3514040 & 3.28 & 18513 & 0.18 & 14803 & 0.23 & 48.99 & 70.13 \\
\midrule
\bfseries Average & \multicolumn{1}{c}{\bfseries } & \bfseries 1.54 & \bfseries 10852 & \bfseries 0.33 & \bfseries 13299 & \bfseries 1.75 & \bfseries 326.65 & \bfseries 342.20 \\
\midrule
\multicolumn{2}{l}{Processor}
  & \multicolumn{2}{r}{Xeon E5-2670}
  & \multicolumn{2}{r}{Xeon E-2670}
  & \multicolumn{3}{r}{Xeon Gold 6248R} \\
\multicolumn{2}{l}{GHz}
  & \multicolumn{2}{r}{2.6}
  & \multicolumn{2}{r}{2.5}
  & \multicolumn{3}{r}{3.0}\\
\bottomrule
\end{longtable}
\endgroup

\subsection{Routing Solver Ablation}
Table~\ref{tab:detailed_routing_solver} reports detailed results for Section~\ref{sec:routing_solver}.

\subsection{Normalization Scheme Ablation}
Table~\ref{tab:detailed_normalization} reports detailed results for Section~\ref{sec:target_normalization}.

\subsection{Neural Architecture Ablation}
Table~\ref{tab:optimization_neods_neogt} reports detailed results for Section~\ref{sec:ablation_architecture}.

\begin{table}[H]
\centering
\scriptsize
\caption{Comparison of NEO-DS and NEO-GT on the $\mathbb{P}$ benchmark of \citet{prins2004nouveaux}. BKS values are taken from \citet{loffler2023conceptually} and \citet{he2025hybrid}.}
\label{tab:optimization_neods_neogt}
\renewcommand{\arraystretch}{0.9}%
\resizebox{\ifdim\width>\linewidth \linewidth\else \width\fi}{!}{%
\begin{tabular}{lrrrrrrr}
\toprule
 & & \multicolumn{3}{c}{NEO-DS}
   & \multicolumn{3}{c}{NEO-GT} \\
\cmidrule(r){3-5}\cmidrule(l){6-8}
Instance & BKS
  & $E^{\text{gap}}_{\text{BKS}}$ (\%) & $T^{\text{LA}}$ (s) & $T_{\text{total}}$ (s)
  & $E^{\text{gap}}_{\text{BKS}}$ (\%) & $T^{\text{LA}}$ (s) & $T_{\text{total}}$ (s) \\
\midrule
20-5-1a & 54793 &  3.24  &  0.19  &  0.20  &  7.79  &  0.38  &  0.40  \\
20-5-1b & 39104 &  3.65  &  0.17  &  0.19  &  6.58  &  0.24  &  0.26  \\
20-5-2a & 48908 &  5.82  &  0.15  &  0.16  &  13.26  &  0.40  &  0.41  \\
20-5-2b & 37542 &  5.49  &  0.16  &  0.19  &  20.48  &  0.37  &  0.40  \\
\midrule
\textbf{Average} &  & \textbf{ 4.55 } & \textbf{ 0.17 } & \textbf{ 0.18 } & \textbf{ 12.03 } & \textbf{ 0.35 } & \textbf{ 0.37 } \\
\midrule
50-5-1a & 90111 &  1.63  &  0.26  &  0.37  &  13.36  &  0.46  &  0.72  \\
50-5-1b & 63242 &  1.18  &  0.28  &  0.45  &  6.79  &  0.60  &  0.78  \\
50-5-2a & 88293 &  2.74  &  0.32  &  0.45  &  15.17  &  5.79  &  5.93  \\
50-5-2b & 67308 &  5.80  &  0.28  &  0.38  &  24.87  &  0.99  &  1.09  \\
50-5-2bBIS & 51822 &  4.86  &  0.46  &  0.53  &  19.30  &  6.72  &  6.82  \\
50-5-2BIS & 84055 &  1.62  &  0.25  &  0.45  &  34.13  &  0.59  &  0.77  \\
50-5-3a & 86203 &  0.96  &  0.37  &  0.59  &  8.15  &  1.25  &  1.54  \\
50-5-3b & 61830 &  4.32  &  0.40  &  0.60  &  2.23  &  0.77  &  0.93  \\
\midrule
\textbf{Average} &  & \textbf{ 2.89 } & \textbf{ 0.33 } & \textbf{ 0.48 } & \textbf{ 15.50 } & \textbf{ 2.15 } & \textbf{ 2.32 } \\
\midrule
100-5-1a & 274814 &  1.49  &  1.40  &  2.39  &  6.66  &  1.99  &  3.10  \\
100-5-1b & 213568 &  1.25  &  1.36  &  2.25  &  11.92  &  6.35  &  7.29  \\
100-5-2a & 193671 &  0.28  &  0.45  &  2.22  &  9.36  &  1.07  &  2.45  \\
100-5-2b & 157095 &  0.79  &  0.51  &  2.18  &  8.52  &  6.07  &  7.49  \\
100-5-3a & 200079 &  1.44  &  0.43  &  1.94  &  2.78  &  0.93  &  2.35  \\
100-5-3b & 152441 &  3.09  &  0.45  &  1.81  &  6.21  &  1.68  &  2.82  \\
100-10-1a & 287661 &  1.52  &  3.78  &  4.60  &  5.33  &  36.18  &  37.00  \\
100-10-1b & 230989 &  2.35  &  1.80  &  2.38  &  5.25  &  44.94  &  45.56  \\
100-10-2a & 243590 &  1.49  &  2.95  &  3.75  &  8.57  &  4.36  &  5.23  \\
100-10-2b & 203988 &  2.54  &  2.83  &  3.56  &  4.72  &  36.15  &  36.67  \\
100-10-3a & 250882 &  1.79  &  3.25  &  4.18  &  14.56  &  3.07  &  3.92  \\
100-10-3b & 203114 &  3.08  &  25.65  &  26.52  &  11.72  &  5.01  &  5.50  \\
\midrule
\textbf{Average} &  & \textbf{ 1.76 } & \textbf{ 3.74 } & \textbf{ 4.81 } & \textbf{ 7.97 } & \textbf{ 12.32 } & \textbf{ 13.28 } \\
\midrule
200-10-1a & 474702 &  1.03  &  5.89  &  10.85  &  6.21  &  61.00  &  65.52  \\
200-10-1b & 375177 &  0.74  &  60.58  &  64.80  &  12.31  &  61.07  &  65.36  \\
200-10-2a & 448077 &  0.65  &  2.28  &  6.40  &  15.42  &  47.48  &  51.88  \\
200-10-2b & 373696 &  0.71  &  3.02  &  6.65  &  10.47  &  42.67  &  46.17  \\
200-10-3a & 469433 &  1.67  &  5.58  &  10.51  &  6.86  &  21.95  &  25.93  \\
200-10-3b & 362320 &  1.30  &  2.96  &  7.33  &  9.86  &  7.85  &  11.81  \\
\midrule
\textbf{Average} &  & \textbf{ 1.02 } & \textbf{ 13.38 } & \textbf{ 17.76 } & \textbf{ 10.19 } & \textbf{ 40.34 } & \textbf{ 44.44 } \\
\bottomrule
\end{tabular}%
}
\end{table}

\end{appendices}
\end{document}